\documentclass[10pt]{amsart}

\usepackage{graphicx} 
\usepackage[svgnames]{xcolor}
\usepackage[all]{xy}
\usepackage{tikz-cd}
\usepackage{stmaryrd}
\usepackage{multicol}
\usepackage{verbatim}

\usepackage{float}

\usepackage[table]{xcolor}

\usepackage{cancel}

\usepackage{amsmath, amssymb, amsfonts, mathtools}
\usepackage{amsthm}
\usepackage[svgnames]{xcolor}
\usepackage{tikz}
\usetikzlibrary{graphs}
\usepackage{ulem}
\usepackage{tikz,tikz-cd}
\usepackage{tensor}
\usepackage{mathrsfs}
\usepackage{enumitem}
\usepackage{dsfont}

\newcommand{\sbl}{\scalebox{0.5}{$\blacklozenge$}}

\newcommand{\xto}{\xrightarrow}

\newcommand{\circnum}[1]{%
  \tikz[baseline=(char.base)]{
    \node[
      draw,
      circle,
      inner sep=1.2pt,
      line width=0.4pt
    ] (char) {\textnormal{\footnotesize #1}}; 
  }%
}

\newlist{cenum}{enumerate}{1}
\setlist[cenum,1]{%
  label=\protect\circnum{\arabic*}, 
  labelsep=0.5em,
  leftmargin=2em,
  itemsep=0.2em,
  align=left
}

\newcommand{\cnum}[1]{\raisebox{0.3pt}{\textcircled{\scriptsize #1}}}

\usepackage[colorlinks=true, linkcolor=blue]{hyperref}

\usepackage{setspace}

\newcommand{\llangle}[1]{\langle\!\langle #1 \rangle\!\rangle}

\renewcommand{\emph}[1]{\textit{#1}}

\theoremstyle{plain}
\newtheorem{Thm}{Theorem}[section]
\newtheorem{Lem}[Thm]{Lemma}
\newtheorem{Cor}[Thm]{Corollary}
\newtheorem{Prop}[Thm]{Proposition}
\newtheorem{Prop-Def}[Thm]{Proposition-Definition}

\theoremstyle{definition}
\newtheorem{Def}[Thm]{Definition}
\newtheorem{Ex}[Thm]{Example}
\newtheorem{Not}[Thm]{Notation}
\newtheorem{Rmk}[Thm]{Remark}
\newtheorem{Fact}[Thm]{Fact}

\numberwithin{equation}{section}

\newenvironment{Proof}[1][Proof]{\begin{trivlist}
		\item[\hskip \labelsep {\bfseries #1}]}{\flushright
		$\Box$\end{trivlist}}

 \normalbaselines
\usepackage{tikz}
\usetikzlibrary{cd}
\usetikzlibrary{calc}
\usetikzlibrary{arrows}
\usetikzlibrary{arrows.meta}
\usetikzlibrary{decorations}
\usetikzlibrary{backgrounds}
\usetikzlibrary{decorations.pathmorphing}
\usetikzlibrary{decorations.markings}

\tikzset{
  snake it/.style={decorate, decoration=snake},
  bl/.style={circle, draw=white, thin,fill=black!100, scale=0.5},
  cg/.style={circle, draw, thin,fill=black!10, scale=0.8},
  cw/.style={circle, draw, thin,fill=white, scale=0.8},
   bw/.style={circle, draw=black, fill=white, scale=0.3},
  B/.style={circle, draw,fill=black!30, scale=1.6},
  b/.style={circle, draw,fill=black, scale=0.3},
  r/.style={circle, draw=green!10!red,fill=green!10!red, scale=0.3},
}

\newcommand{\ot}{\otimes}
\newcommand{\Char}{\mathrm{char}}

\newcommand{\Sp}{\mathrm{Sp}}

\newcommand{\rmHom}{\mathrm{Hom}}
\newcommand{\rmrot}{\mathrm{rot}}
\newcommand{\rmsrot}{\mathrm{srot}}
\newcommand{\HH}{\mathrm{HH}}
\newcommand{\rmH}{\mathrm{H}}
\newcommand{\basic}{\mathrm{basic}}
\newcommand{\Ker}{\mathrm{Ker}}
\newcommand{\rmIm}{\mathrm{Im}}
\newcommand{\rt}{\mathsf{rt}}

\newcommand{\calB}{{\mathcal B}}
\newcommand{\calH}{{\mathcal H}}

\newcommand{\bbZ}{{\mathbb Z}}
\newcommand{\bbK}{{\mathbb K}}
\newcommand{\bbN}{{\mathbb N}}
\newcommand{\bbP}{{\mathbb P}}
\newcommand{\bbB}{{\mathbb B}}
\newcommand{\oA}{{\overline{A}}}
\newcommand{\oC}{{\overline{\mathscr C}}}

\newcommand{\scrG}{{\mathscr G}}
\newcommand{\scrR}{{\mathscr R}}
\newcommand{\C}{{\mathscr C}}

\newcommand\lfact[2]{{#2}\if1#1'\fi_{\mathsf{lt}}}
\newcommand\rfact[2]{{#2}\if2#1'\fi_{\mathsf{rt}}}

\title[Hochschild cohomology of graded skew-gentle algebras]{Hochschild cohomology of graded skew-gentle algebras: Gerstenhaber algebra structure and geometric interpretation}

\author[Xiuli Bian, Sibylle Schroll, Andrea Solotar, Xiao-chuang Wang, Can Wen]%
       {Xiuli Bian, Sibylle Schroll, Andrea Solotar, Xiao-chuang Wang, Can Wen}
\date{}



\address{School of Mathematical Sciences, East China Normal University, Shanghai 200241, China
\newline Department of Mathematics, Universit\"at zu K\"oln, Weyertal 86-90, 50931 K\"oln, Germany }
\email{52205500018@stu.ecnu.edu.cn}

\address{Department of Mathematics, Universit\"at zu K\"oln, Weyertal 86-90, 50931 K\"oln, Germany}
\email{schroll@math.uni-koeln.de}

\address{  IMAS-CONICET and Departamento de Matemática, Facultad de Ciencias Exactas y Naturales, Universidad de Buenos Aires, Pabellon I, Ciudad Universitaria, Buenos Aires, 1428, Argentina.
\newline Guangdong Technion Israel Institute of Technology, Shantou, Guangdong Province, China.} 

\email{asolotar@dm.uba.ar }

\address{ School of Mathematical Sciences, University of Science and Technology of China, Hefei, Anhui, 230026, CHINA.
\newline Department of Mathematics, Universit\"at zu K\"oln, Weyertal 86-90, 50931 K\"oln, Germany.}
\email{wxchuang@mail.ustc.edu.cn}

\address{School of Mathematical Sciences, Laboratory of Mathematics and Complex Systems, Beijing Normal University, Beijing 100875, China.
\newline Department of Mathematics, Universit\"at zu K\"oln, Weyertal 86-90, 50931 K\"oln, Germany.
}
\email{cwen@mail.bnu.edu.cn}

\begin{document}

\begin{abstract}
In this paper we calculate the Hochschild cohomology of graded skew-gentle algebras, together with its structure as graded commutative algebra under the cup product and its Lie algebra structure given by the Gerstenhaber bracket.
One of the results of this paper is that for graded skew-gentle algebras and thus also for partially wrapped Fukaya categories of orbifold surfaces with stops, their Hochschild cohomology is encoded in the underlying graded surface.

\end{abstract}

\maketitle

\tableofcontents

\section{Introduction}

Although explicit calculations of Hochschild cohomology and its Gerstenhaber algebra structure are difficult, there is an increasing interest in such calculations. In particular, suitable projective bimodule resolutions such as the Bardzell resolution and the more general resolution in \cite{CS2015} make this possible, and much work has been done in recent years, see, for example, \cite{Sol}, \cite{VW} and the references therein. However, calculations in the case of graded algebras are more rare. In this paper, we do such calculations for the class of graded skew-gentle algebras. More precisely, we calculate the Hochschild cohomology, its structure as graded commutative algebra under the cup product and its Lie algebra structure given by the Gerstenhaber bracket. 

There is a second motivation for calculating the Hochschild cohomology of graded skew-gentle algebras. Namely, in \cite{HKK17}, the (topological) partially wrapped Fukaya category $\mathcal{W}(S)$ of a graded marked surface $S$ is explicitly determined and it is shown that its derived category is equivalent to the perfect derived category of a graded gentle algebra. In \cite{AP24, CK24, BSW24} this is extended to the partially wrapped Fukaya category $\mathcal{W}(S)$ of a graded marked surface $S$ with orbifold points of order 2. It is shown that the derived category of such a Fukaya category is  equivalent to the perfect derived category of a graded skew-gentle algebra. Moreover, in \cite{BSW25}  it is shown that any $A_\infty$-deformation of the partially wrapped Fukaya category of a smooth surface with stops is in fact of this form, that is, it deforms to a partially wrapped Fukaya category of a surface with orbifold points of order 2. We also note that in order to show this, in \cite{BSW25} the second Hochschild cohomology of a graded gentle algebra was calculated and related to certain boundary components and their winding numbers. 

Coming back to Hochschild cohomology: 
in the case of  wrapped Fukaya categories it is shown in \cite{Ganatra1, Ganatra2} that there is an isomorphism between the symplectic cohomology of the surface and the Hochschild cohomology of the wrapped Fukaya category. A similar result is expected to hold in the partially wrapped case, that is in the presence of stop data corresponding in this paper to red marked points.

For this reason, a second motivation  for us was to determine whether and how the Hochschild cohomology of graded skew-gentle algebras can be `read off' the associated surface models. 
An indication that this should be possible comes from the case of ungraded gentle algebras, where in \cite{CSSS24} it is shown that the Hochschild cohomology can indeed be read off the surface. One of the results of this paper is that the same is true for graded skew-gentle algebras and thus also for partially wrapped Fukaya categories of orbifold surfaces, namely, their Hochschild cohomology is encoded in the underlying graded surface.

More precisely, after determining a basis of the Hochschild cohomology of a graded skew-gentle algebra (Theorem~\ref{Thm: coho-basis-sg}) and its graded commutative algebra structure under the cup product $\smile$ (Table~\ref{table: cup product}), we explicitly present this structure by  generators (Theorem~\ref{Thm: basis-as-alg}) and relations (Theorem~\ref{Thm: relation of HH}). 
We also recover a result from \cite{BSW25}, in that we explicitly show how the Hochschild cohomology of a graded skew-gentle algebra differs from the Hochschild cohomology of a graded gentle algebra of which it is an (associative) deformation (Theorem~\ref{gentle-vs-sg}).
In Section~\ref{Sec: geometric}, we show that the generators of the Hochschild cohomology of a graded skew-gentle algebra (and indeed all basis elements), can be `read off' the graded marked surface with orbifold points of order 2 associated to the graded skew-gentle algebra. We summarise this latter result in the following theorem, see also Table~\ref{table: Hochschild} for a pictorial description of this correspondence.

\begin{Thm}[Theorem~\ref{geom.cohom}] Let $A$ be a graded skew-gentle algebra and  let $S$ be the graded marked surface with orbifold points of order 2 associated to $A$. Denote by $\eta$ the grading structure of $S$ which is given in form of a line field $\eta$ on $S$.  Then the generators of the Hochschild cohomology algebra $(\HH^*(A), \smile)$ of $A$ are in bijection with the set consisting of 
\begin{enumerate}
\item[$(1)$] the boundary components with exactly one marked point, 
\item[$(2)$] the red and green marked points in the interior of the surface (corresponding to unstopped and fully stopped boundary components in \cite{BSW24}), 
\item[$(3)$] the generators of the fundamental group of the smooth surface  which is obtained from $S$ by replacing a small disk around each orbifold point by a smooth disk.
\end{enumerate}
In particular, we show that the winding numbers with respect to $\eta$ of the closed curves around the boundary components and the interior points in the surface determine the total degrees of the  corresponding generators of $\HH^*(A)$. 
\end{Thm}

 Finally, we calculate the Gerstenhaber bracket of the basis elements of the Hochschild cohomology of a graded skew-gentle algebra (Theorem~\ref{thm: bracket of HH}) and give a possible geometric interpretation. 

For ease of navigation, we succinctly summarise the contents of this paper.
In Section \ref{Sec:Preliminaries} 
we recall some notations and basic facts. In Section \ref{Sec3} we compute the Hochschild cohomology of graded skew-gentle algebras, obtaining an explicit description of how the Hochschild cohomology of these deformations of graded gentle algebras differs from the Hochschild cohomology of the latter, see Theorem \ref{gentle-vs-sg}. We also  describe completely the  Gerstenhaber algebra structure in Theorem \ref{Thm: relation of HH} (graded commutative algebra structure) and 
Theorem \ref{thm: bracket of HH} 
(Gerstenhaber bracket). 
Section \ref{Sec: geometric} contains the geometric interpretation of this structure which is summarised in Theorem \ref{geom.cohom}, see also Table~\ref{table: Hochschild}. 

While working on this paper, we learned that Sebastian Opper has also been working on the Gerstenhaber algebra structure of the Hochschild cohomology of graded gentle algebras.

\textbf{Conventions:} throughout this paper $\bbK$ is a field with characteristic denoted by $\Char(\bbK)$. For a bound quiver algebra  $\bbK Q/I$, let $E=\bbK Q_0$ denote the semisimple  $\bbK$-algebra with basis the set of vertices of $Q$. In the differential graded (dg for short) situation, we always use the \emph{Koszul sign rule}, i.e., we add the sign $(-1)^{|x|\cdot|y|}$ whenever we exchange the positions of two graded objects $x$ and $y$ of degrees $|x|$ and $|y|$, respectively. Paths will be written from right to left, as compositions of functions.

\section{Preliminaries}\label{Sec:Preliminaries}

In this section we recall some basic notions and facts that are necessary throughout the paper.

\subsection{Graded skew-gentle algebras}\label{section:sg-algebra}\

Skew-gentle algebras were introduced in \cite{GP99}. 
These algebras can be presented by generators and relations as quotients of path algebras by ideals generated by quadratic mononials and quadratic commutation relations. They are thus graded by the length of the paths. It is worth noting that they are deformations of gentle algebras by $2$-cocycles in the Hochschild cohomology \cite{BSW25}.

Now we fix some notations while recalling some basic facts. A \emph{quiver} $Q$ is  a tuple $(Q_0,Q_1,s,t)$, where $Q_0$ is the set of \emph{vertices}, $Q_1$ is the set of \emph{arrows}, and $s,t: Q_1\to Q_0$ are the \emph{source} and \emph{target} functions of the arrows, respectively. We will only consider quivers with a finite number of vertices and a finite number of arrows. We call an arrow $\alpha\in Q_1$ a \emph{loop} if $s(\alpha)=t(\alpha)$. A \emph{path of length $l$} in $Q$ is a sequence $a_l\cdots a_2 a_1$ of arrows $a_i$, $1\le i\le l$ with $t(a_i)=s(a_{i+1})$ for any $1\le i\le l-1$. We denote by $l(p)$ the \emph{length} of the path $p$. For each vertex $i\in Q_0$, we denote by $e_i$ the trivial path of length $0$ associated to it.  A path $C=c_l\cdots c_1$ in $Q$ is a \emph{cycle} if its length $l$ is positive and it ends at the vertex where it begins. A \emph{relation set} of $Q$ is a finite set of linear combinations of paths in $Q$. We denote by $\bbK Q$ the path algebra of $Q$.

A \emph{graded quiver} is a quiver $Q$ with a grading, i.e., a function $|\cdot|:\ Q_1\to \bbZ$. For any arrow $\alpha\in Q_1$, denote by $|\alpha|$ the \emph{degree} of $\alpha$. For a path $p=\alpha_n\cdots\alpha_1$ of length $n\ge 1$, the degree of $p$ is defined as $|p|:=\sum_{i=1}^n |\alpha_i|$, while for a trivial path $e_i$, $i\in Q_0$, the degree of $e_i$ is defined as $|e_i|:=0$. The path algebra of a graded quiver is naturally a $\bbZ$-graded algebra. A quiver can always be graded by the length of the paths.\\

\begin{Def}\label{graded-gentle-pair}
A pair $(Q,R)$ of a quiver $Q$ and a relation set $R$ is called a \emph{gentle pair} if the following hold.
\begin{itemize}
        \item [(G1)] Any element in $R$ is a path of length $2$,

        \item [(G2)] for each vertex $i$ in $Q_0$, there are at most two arrows with source $i$ 
        and at most two arrows with target $i$,

        \item [(G3)] for any arrow $\alpha$, there is at most one arrow $\beta$ (resp. $\gamma$) with $\alpha\beta\in R$ (resp. $\gamma\alpha\in R$),

        \item [(G4)] for any arrow $\alpha$, there is at most one arrow $\beta$ (resp. $\gamma$) with $\alpha\beta\not\in R$ (resp. $\gamma\alpha\not\in R$).
    \end{itemize}
    
A gentle pair $(Q,R)$ is called a \emph{graded gentle pair} if, moreover $Q$ is a graded quiver. A graded $\bbK$-algebra $A$ is called \emph{graded gentle} if it is Morita equivalent to $\bbK Q/\langle R \rangle$ with $(Q,R)$ a graded gentle pair.
\end{Def}

\begin{Def}
A triple $(Q,R,\Sp)$ is called \emph{graded skew-gentle triple} if it satisfies
\begin{itemize}
		\item $Q$ is a graded quiver, 
            
		\item $\Sp$ is a subset of loops in $Q$ (called \emph{special loops}), and
            
		\item $( Q, R\cup \{\varepsilon^2\mid \varepsilon \in \Sp\} ) $ is a graded gentle pair.
\end{itemize}

A graded $\bbK$-algebra $A$ is called \emph{graded skew-gentle} if it is Morita equivalent to 
\[
\bbK Q/\langle R\cup \{\varepsilon^2-\varepsilon\ |\ \varepsilon\in \Sp\} \rangle
.\]
When $\Sp=\emptyset$, the algebra $A$ is graded gentle. A vertex $i\in Q_0$ is called \emph{special} if there is a loop $\varepsilon_i\in \Sp$ such that $i=s(\varepsilon_i)$. By abuse of notation, we still write the set of special vertices as $\Sp$.
\end{Def}

\medskip 

\begin{Rmk}
  Note that a skew-gentle algebra $A$ can be graded only if  the degree of the special loops is $0$. 
\end{Rmk}

We now fix a graded skew-gentle triple $(Q,R,\Sp)$, and consider the graded skew-gentle algebra 
\[
A=\bbK Q/\langle R\cup \{\varepsilon^2-\varepsilon\ |\ \varepsilon\in \Sp\} \rangle.
\]
Let $\pi:\bbK Q \to A$ be the canonical projection. Denote by $\calB$ the set of paths in $Q$ that do not have subpaths in $S= R\cup \{\varepsilon^2\mid \varepsilon \in \Sp\}$. The set $\calB$ is clearly a $\bbK$-basis of $A$. 

Hereafter, the graded skew-gentle algebra $A$ we refer to is always defined by a graded skew-gentle triple $(Q,R,\Sp)$; that is, $A=\bbK Q/I$, where $I$ denotes the two-sided ideal of $\bbK Q$ generated by $R\cup \{ \varepsilon^2-\varepsilon\ |\ \varepsilon\in\Sp \}$.

\subsection{Gerstenhaber algebra structure in Hochschild cohomology of graded algebras} 
\label{section: TT calculus of dg algebra}\

In this section, 
we recall the Gerstenhaber structure of the Hochschild cohomology of graded $\bbK$-algebras. Even if the action of $E=\bbK Q_0$ on the algebras we consider is not symmetric, we will still call them $E$-algebras like in \cite{HLW23}.
The primary reference for this section is \cite[Chapter 3]{HLW23}. However, we note that our sign occasionally differs from that in \cite{HLW23}.

Hochschild (co)homology of an associative algebra $A$ with coefficients in an $A$-bimodule $M$ was originally introduced in \cite{H45} and its Gerstenhaber structure was defined in \cite{G63}. This definition has been extended to differential graded (dg) $\bbK$-algebras \cite{FMT02, Abb15} and dg $E$-algebras \cite{HLW23}.

We always view the graded algebra $A$ as a dg algebra with trivial differential of degree 1 (i.e., $d_A=0$ and $|d_A|=+1$). Denote by $s$  the suspension of $A$ with $|s|=-1$ and $s^{-1}\circ s=\mathrm{id}$. Let $A=E\oplus \oA$ and $\ot:=\ot_E$. \\

\textbf{Two-sided bar resolution.} In \cite{HLW23}, the Hochschild (co)homology of dg $E$-rings
was originally defined using the two-sided bar resolution. 
We recall this definition for the case of graded algebras here.

The \emph{two-sided bar resolution} \cite{CE56, ML63, HLW23} is given by $B(A,A,A): = A\otimes_E T^c(s \oA)\otimes_E A$ equipped with differential $\delta=\delta_0+\delta_1$, where $T^c(s \oA)$ is the tensor dg coalgebra of $s\oA$ over $E$, 
that is,
\[
T^c(s \oA)=\bigoplus_{n=0}^{\infty}(s\oA)^{\otimes n}=E \oplus s\oA\oplus (s\oA)^{\otimes 2}\oplus\cdots,
\]

and $\delta_0$ is the \emph{internal differential} given by 
\begin{align*}
    \delta_0(a\ot sa_n\ot \cdots sa_1\ot b)=
    & d_A(a) \ot sa_n\ot \cdots \ot sa_1\ot b\\
    &+\sum_{i=1}^{n} (-1)^{|a|+\epsilon_{i+1,n}-1} a\ot sa_n\ot\cdots \ot sa_{i+1}\ot sd_A(a_i)\ot sa_{i-1}\ot\cdots\ot sa_1\ot b\\
    &+(-1)^{|a|+\epsilon_{1,n}} a\ot sa_n\ot \cdots sa_1\otimes d_A(b)
\end{align*}
where $\epsilon_{i,n}=\sum\limits_{j=i}^n(|a_j|-1)$ for $1\le i\le n$, with the convention $\epsilon_{n+1,n}=0$. The map $\delta_1$ denotes the \emph{external differential}, defined by 
\begin{align*}
\delta_1(1\ot sa_n\ot \cdots sa_1\ot 1)=
& a_n\otimes sa_{n-1} \ot \cdots\ot sa_1\otimes 1\\
&+\sum_{i=2}^n (-1)^{\epsilon_{i,n}} 1\ot sa_n\ot \cdots \ot sa_{i+1}\ot s(a_ia_{i-1})\ot sa_{i-2}\ot \cdots\ot sa_1\otimes 1\\
&+(-1)^{\epsilon_{2,n}-1} 1\otimes sa_n\otimes\cdots\otimes sa_2\otimes a_1.
\end{align*}
Here, $a, b$ and $a_i(1\le i\le n)$ are homogeneous elements of $A$ and $\delta_1$ is an $A$-bimodule morphism. Note that for a graded algebra $A$, the internal differential $\delta_0$ induced from $d_A=0$ is zero, hence the differential $\delta=\delta_1$ is given by the external differential.

\textbf{Hochschild cochain complex.} The \emph{Hochschild cochain complex} of $A$ is
\[
C^*(A):=\rmHom_{A^e}(B(A,A,A),A)\cong\rmHom_{E^e}(T^c(s\oA), A),
\]
where $A^{op}$ is the algebra whose underlying graded vector space is $A$ with the opposite multiplication of $A$, i.e. $a\cdot^{op}b=(-1)^{|a|\cdot |b|}b\cdot a$ and $A^e:=A\otimes A^{op}$ (resp. $E^e$) denotes the enveloping algebra of $A$ (resp. $E$).
Its cohomology is called the \emph{Hochschild cohomology} of $A$, denoted by $\mathrm{HH}^*(A)$.  More precisely, given a homogeneous cochain $f\in\rmHom_{E^e}((s\oA)^{\otimes n},A)$, we say that $f$ is of degree $p$, denoted by $|f|=p$, if $f\Big(\big((s\oA)^{\otimes n}\big)^i\Big)\subseteq \big((s\oA)^{\otimes n}\big)^{i+p}$ for all $i$. From now on, 
$\rmHom_{E^e}((s\oA)^{\otimes n},A)^p$ will denote the subspace of 
$\rmHom_{E^e}((s\oA)^{\otimes n},A)$ of maps of degree $p$, and 
 \[
 (C^*(A))^p=\bigoplus_{n\in \bbN}\rmHom_{E^e}((s\oA)^{\otimes n},A)^p
\]
with differential given by $\partial=\partial_0+\partial_1$, and
for each $f\in\rmHom_{E^e}((s\oA)^{\otimes n},A)^p$, the \emph{internal differential} $\partial_0$ maps $f$ to
\begin{align*}
\partial_0(f)(sa_n\ot\cdots\ot sa_1):=
& d_Af(sa_n\ot\cdots\ot sa_1)\\
&+\sum_{i=1}^n (-1)^{|f|+\epsilon_{i+1,n}} f(sa_n\ot\cdots \ot sa_{i+1}\ot sd_A(a_i)\ot sa_{i-1}\ot\cdots\ot sa_1),
\end{align*}
which is zero when $d_A=0$, and the \emph{external differential} $\partial_1$ maps $f$ to
\begin{align*}
\partial_1(f)(sa_{n+1}\ot \cdots\ot sa_1)
:=&-(-1)^{(|a_{n+1}|-1)|f|}a_{n+1} f(sa_n\ot \cdots \ot sa_1)\\
&+\sum_{i=2}^{n+1} (-1)^{|f|+\epsilon_{i,n+1}-1}f(sa_{n+1}\ot \cdots\ot sa_{i+1}\ot s(a_i a_{i-1})\ot sa_{i-2}\ot \cdots\ot sa_1)\\
    &+(-1)^{|f|+\epsilon_{2,n+1}}f(sa_{n+1}\ot \cdots \ot sa_2)a_1.
\end{align*}\\

\textbf{Cup product and Gerstenhaber bracket.} For $f\in \rmHom_{E^e}((s\oA)^{\otimes m},A)$ and $g\in \rmHom_{E^e}((s\oA)^{\otimes n},A)$, the \emph{cup product} $f\smile g\in \rmHom_{E^e}((s\oA)^{\otimes m+n},A)$ is defined by
\[
(f\smile g) (s a_{m+n}\ot\cdots\ot sa_1):=(-1)^{|g|\epsilon_{n+1,m+n}} f(sa_{m+n}\ot \cdots sa_{n+1}) g(sa_n\ot \cdots sa_1),
\]
where $\epsilon_{n+1,m+n}:=\sum\limits_{j=n+1}^{m+n}(|a_j|-1)$ as before.

The \emph{Gerstenhaber bracket} $[f,g]\in  \rmHom_{E^e}((s\oA)^{\otimes m+n-1},A)$ is defined by
\[
[f,g]:= f\circ g-(-1)^{(|f|-1)(|g|-1)} g\circ f,  
\]
where $f\circ g=\sum\limits_{i=1}^m f\circ_i g$ with 
\[
(f\circ _i g)(sa_{m+n-1}\ot\cdots\ot s a_1):=(-1)^{\epsilon_{i+n,m+n-1}(|g|+1)} f(sa_{m+n-1}\ot\cdots\ot s g(s a_{i+n-1}\ot\cdots\ot sa_i)\ot\cdots\ot sa_1).
\]
In particular, if $n=0$, then
\[
(f\circ_i g)(sa_{m-1}\ot \cdots \ot sa_1):=(-1)^{\epsilon_{i,m-1}(|g|+1)} f(sa_{m-1}\ot \cdots sa_i\ot sg\ot sa_{i-1}\ot\cdots \ot sa_1).
\]

The cup product satisfies the following properties.
\begin{Lem}\label{Lemma: properties of cup product}
    For all $f\in \rmHom_{E^e}((s\oA)^{\otimes m},A)$, $g\in \rmHom_{E^e}((s\oA)^{\otimes n},A)$ and $h\in \rmHom_{E^e}((s\oA)^{\otimes r},A)$: 
    \begin{itemize}
        \item [$(1)$] $\partial_1(f\smile g)=\partial_1(f)\smile g+(-1)^{|f|} f\smile \partial_1(g)$,
        \item [$(2)$] $(f\smile g)\smile h=f\smile (g\smile h)$,
        \item [$(3)$] $f\smile g-(-1)^{|f|\cdot|g|} g\smile f=(-1)^{|f|}\partial_1(f)\circ g+(-1)^{|f|+1} \partial_1(f\circ g)-f\circ \partial_1(g)$.
    \end{itemize}
\end{Lem}
\medskip

According to Lemma~\ref{Lemma: properties of cup product} (1), the cup product on $C^*(A)$ induces the cup product on $\HH^*(A)$. It follows from Lemma~\ref{Lemma: properties of cup product} (2) and (3) that $(\HH^*(A),\smile)$ is a graded commutative associative algebra.

The Gerstenhaber bracket $[-,-]$ satisfies the following properties.
\begin{Lem}\label{lemma: properties of Ger bracket}
    For all $f\in \rmHom_{E^e}((s\oA)^{\otimes m},A)$, $g\in \rmHom_{E^e}((s\oA)^{\otimes n},A)$ and $h\in \rmHom_{E^e}((s\oA)^{\otimes r},A)$, one has
    \begin{itemize}
        \item [$(1)$] $\partial_1 [f,g]=[\partial_1(f),g]+(-1)^{|f|+1} [f,\partial_1 (g)]$,

        \item [$(2)$] $[f,g]=-(-1)^{(|f|+1)(|g|+1)}[g,f]$,

        \item [$(3)$] \emph{graded Jacobi identity}:
        \[
        (-1)^{(|f|+1)(|h|+1)}[[f,g],h]+(-1)^{(|g|+1)(|f|+1)}[[g,h],f]+(-1)^{(|h|+1)(|g|+1)}[[h,f],g]=0.
        \]
    \end{itemize}
\end{Lem}

According to Lemma~\ref{lemma: properties of Ger bracket} (1), the Gerstenhaber bracket on $C^*(A)$ induces a Gerstenhaber bracket on $\HH^*(A)$. It follows from Lemma~\ref{lemma: properties of Ger bracket} (2) and (3) that $(s\HH^*(A),[-,-])$ is a graded Lie algebra.\\

The cup product and the Gerstenhaber bracket are related by the following \emph{Leibniz rule}.
\begin{Lem}
    For all $f\in \rmHom_{E^e}((s\oA)^{\otimes m},A)$, $g\in \rmHom_{E^e}((s\oA)^{\otimes n},A)$ and $h\in \rmHom_{E^e}((s\oA)^{\otimes r},A)$, one has
    \[
    [f,g\smile h]=[f,g]\smile h+(-1)^{(|f|+1)|g|}g\smile [f,h].
    \]
\end{Lem}

\medskip

\textbf{Gerstenhaber algebras.} A \emph{Gerstenhaber algebra} 
$H^{\bullet}=(H^{\bullet},\cdot,[-,-])$ is both a graded commutative algebra $(H^{\bullet},\cdot)$ and a graded Lie algebra $(sH^{\bullet},[-,-])$ such that the graded Leibniz rule holds.
\begin{Thm}
    Let $A$ be an augmented dg $E$-algebra. Then $(\HH^*(A),\smile,[-,-])$ is a Gerstenhaber algebra.
\end{Thm}

\smallskip

\smallskip

\subsection{Graded and two-sided CS-projective resolution}\label{section: cs-projective resolution}\

\medskip

Recall that we fix a graded skew-gentle triple $(Q,R,\Sp)$, and set $A=\bbK Q/I$ the corresponding graded skew-gentle algebra, where $I$ denotes the two-sided ideal of $\bbK Q$ generated by $R\cup \{\varepsilon^2-\varepsilon\ |\ \varepsilon\in \Sp\}$. We present a graded version of the projective resolution of skew-gentle algebras introduced by S.~Chouhy and the third author  in~\cite{CS2015}. 
The CS-resolution can be viewed as a deformation of Bardzell resolution for monomial algebras and is applicable to graded skew-gentle algebras.

As recalled in Section~\ref{section: TT calculus of dg algebra}, the Hochschild (co)homology of graded skew-gentle algebras can be computed using the  Hochschild (co)chain complex induced by the two-sided bar resolution, where the operations and algebraic structures are naturally defined. However, this resolution is often too large to be practical for explicit computations. To address this issue, we instead employ the graded CS-projective resolution of $A$, which is significantly smaller and more suitable for concrete calculations.

Let $\Gamma=\cup_{n\in \mathbb{N}}\Gamma_n$, where $\Gamma_n $ is the set consisting of the paths of length $n$ in $Q$ all of whose subpaths of length $2$ are in $S=R\cup \{ \varepsilon^2\ |\ \varepsilon\in \Sp \}$. So, we have $\Gamma_0=Q_0$ and $\Gamma_1=Q_1$. 
For simplicity, we write $\{\varepsilon^n\ |\ \varepsilon\in \Sp\}$ as $\Sp^n$, note that $\Sp^n$ is a subset of $\Gamma_n$ for $n\ge 1$. 
In case $\Sp=\emptyset$, we have that $A$ is a graded gentle algebra and the definition of $\Gamma_n$ coincides with the definition of $\Gamma_n$ in \cite{CSSS24}. 

\begin{Prop}\label{CS-reso-skew}
    Let $(Q,R,\Sp)$ be a (graded) skew-gentle triple and $A=\bbK Q/\langle R\cup \{\varepsilon^2-\varepsilon\ |\ \varepsilon\in \Sp\} \rangle$ be the corresponding (graded) skew-gentle algebra. Let $\Gamma_n$ be defined as in the previous paragraph for each $n\geq 0$. The \emph{graded CS-projective resolution} $\bbP$ of $A$ as graded $A$-bimodule is as follows
\[
\cdots\xrightarrow{}A\otimes \bbK \Gamma_n \otimes A \xrightarrow{d_n} A\otimes \bbK\Gamma_{n-1} \otimes A \xrightarrow{}  \cdots\xrightarrow{} A\otimes \bbK \Gamma_1\otimes A \xrightarrow{d_1} A\otimes \bbK \Gamma_0\otimes A \xrightarrow{d_0} A\xrightarrow{} 0,
\]
where the differentials $d_n$ are $A$-bimodule morphisms given by
\[d_n(1\otimes q\otimes 1)=
\left\{\begin{array}{ll}
	 a_n\otimes a_{n-1} \cdots a_1\otimes 1+(-1)^n 1\otimes a_n\cdots a_2\otimes a_1 &\text{if $q=a_n\cdots a_1\in \Gamma_n\setminus\Sp^n$,}\\
     
    \varepsilon\otimes \varepsilon^{n-1}\otimes 1-1\otimes \varepsilon^{n-1}\otimes \varepsilon &\text{if $q=\varepsilon^n\in \Sp^n$ and $n$ is odd,}\\
    
    \varepsilon\otimes \varepsilon^{n-1}\otimes 1+1\otimes \varepsilon^{n-1}\otimes \varepsilon-1\otimes \varepsilon^{n-1}\otimes 1 & \text{if $q=\varepsilon^n\in \Sp^n$ and $n$ is even}.
  \end{array}
  \right.\]
In particular, $d_0(1\otimes e\otimes 1)=e$ for all $e\in \Gamma_0$.
\end{Prop}

\begin{Proof}
  By \cite[Theorem 4.1]{CS2015} we know that for any $q\in \Gamma_n$, $d_n(1\otimes q\otimes 1)$ consists of two parts, one  induced from Bardzell resolution and another one due to the reduction. Since the algebra we are considering is skew-gentle, that reduction has a particularly simple form, which we now describe. 
  
  Let $r$ be a basic reduction, if $q=a_n\cdots a_1\in \Gamma_n\setminus\Sp^n$ then $r(q)=0$, hence the only non zero terms of $d_n(1\otimes q\otimes 1)$ come from the Bardzell resolution, that is, $d_n(1\otimes q\otimes 1)=a_n\otimes a_{n-1} \cdots a_1\otimes 1+(-1)^n 1\otimes a_n\cdots a_2\otimes a_1$.
  Otherwise $q\in \Sp^n$, say $q=\varepsilon^n$ for $\varepsilon\in \Sp$, thus $r(q)=\varepsilon^{n-1}$. By \cite[Theorem 4.1]{CS2015}, we can assume that $d_n(1\otimes \varepsilon^n\otimes 1)=\varepsilon\otimes \varepsilon^{n-1}\otimes 1+(-1)^n 1\otimes \varepsilon^{n-1}\otimes \varepsilon+\lambda_n 1\otimes\varepsilon^{n-1}\otimes 1$ where $\lambda_n\in \bbK$. 
  
  From the equality $d_{n-1}\circ d_n=0$ to $1\otimes \varepsilon^n\otimes 1$, it is easy to get \[
\lambda_n=
\left\{
\begin{array}{ll}
-1 &\text{if $n$ is even},\\
0 &\text{otherwise}.
\end{array}
\right.
\] 
\end{Proof}

After applying $\rmHom_{A^e}(-,A)$ to  $\bbP$ we obtain the following cochain complex $\rmHom_{A^e}(\bbP,A)$,
\[
0 \to\rmHom_{A^e}(A\otimes \bbK \Gamma_0\otimes A, A)\xto{d^0}  \cdots \to \rmHom_{A^e}(A \ot\bbK \Gamma_n \ot A, A)
\xto{d^{n}} \rmHom_{A^e}(A\ot \bbK \Gamma_{n+1} \ot A, A)\to \cdots,
\]
with differentials given by $d^n(f)=f\circ d_{n+1}$. Since $d_*$ preserves the grading, the same is true for $d^*$. As a result, the complex can be decomposed as 
\[
\rmHom_{A^e}(\bbP,A)=\bigoplus_{j\in \bbZ} \rmHom_{A^e}(\bbP,A)^j.
\]
This decompositon provides, after computing Hochschild cohomology, a decomposition of each $n$-th Hochschild cohomology  space $\rmH^n(\rmHom_{A^e}(\bbP,A))\cong \Ker(d^n)/\mathrm{Im}(d^{n-1})$. More precisely, the   
\emph{graded Hochschild cohomology 
of $A$ is given by \[\rmH^n(\rmHom_{A^e}(\bbP,A))=\bigoplus_{j\in \bbZ}\rmH^n(\rmHom_{A^e}(\bbP,A)^j)=:\bigoplus_{j\in \bbZ} \HH^{n,j}(A), \] where $\HH^{n,\sbl}(A)$ denotes $\rmH^n(\rmHom_{A^e}(\bbP,A)^{\sbl})$.} Moreover, for each element $f\in \rmHom_{A^e}(\bbP_n,A)^j$, we will call $j$ the \emph{internal degree} of $f$, denoted by $|f|=j$, and $n$ the \emph{external degree} of $f$.

Using parallel paths, introduced by C. Cibils \cite{Cib90}, the complex $\rmHom_{A^e}(\bbP,A)$ is better understood. Two paths $p$ and $q$ are called \emph{parallel}, denoted by $p\| q$, if they have the same source and target. Let $X$ and $Y$ be the set of paths in $Q$, define $X\| Y:=\{(p, q)\in X\times Y \ | \ p \| q\}$ and $\bbK (X\| Y)$ to be the $\bbK$-vector space with basis $X\| Y$.

\begin{Lem}\cite{Cib90}
For each $n\ge 0$ there is isomorphism of $\bbK$-vector spaces: 
\[
\rmHom_{A^e}(\bbP_n,A)=\rmHom_{A^e}(A\ot \bbK \Gamma_n \ot A,A)\cong \bbK( \Gamma_n\| \mathcal \calB),
\]
where for $(\gamma, b)\in \Gamma_n\| \mathcal \calB$, the corresponding map $f_{(\gamma, b)}: \bbP_n=A\ot \bbK \Gamma_n\ot A\to A$ sends $p\ot \gamma'\ot q\in \bbP_n$ to $(-1)^{|p|(|b|-|\gamma|)}\ \delta_{\gamma, \gamma'} \cdot \pi(pbq)$, where $p,\, q\in \calB$ and $\gamma'\in \Gamma_n$.
\end{Lem}

After this lemma, we give a new description of the complex $\rmHom_{A^e}(\bbP,A)$.
\begin{Prop}\label{CS-cochain-complex}
Let $(Q,R,\Sp)$ be a graded skew-gentle triple and let $A=\bbK Q/\langle R\cup \{\varepsilon^2-\varepsilon\ |\ \varepsilon\in \Sp\} \rangle$ be the corresponding graded skew-gentle algebra. Let $\Gamma_n$ be defined as above for each $n\geq 0$. The graded Hochschild cohomology of $A$ is the cohomology of the following complex
\[
\bbK(\Gamma\| \calB):\; 0\xrightarrow{}\bbK (\Gamma_0\| \calB)\xrightarrow{d^0}\bbK (\Gamma_1\|  \calB)\xrightarrow{}\cdots \xrightarrow{} \bbK (\Gamma_n\|  \calB) \xrightarrow{d^n} \bbK (\Gamma_{n+1}\|  \calB)\xrightarrow{} \cdots
\]
with differential $d^n$ given by 
\[ d^n(q, \alpha)=
\begin{cases}
    \sum\limits_{\substack{b\in Q_1\\ bq\in \Gamma_{n+1}}} (-1)^{|b|(|\alpha|-|q|)}(bq,\pi(b\alpha))-(-1)^n \sum\limits_{\substack{a\in Q_1\\ q a\in \Gamma_{n+1}}} (q a,\pi(\alpha a)) & \text{if $q\in \Gamma_n\setminus \Sp^n$},\\
    (\varepsilon^{n+1},\pi(\varepsilon\alpha))-(\varepsilon^{n+1},\pi(\alpha\varepsilon)) & \text{if $q=\varepsilon^n\in {\Sp}^n$ and $n$ is even},\\
    (\varepsilon^{n+1},\pi(\varepsilon\alpha))+(\varepsilon^{n+1},\pi(\alpha\varepsilon))-(\varepsilon^{n+1},\pi(\alpha)) & \text{if $q=\varepsilon^n\in {\Sp}^n$ and $n$ is odd},
\end{cases}\]
for all $(q,\alpha)\in \Gamma_n\|  \calB$.
\end{Prop}

Therefore, $\HH^{n,\sbl}(A)\cong (\Ker(d^n)/\mathrm{Im}(d^{n-1}))^{\sbl}$. 
However, if we view $A$ as a dg algebra with trivial differential, then we can also define the Hochschild cohomology $\HH^*(A)$ of $A$ as in Section~\ref{section: TT calculus of dg algebra}. We will explore the relationship between $\HH^{\bullet,\sbl}(A)$ and $\HH^*(A)$.

\begin{Def}\label{def: two-sided cs-proj. resol}

Consider the following bicomplex
\[
\xymatrix@R=1.2pc{
&{} &-3 & -2 & -1 & 0\\
&{} & \vdots & \vdots & \vdots & \vdots & {} \\
1& \cdots \ar[r] & \bbP_3^1\ar[u] \ar[r]^{d_{3,1}} & \bbP_2^1 \ar[u] \ar[r]^{d_{2,1}} & \bbP_1^1 \ar[u] \ar[r]^{d_{1,1}} & \bbP_0^1\ar[u] \ar[r] & 0  \\
0& \cdots \ar[r] & \bbP_3^0\ar[u] \ar[r]^{d_{3,0}} & \bbP_2^0 \ar[u] \ar[r]^{d_{2,0}}  & \bbP_1^0 \ar[u] \ar[r]^{d_{1,0}} & \bbP_0^0\ar[u] \ar[r] & 0\\
-1& \cdots \ar[r] & \bbP_3^{-1}\ar[u] \ar[r]^{d_{3,-1}} & \bbP_2^{-1} \ar[u] \ar[r]^{d_{2,-1}} & \bbP_1^{-1} \ar[u] \ar[r]^{d_{1,-1}} & \bbP_0^{-1} \ar[u] \ar[r] & 0\\
&{}&\vdots \ar[u] & \vdots \ar[u] & \vdots \ar[u] & \vdots \ar[u] & {}
}
\]
where the superscript $j$ in $\bbP_n^j$ denotes the degree of the elements, with the horizontal differential given by the restrictions $d_{n,j}:=d_n|_{\bbP_n^j}$ of $d_n$ for $j\in \bbZ$ and the vertical differentials are identically $0$. We denote by $(T_*,t_*)$ the total complex associated to this bicomplex, called \emph{two-sided CS-projective resolution}, that is,
\[
T_N=\bigoplus_{\substack{N=j-n\\ n\ge 0,j\in\bbZ}} \bbP_n^j, \text{ and } t_N=\sum_{\substack{N=j-n\\ n\ge 0,j\in\bbZ}} d_{n,j}.
\]
In particular, for the ungraded quiver algebra, the complex $T_*$ coincides with the CS-projective resolution.

The complex $(\calH om_{A^e}(T_*,A),\varkappa)$ is defined as follows: for any $N\in \bbZ$
\begin{align*}
\calH om_{A^e}^N(T_*,A)
&=\prod_{M\in\bbZ} \rmHom_{A^e}(T_M,A^{M+N})\\
&=\prod_{\substack{N=n+j \\n\ge 0,\, j\in\bbZ}} \rmHom_{A^e}(\bbP_n,A)^{j},
\end{align*}
and the differential $\varkappa$ is defined by
\[
\varkappa(f)=(-1)^{n+j-1}f\circ d_{n+1},
\]
for $f\in \rmHom_{A^e}(\bbP_n,A)^j$.

\end{Def}

In the next section, we will construct a comparison map between the two-sided bar resolution and the two-sided CS-projective resolution, allowing us to transfer algebraic structures such as the cup product and Gerstenhaber bracket defined on the former to the latter. 
According to Definition~\ref{def: two-sided cs-proj. resol}, we have
\[
\HH^N(A)=\rmH^N(\calH om_{A^e}(T_*,A))\cong \prod_{\substack{N=n+j\\ n\ge 0,\, j\in \bbZ}} \rmH^n(\rmHom_{A^e}(\bbP_n,A)^j)= \prod_{\substack{N=n+j\\ n\ge 0,\, j\in \bbZ}}\HH^{n,j}(A),
\]
where 
$\HH^{n,j}(A)\cong \Ker(d^{n,j})/\mathrm{Im}(d^{n-1,j})\cong(\Ker(d^{n})/\mathrm{Im}(d^{n-1}))^j$ with $d^{n,j}:=\rmHom_{A^e}(d_{n,j},A)\cong d^n|_{(\bbK(\Gamma_n\|\calB))^j}.$

For an element $f\in \rmHom_{A^e}(\bbP_n,A)^j$, we call 
$N=n+j$ the \emph{total degree} of $f$ denoted by $\deg(f)$, that is, the sum of the internal degree 
$|f|=j$ and the external degree $n$ of $f$.

However, we observe that the signs of the differentials are different, that is, we have 
\[
\varkappa(f)= (-1)^{\deg(f)+1} f\circ d_{n+1}= (-1)^{\deg(f)+1} d^n(f),
\]
this difference does not affect the computation of the generators of Hochschild cohomology as vector spaces, but only impacts on its algebra structure. Therefore, we will use the complex \( (\rmHom_{A^e}(\bbP, A), d^*) \) to compute the generators of \( \mathrm{HH}^*(A) \) as a graded vector space and the complex \( (\mathcal{H}om_{A^e}(T_*, A), \varkappa) \) when considering the algebra structure.

\subsection{Comparison morphisms}\label{section: comparison-map}\

In this section, let $A=\bbK Q/\langle R\cup \{\varepsilon^2-\varepsilon\ |\ \varepsilon\in \Sp\} \rangle$ be a graded skew-gentle algebra, we will construct the comparison morphisms between the two-sided bar resolution $B(A,A,A)$ and the two-sided CS-projective resolution $T_*$ of $A$.

A comparison morphism between two projective resolutions of a graded algebra $A$, is a morphism of chain complexes that lifts the identity map on $A$. Applying $\rmHom_{A^e}(-,A)$ to the comparison morphism yields a quasi-isomorphism. For graded algebras, the comparison morphisms preserve the grading.

We will define maps $F:\; T_* \to B(A,A,A)$ and $G: \; B(A,A,A) \to T_*$ and then prove that they are in fact comparison morphisms.

For any $N\in\bbZ$, let
\[
F_N=(F_N^n)_{n\ge 0}:\; T_N=\bigoplus_{\substack{N=r-n\\ n\ge 0,r\in\bbZ}} \bbP_n^r\to B(A,A,A)_N=\bigoplus_{n\ge 0} (A\ot (s\oA)^{\ot n} \ot A)^N
\]
be such that  $F_N^n:\; \bbP_n^{N+n}=(A\ot \bbK \Gamma_n\ot A)^{N+n} \to (A\ot (s\oA)^{\ot n} \ot A)^N$ is the graded $A$-bimodule extension (i.e., for $a\ot q\ot b\in (A\ot \bbK\Gamma_n\ot A)^{N+n}$, $F_N^n(a\ot q\ot b)=(-1)^{n|a|}a F_{|q|-n}^n(1\ot q\ot 1)b$) of the map defined as follows.

Let $F_N^0=\mathrm{id}_{(A\otimes A)^N}$ and, for $n\ge 1$, let $1\ot q\ot 1$ be an element of $1\otimes \Gamma_n\otimes 1$. Let $q=a_n\cdots a_1\in \Gamma_n$ and $\sigma(a)=\sum\limits_{i=2}^{n}(|a_i|+\cdots+|a_n|)$, define 
\[
F_{|q|-n}^n(1\otimes q \otimes 1):=(-1)^{\sigma(a)}\ 1\otimes s a_n\otimes \cdots \otimes s a_1\otimes 1.
\]

\begin{Prop}
The map  $F:\; T_* \to B(A,A,A)$ is a comparison morphism.
\end{Prop}
\begin{Proof}
    In order to prove that $F$ is a comparison morphism, one should check that the diagram
         \[
         \xymatrix@R=1.2pc{
         &(A\otimes_E \bbK\Gamma_n\otimes_E A)^{N+n} \ar[dd]^{F_N^n} \ar[rr]^{d_n} && (A\otimes_E \bbK\Gamma_{n-1}\otimes_E A)^{N+n} \ar[dd]^{F_{N+1}^{n-1}}\\
         \\
         &(A\otimes_E (s\oA)^n \otimes_E A)^N \ar[rr]^{\delta_1} && (A\otimes_E (s\oA)^{n-1}\otimes_E A)^{N+1}}
        \]
is commutative for $n\ge 1$. This is not difficult for elements in $\Gamma_n\setminus {\Sp}^n$. 

Let $\varepsilon^n\in {\Sp}^n$, we need to consider two cases: $n$ odd and $n$ even. We only provide the proof for the case where $n$ is odd since the other case ($n$ even) can be derived in a similar manner.

Assume that $n$ is odd. We have
    \[
    F_{-n}^n(1\otimes \varepsilon^n\otimes 1)=1\otimes\underbrace{ s\varepsilon \otimes \cdots \otimes s\varepsilon}_{n} \otimes 1
    \]
    and
    \begin{align*}
    \delta_1(F_{-n}^n(1\otimes \varepsilon^n\otimes 1)) 
    &= \varepsilon\otimes\underbrace{ s\varepsilon \otimes \cdots \otimes s\varepsilon}_{n-1} \otimes 1 
    + \sum_{i=1}^{n-1}(-1)^i\  1\otimes\underbrace{ s\varepsilon \otimes \cdots \otimes s\varepsilon}_{n-1} \otimes 1 \\
    &\quad - 1\otimes\underbrace{ s\varepsilon \otimes \cdots \otimes s\varepsilon}_{n-1} \otimes \varepsilon \\
    &= \varepsilon\otimes\underbrace{ s\varepsilon \otimes \cdots \otimes s\varepsilon}_{n-1} \otimes 1 
    - 1\otimes\underbrace{ s\varepsilon \otimes \cdots \otimes s\varepsilon}_{n-1} \otimes \varepsilon.
    \end{align*}

    On the other hand, $n$ odd implies 
    \[
    d_n(1\otimes \varepsilon^n\otimes 1)=\varepsilon\otimes \varepsilon^{n-1}\otimes 1-1\otimes \varepsilon^{n-1}\otimes \varepsilon
    \]
    and
    \[
    F_{-n+1}^{n-1}(d_n(1\otimes \varepsilon^n\otimes 1))=\varepsilon\otimes\underbrace{ s\varepsilon \otimes \cdots \otimes s\varepsilon}_{n-1} \otimes 1 
    - 1\otimes\underbrace{ s\varepsilon \otimes \cdots \otimes s\varepsilon}_{n-1} \otimes \varepsilon.
    \]
    Therefore, $F_{-n+1}^{n-1}(d_n(1\otimes \varepsilon^n\otimes 1))= \delta_1(F_{-n}^n(1\otimes \varepsilon^n\otimes 1))$. This completes the proof.
\end{Proof}

Next, we define the map $G: B(A,A,A) \to T_*$ and show that it is a comparison map. Firstly, we introduce the following  $\bbK$-basis of $A\otimes (s\oA)^{\otimes n}\otimes A$:
\[
\calB_{\bbB_n} =
\left\{
\alpha_{n+1} \ot s\alpha_n \ot \dots \ot s\alpha_1\ot \alpha_0
\;\middle|\;
\begin{array}{ll}
\text{(i)} & \alpha_j \in \calB, \text{ for all } j;\ l(\alpha_j) \geq 1,\ \text{for } j = 1, \dots, n; \\
\text{(ii)} & t(\alpha_j) = s(\alpha_{j+1}),\ \text{for } j = 0, \dots, n
\end{array}
\right\}.
\]
Secondly, let 
\[
\calB'_{\bbB_n} =
\left\{
1  \ot s\alpha_n \ot \dots \ot s\alpha_1\ot  1
\;\middle|\;
\begin{array}{ll}
\text{(i)} & \alpha_j \in \calB \text{ and }  l(\alpha_j) \geq 1,\ \text{for } j = 1, \dots, n; \\
\text{(ii)} & t(\alpha_j) = s(\alpha_{j+1}),\ \text{for } j = 1, \dots, n-1
\end{array}
\right\}.
\]
Clearly, $\calB'_{\bbB_n}\subset \calB_{\bbB_n}$ and $\calB'_{\bbB_n}$ generates $A\otimes (s\oA)^{\otimes n}\otimes A$ as a graded $A$-bimodule.

For any $N\in\bbZ$, let
\[
G_N=(G_N^n)_{n\ge 0}:\; B(A,A,A)_N=\bigoplus_{n\ge 0} (A\ot (s\oA)^{\ot n} \ot A)^N \to T_N=\bigoplus_{\substack{N=r-n\\ n\ge 0,r\in\bbZ}} \bbP_n^r
\]
where $G_N^n:\; (A\ot (s\oA)^{\ot n} \ot A)^N \to \bbP_n^{N+n}=(A\ot \bbK \Gamma_n\ot A)^{N+n}$ is the graded $A$-bimodule extension of the following map.

For $n=0$, define $G_N^0=\mathrm{id}_{(A\otimes A)^N}$ and, for $n=1$ and $1\otimes s\alpha\otimes 1\in \calB'_{\bbB_1}$, with $\alpha=a_l\cdots a_1$ and $a_i\in Q_1$, let
\[
G_{|\alpha|-1}^1(1\otimes s\alpha\otimes 1):=\sum_{i=1}^l  a_l\cdots a_{i+1}\otimes  a_i\otimes  a_{i-1}\cdots  a_1.
\]
For $n>1$ and $q=1  \ot s\alpha_n \ot \dots \ot s\alpha_1\ot  1\in \calB'_{\bbB_n}$ with $|q|=\sum\limits_{i=1}^n(|\alpha_i|-1)$, assume $\alpha_i=a_{l(\alpha_i)}^i\cdots a_1^i$ with each $a_j^i\in Q_1$, if $ a_1^n \alpha_{n-1} \dots \alpha_2 a_{l(\alpha_1)}^1\in \Gamma_n$, let 
\[
G_{|q|}^n(q):=(-1)^{\sigma(\alpha)} a_{l(\alpha_n)}^n \cdots a_2^n \otimes a_1^n \alpha_{n-1} \dots \alpha_2 a_{l(\alpha_1)}^1 \otimes  a_{l(\alpha_1)-1}^1 \cdots a_1^1,
\]
where $\sigma(\alpha)=\sum\limits_{i=2}^n (|\alpha_i|+\cdots+|\alpha_n|)$, in particular, if $a_1^n \alpha_{n-1} \dots \alpha_2 a_{l(\alpha_1)}^1=\varepsilon^n \in {\Sp}^n$, by definition, we know that
\[
G_{|\alpha_1|+|\alpha_n|}^n(q):=(-1)^{(n-1)|\alpha_n|} a_{l(\alpha_n)}^n \cdots a_2^n \otimes \varepsilon^n \otimes a_{l(\alpha_1)-1}^1 \cdots a_1^1,
\]
for all other cases, $G_{|q|}^n(q):=0$.

\begin{Prop}
    $G:  B(A,A,A) \to T_*$ is a comparison morphism. Moreover, $GF=\mathrm{id}$ and $FG$ is homotopic  to $\mathrm{id}$.
\end{Prop}
\begin{Proof}
    The commutativity of the first diagram
    \[
    \xymatrix@R=1.2pc{
    (A\otimes_E s\oA \otimes_E A\ar[dd]^{G_N^1})^N \ar[rr]^{\delta_1} && (A\otimes A)^{N+1}\ar@{=}[dd]^{} \ar[r] & A\\
    \\
   ( A\otimes_E \bbK \Gamma_1 \otimes_E A)^{N+1} \ar[rr]^{d_1} && (A\otimes A)^{N+1} \ar[r] & A
    }
    \]
    
    follows immediately by calculating $d_1\circ G_N^1$. We also need to  prove the commutativity of the following diagram:
    
    \[
     \xymatrix@R=1.2pc{
    (A\otimes_E (s\oA)^{n} \otimes_E A)^N\ar[dd]^{G_N^n} \ar[rr]^{\delta_1} 
    && (A\otimes_E (s\oA)^{n-1} \otimes_E A)^{N+1}\ar[dd]^{G_{N+1}^{n-1}} \\
    \\
    (A\otimes_E \bbK \Gamma_{n} \otimes_E A)^{N+n} \ar[rr]^{d_{n}} 
    &&  (A\otimes_E \bbK \Gamma_{n-1} \otimes_E A)^{N+n}
    }
    \]
    Let $q=1  \ot s\alpha_n \ot \cdots \ot s\alpha_1\ot  1\in \calB'_{\bbB_n}$, with $\alpha_i=a_{l(\alpha_i)}^i \cdots a_1^i$ for each $a_j^i\in Q_1$, 
    \begin{itemize}
    
    \item [(1)] if $ a_1^n \alpha_{n-1} \dots \alpha_2 a_{l(\alpha_1)}^1 \in \Gamma_n\setminus {\Sp}^n$, then 
    \begin{align*}
        d_n (G_{|q|}^n(q))
        =& d_n((-1)^{\sigma(\alpha)} a_{l(\alpha_n)}^n \cdots a_2^n \ot a_1^n \alpha_{n-1} \dots \alpha_2 a_{l(\alpha_1)}^1 \ot  a_{l(\alpha_1)-1}^1 \cdots a_1^1)\\
        =& (-1)^{\sigma(\alpha)}\alpha_n\ot \alpha_{n-1}\cdots \alpha_2 a_{l(\alpha_1)}^1 \ot a_{l(\alpha_1)-1}^1 \cdots a_1^1\\
        &+(-1)^{\sigma(\alpha)-n} a_{l(\alpha_n)}^n \cdots a_2^n \ot a_1^n\alpha_{n-1}\cdots \alpha_2\ot \alpha_1
    \end{align*}
    and
    \begin{align*}
        G_{|q|+1}^{n-1}(\delta_1 (q))
        =& G_{|q|+1}^{n-1}(\alpha_n\ot s\alpha_{n-1}\ot\cdots \ot s\alpha_1\ot 1+(-1)^ {\sum\limits_{i=2}^n|\alpha_i|-n} 1\ot s\alpha_n\ot\cdots\ot s\alpha_2\ot \alpha_1)\\
        =& (-1)^{(n-1)|\alpha_n|}\alpha_n
        G_{|q|-|\alpha_n|+1}^{n-1}(1 \ot s\alpha_{n-1}\ot\cdots \ot s\alpha_1\otimes 1)\\
        &+(-1)^{\sum\limits_{i=2}^n|\alpha_i|-n} G_{|q|-|\alpha_1|+1}^{n-1} (1\ot s\alpha_n\ot\cdots\ot s\alpha_2\otimes 1) \alpha_1\\
        =& (-1)^{\sigma(\alpha)}\alpha_n\ot \alpha_{n-1}\cdots \alpha_2 a_{l(\alpha_1)}^1 \ot a_{l(\alpha_1)-1}^1 \cdots a_1^1\\
        &+(-1)^{\sigma(\alpha)-n} a_{l(\alpha_n)}^n \cdots a_2^n \ot a_1^n\alpha_{n-1}\cdots \alpha_2\ot \alpha_1\\
        =& d_n (G_{|q|}^n(q));
    \end{align*}

    \item [(2)] if $a_1^n \alpha_{n-1} \dots \alpha_2 a_{l(\alpha_1)}^1=\varepsilon^n \in {\Sp}^n$, then
    \begin{align*}
        d_n (G_{|q|}^n(q))
        =& d_n((-1)^{(n-1)|\alpha_n|} a_{l(\alpha_n)}^n \cdots a_2^n \ot \varepsilon^n \ot a_{l(\alpha_1)-1}^1 \cdots a_1^1)\\
        =& (-1)^{(n-1)|\alpha_n|} a_{l(\alpha_n)}^n \cdots a_2^n
        (\varepsilon\ot \varepsilon^{n-1}\ot 1+(-1)^n1\ot \varepsilon^{n-1}\ot \varepsilon+\lambda_n\; 1\ot \varepsilon^{n-1}\ot 1)
        a_{l(\alpha_1)-1}^1 \cdots a_1^1\\
        =&(-1)^{(n-1)|\alpha_n|} \alpha_n\ot \varepsilon^{n-1}\ot a_{l(\alpha_1)-1}^1 \cdots a_1^1\\
        &+(-1)^{(n-1)|\alpha_n|+n} a_{l(\alpha_n)}^n \cdots a_2^n\ot\varepsilon^{n-1}\ot \alpha_n\\
        &+ (-1)^{(n-1)|\alpha_n|}\lambda_n a_{l(\alpha_n)}^n \cdots a_2^n \ot\varepsilon^{n-1}\ot a_{l(\alpha_1)-1}^1 \cdots a_1^1,
    \end{align*}
    and
    \begin{align*}
    G_{|q|+1}^{n-1}(\delta_1 (q))
    =& G_{|q|+1}^{n-1}(\alpha_n\ot \underbrace{s\varepsilon \ot \cdots \ot s\varepsilon}_{n-2} \ot s\alpha_1\ot 1
    +(-1)^{|\alpha_n|-n} 1\ot s\alpha_n\ot \underbrace{s\varepsilon\ot \cdots\ot s\varepsilon}_{n-2}\ot \alpha_1\\
    &\quad \quad\quad +(-1)^{|\alpha_n|}\lambda_n\; 1\ot s\alpha_n\ot \underbrace{s\varepsilon\ot \cdots\ot s\varepsilon}_{n-3}\ot s\alpha_1\ot 1)\\
    =& (-1)^{(n-1)|\alpha_n|}\alpha_n
    G_{|q|-|\alpha_n|+1}^{n-1}(1 \ot \underbrace{s\varepsilon \ot \cdots \ot s\varepsilon}_{n-2} \ot s\alpha_1\ot 1)\\
    &+(-1)^{|\alpha_n|+n} G_{|q|-|\alpha_1|+1}^{n-1} (1\ot s\alpha_n\ot \underbrace{s\varepsilon\ot \cdots\ot s\varepsilon}_{n-2}\ot 1) \alpha_1\\
    &+(-1)^{|\alpha_n|}\lambda_n G_{|q|+1}^{n-1}(1\ot s\alpha_n\ot \underbrace{s\varepsilon\ot \cdots\ot s\varepsilon}_{n-3}\ot s\alpha_1\ot 1) \\
    =&(-1)^{(n-1)|\alpha_n|} \alpha_n\ot \varepsilon^{n-1}\ot a_{l(\alpha_1)-1}^1 \cdots a_1^1\\
    &+(-1)^{(n-1)|\alpha_n|+n} a_{l(\alpha_n)}^n \cdots a_2^n\ot\varepsilon^{n-1}\ot \alpha_n\\
    &+ (-1)^{(n-1)|\alpha_n|}\lambda_n a_{l(\alpha_n)}^n \cdots a_2^n \ot\varepsilon^{n-1}\ot a_{l(\alpha_1)-1}^1 \cdots a_1^1\\
    =& d_n (G_{|q|}^n(q))
    \end{align*}  
    where 
    \[
    \lambda_n=
    \left\{
    \begin{array}{ll}
    -1 &\text{if $n$ is even},\\
    0 &\text{otherwise}.
    \end{array}
    \right.
    \]

     \item [(3)] if $a_1^n \alpha_{n-1} \dots \alpha_2 a_{l(\alpha_1)}^1\not \in \Gamma_n$, by definition of $G$,
     \[
     G_{|q|}^n(q)=0, \text{ then } d_n(G_{|q|}^n(q))=0,
     \]
     it is enough to prove that  $G_{|q|+1}^{n-1}(\delta_1(q))=0$. In fact,
     \begin{align*}
        \delta_1 (q)
         =&  \alpha_n\ot s\alpha_{n-1}\ot\cdots\ot s\alpha_1\ot 1+ (-1)^{|\alpha_n|-1} 1\ot s(\alpha_n \alpha_{n-1})\ot s\alpha_{n-2}\ot\cdots\ot s\alpha_1\ot 1\\
         &+ \sum_{i=3}^{n-1} (-1)^{\epsilon_{i,n}} 1\ot s\alpha_n\ot\cdots\ot s\alpha_{i+1}\ot s(\alpha_i\alpha_{i-1})\ot s\alpha_{i-2}\ot \cdots\ot s\alpha_1\ot 1\\
         &+(-1)^{\epsilon_{2,n} }1\ot s\alpha_n\ot \cdots\ot s\alpha_3\ot s(\alpha_2\alpha_1)\ot 1+(-1)^{\epsilon_{2,n}-1} 1\ot s\alpha_n\ot \cdots\ot s\alpha_2\ot \alpha_1,
     \end{align*}
    for $2\le i\le n-2$, $G_{|q|+1}^{n-1}(1\ot s\alpha_n\ot\cdots\ot s\alpha_{i+1}\ot s(\alpha_i\alpha_{i-1})\ot s\alpha_{i-2}\ot \cdots\ot s\alpha_1\ot 1)=0$ by the definition of $G$ and $a_1^n \alpha_{n-1} \dots \alpha_2 a_{l(\alpha_1)}^1\not \in \Gamma_n$.

    Now we will prove that $G_{|q|+1}^{n-1}( \alpha_n\ot s\alpha_{n-1}\ot\cdots\ot s\alpha_1\ot 1) = (-1)^{|\alpha_n|} G_{|q|+1}^{n-1}(1\ot s(\alpha_n \alpha_{n-1})\ot s\alpha_{n-2}\ot\cdots\ot s\alpha_1\ot 1)$. There are two cases. 
    
    If $a_1^{n-1} \alpha_{n-2}\cdots\alpha_2 a_{l(\alpha_1)}^1\not\in \Gamma_{n-1}$, then
    \[
    G_{|q|+1}^{n-1}( \alpha_n\ot s\alpha_{n-1}\ot\cdots\ot s\alpha_1\ot 1)=0 = (-1)^{|\alpha_n|} G_{|q|+1}^{n-1}(1\ot s(\alpha_n \alpha_{n-1})\ot s\alpha_{n-2}\ot\cdots\ot s\alpha_1\ot 1),
    \]
    while if $a_1^{n-1} \alpha_{n-2}\cdots\alpha_2 a_{l(\alpha_1)}^1 \in \Gamma_{n-1}$, then
    \begin{align*}
        G_{|q|+1}^{n-1}( \alpha_n\ot s\alpha_{n-1}\ot\cdots\ot s\alpha_1\ot 1)
        =& (-1)^{(n-1)|\alpha_n|} \alpha_n G_{|q|+1-|\alpha_n|}^{n-1}( 1\ot s\alpha_{n-1}\ot\cdots\ot s\alpha_1\ot 1)\\
        =& (-1)^{\sigma(\alpha)}\alpha_n a_{l(\alpha_{n-1})}^{n-1}\cdots a_2^{n-1} \ot a_1^{n-1}\alpha_{n-1} \cdots \alpha_2 a_{l(\alpha_1)}^1\ot a_{l(\alpha_1)-1}^1 \cdots a_1^1\\
        =&(-1)^{|\alpha_n|} G_{|q|+1}^{n-1}(1\ot s(\alpha_n \alpha_{n-1})\ot s\alpha_{n-2}\ot\cdots\ot s\alpha_1\ot 1).
    \end{align*}
   Similarly, we can prove that 
   \[
   G_{|q|+1}^{n-1}(1\ot s\alpha_n\ot \cdots\ot s\alpha_3\ot s(\alpha_2\alpha_1)\ot 1) = G_{|q|+1}^{n-1}(1\ot s\alpha_n\ot \cdots\ot s\alpha_2\ot \alpha_1).
   \] 
   
    \end{itemize}

    The remaining identities for $GF=\mathrm{id}$ can be deduced directly from the definitions of $F$ and $G$ and are therefore omitted.
\end{Proof}

\section{Hochschild cohomology of graded skew-gentle algebras}\label{Sec3}

In this section, as usual, we fix a graded skew-gentle triple $(Q,R,\Sp)$ with $Q$ a connected quiver, and set $A=\bbK Q/I$ the corresponding graded skew-gentle algebra, where $I$ denotes the two-sided ideal of $\bbK Q$ generated by $R\cup \{\varepsilon^2-\varepsilon\ |\ \varepsilon\in \Sp\}$. Throughout, we keep the notation $\ot:=\ot_E$ with $E:=\bbK Q_0$ as mentioned in Subsection \ref{section: TT calculus of dg algebra}. We also keep the notation  
$\Gamma_n$ already defined in Section \ref{section: cs-projective resolution}.

To state our main results, we need to extend some notation in \cite[Chapter 3]{CS2015} to the graded skew-gentle version. In the presentation $(Q,R,\Sp)$, a cycle $C=c_m\dots c_1$ is called \emph{primitive} if it is not a proper power of another cycle. Set 
\[
\mathrm{rot}(C):=c_{m-1} \cdots c_1 c_m,
\]
the cycle obtained from $C$ by `rotating it one step to the left'. In case $C=D^k$, where $D$ is a primitive cycle and $k$ is a positive integer, both uniquely determined by $C$, the length of the cycle $D$ will be called the \emph{period} of the cycle $C$.

We say that two cycles in $Q$ are \emph{conjugate} if one can be obtained from the other by repeated rotation. This is clearly an equivalence relation in the set of all cycles in the quiver, and we call its equivalence classes the \emph{circuits} of $Q$ and write $\C$ for the set of all circuits in $Q$. The \emph{length} (resp. \emph{degree}) of a circuit is the length (resp. degree) of any of the cycles it contains, and a circuit is \emph{primitive} if the cycles it contains are all primitive or, equivalently, if any one of them is.

A path $\alpha$ is \emph{$\calB$-maximal} (resp. \emph{$\Gamma$-maximal}) if it belongs to $\calB$ (resp. $\Gamma$) and is not a proper subpath of another element of $\calB$ (resp. $\Gamma$). A cycle $\alpha$ is \emph{cocomplete} (resp. $C$ is \emph{complete}) if $\alpha^2\in \calB$ (resp. $C^2\in \Gamma$). A circuit is \emph{cocomplete} (resp. \emph{complete}) if the cycles it contains are cocomplete (resp. \emph{complete}). We denote by
\begin{itemize}

    \item $\C (\calB)$ the set of all cocomplete circuits $\alpha$ satisfying the following condition:\\ 
    \vspace{1em}\emph{either the degree $|\alpha|$ of $\alpha$ is even or the characteristic $\Char(\bbK)$ of $\bbK$ is 2;}

     \item 
     $\C (\Gamma)$ the set of  all complete circuits $C$ in $\Gamma$ such that:\\ 
    \emph{either $l(C)-|C|$ is even or $\Char (\bbK)$ is 2.}
\end{itemize}

In addition, let $\oC(\calB)$ (resp. $\oC(\Gamma)$) be a set of representatives of the circuits that belong to $\C(\calB)$ (resp.  $\C(\Gamma)$) chosen so that

\emph{if $\alpha$ and $\beta$ are elements of $\oC(\calB)$ (resp. $\oC(\Gamma)$) and $\alpha\beta$ is a cocomplete (resp. complete) cycle in $(Q,I)$, then $\alpha\beta$ also belong to $\oC(\calB)$ (resp. $\oC(\Gamma)$).}

Moreover, for each positive integer $m$ we let $\C_m(\Gamma)$ (resp. $\oC_m(\Gamma)$) be the subset of $\C(\Gamma)$ (resp. $\oC(\Gamma)$) of circuits of length $m$.

A \emph{spanning tree} of the quiver $Q$ is the set of arrows of a connected acyclic subquiver of $Q$ that contains all its vertices.

\subsection{Hochschild cohomology spaces of graded gentle algebras}\

\medskip

The Hochschild cohomology spaces $\HH^m(A)$ ($m\ge 0$) of an ungraded gentle algebra $A$ were computed in \cite{CSSS24}. Now we extend these results to the graded gentle case.

In this section, let $(Q,R)$ be a graded gentle pair and $A=\bbK Q/I$ the corresponding graded gentle algebra, where $I$ is the ideal of $\bbK Q$ generated by $R$. We will view $A$ as a graded skew-gentle algebra with $\Sp=\emptyset$. Note that in this case $\Gamma_n$ is the set of those paths of length $n$ in $Q$ all of whose subpaths of length $2$ are in $R$ for each $n\geq 0$, and $\calB$ denotes the set of paths in $Q$ that do not have subpaths in $R$.

Our aim is to determine the cohomology of the complex $(\bbK(\Gamma\| \calB),d)$.
As already mentioned,  the grading of the quiver $Q$ induces a grading in this complex and the differentials are homogeneous, so the cohomology of the complex decomposes in each cohomological degree as a direct sum whose components are homogeneous with respect to the grading of the gentle algebra. We make use of the results about the Hochschild cohomology of gentle algebras given in \cite[Chapter 3]{CSSS24}, we consider  the complex obtained by forgetting the grading $|\cdot|$ and denote it by $(\bbK(\Gamma\| \calB),d_f)$, that is, for any $(q,\alpha)\in \Gamma_n\| \calB$, we have 
\[
(d^n-d_f^n)(q,\alpha)=\sum\limits_{\substack{b\in Q_1\\ bq\in \Gamma_{n+1}}} ((-1)^{|b|(|\alpha|-|q|)}-1)(bq,\pi(b\alpha)).
\]
We define the {\it weight} of the element $(q,\alpha)$ as $l(\alpha)-l(q)$. The \emph{rank} of a non zero element $u=\sum_{i=1}^mx_i(q_i,\alpha_i)\in\bbK(\Gamma_m\| \calB)$ is the cardinal of the support of the sum.

\begin{Prop}\label{prop: graded HH0} (cf.  \cite[Proposition 3.3]{CSSS24}) 
    The graded vector space $\HH^{0,\sbl}(A)$ is freely generated by the collection of the following elements of $\bbK (\Gamma_0\| \calB)$:
    \begin{itemize}
        \item the unit 
        \[
        \mathds{1}:=\sum_{i\in Q_0} (e_i,e_i)\in \HH^{0,0}(A),
        \]

        \item the pairs $(s(\alpha),\alpha)\in \HH^{0,|\alpha|}(A)$ with $\alpha$ a $\calB$-maximal path in  $(Q,R)$, and
        
        \item the sums 
        \[
         \llangle {\alpha}:=\sum_{i=0}^{r-1} (s(\rmrot^i(\alpha)),\rmrot^i(\alpha))\in \HH^{0,|\alpha|}(A),
        \]
        with $\alpha\in \oC(\calB)$ and $r$ its period.
        
    \end{itemize}
\end{Prop}

\begin{Proof}

    A straightforward calculation shows that the elements listed in the statement of the proposition lie in the kernel of the differential $d^0$, and it is clear that they are linearly independent. Therefore,  it suffices to show that every irreducible element of $\Ker(d^0)$ is a scalar multiple of one of the three listed elements.

    Assume that $u$ is an irreducible $0$-cocycle, then it is homogeneous and weight-homogeneous of some weight $l$. If $l=0$, then $|u|=0$ and  $d_f^0(u)=d^0(u)=0$, thus $u=\mathds{1}$. 
    
    Next, suppose that $l\ge 1$. We can write 
    \[
    u=\sum_{i=1}^rx_i(e_i, \alpha_i),\quad  r\ge 1
    \]
    with $(e_i,\alpha_i)\in \Gamma_0\| \calB$ and $0\neq x_i\in\bbK$ for all $1\le i\le r$. If there is an index $i\in\{1,\ldots,r\}$ such that $\alpha_i$ is $\calB$-maximal, then $d_f^0(e_i,\alpha_i)=d^0(e_i,\alpha_i)=0$, which implies $u=(e_1,\alpha_1)$.
    
    Now suppose that none of $\alpha_1,\ldots,\alpha_r$ appearing in $u$ is $\calB$-maximal, then there is a cocomplete cycle $\alpha=(a_r\cdots a_1)^{m}$ of period $r$ such that $\{\alpha_1,\ldots,\alpha_r\}=\{\alpha,\rmrot(\alpha),\ldots,\rmrot^{r-1}(\alpha)\}$ by the proof of \cite[Proposition 3.3]{CSSS24}. We can relabel the index and write $u=\sum\limits_{i=0}^{r-1}x_i(s(\rmrot^i(\alpha)),\rmrot^i(\alpha))$. Since $d^0(u)=0$, it must satisfy
    \[
    x_i=(-1)^{|a_{r-i+1}|\cdot|\alpha|}x_{i-1}, \qquad 1\le i\le r-1.
    \]
    Hence, $d^0(u)=(-1)^{|a_1|\cdot |\alpha|}(1-(-1)^{|\alpha|})x_0 (a_1,a_1\alpha)$, which vanishes if and only if $\Char(\bbK)=2$ or $|\alpha|$ is even. This completes the proof.
       
\end{Proof}

\begin{Prop}\label{prop:1st-graded-HH(gen)} (cf. \cite[Proposition 3.5]{CSSS24}) 
    Let  $T$ be a spanning tree for the graded quiver $Q$. The graded vector space $\HH^{1,\sbl}(A)$ is freely generated by the following elements of $\bbK(\Gamma\| \calB)$:

    \begin{itemize}
        \item the pairs $(c,c)\in \HH^{1,0}(A)$ with $c\in Q_1\setminus T$,

        \item the pairs $(c,\alpha)\in \HH^{1,|\alpha|-|c|}(A)$ with $c$ a $\Gamma$-maximal arrow and $\alpha$ a path that neither begins nor ends with $c$,

        \item the pairs $(c,c\alpha )\in \HH^{1,|\alpha|}(A)$ with $\alpha\in \oC(\calB)$ and $c$ the first arrow of $\alpha$,

        \item the pairs $(b,s(b))\in \HH^{1,-|b|}(A)$ with $b$ a loop in $Q$ such that $b^2\in R$, $|b|$ odd or $\Char(\bbK)=2$.
    \end{itemize}
\end{Prop}

\begin{Proof}
    Assume $u$ is an irreducible $1$-cocycle and not a coboundary in the  complex $\bbK(\Gamma\| \calB)$. From the proof of \cite[Proposition 3.5]{CSSS24}, the rank of $u$ is $1$, let $u=(c,\alpha)\in \Gamma_1\| \calB$. Since $d^1(u)=0$ one of the following three possibilities occurs:
    \begin{itemize}
        \item [(A)] either there are arrows $a$ and $b$ in $Q$ such that 
        \[
        (bc,b\alpha)=(ca,\alpha a)\in \Gamma_2\| \calB \quad \text{and}\quad (-1)^{(|\alpha|-|c|)|b|}+1=0,
        \]

        \item [(B)] or there is no arrow $b$ such that $(bc,b\alpha)\in \Gamma_2\| \calB$ and there is no arrow $a$ such that $(ca,\alpha a)\in \Gamma_2\| \calB$.
    \end{itemize}
    If $(A)$ holds, then $b=c=a$, $b$ is a loop, $b^2\in R$, $\alpha=s(b)$ and $(-1)^{|b|}+1=0$, that is the following condition holds:
    \begin{itemize}
        \item [(1)] there is a loop $b$ in $Q$ such that $b^2\in R$ and $u=(b,s(b))$, and $|b|$ odd or $\Char(\bbK)=2$.
    \end{itemize}

    Let us suppose that (B) holds. Then by the proof of \cite[Lemma 3.4]{CSSS24}, one of the following conditions holds:
    \begin{enumerate}
        \item [(2)] there is an arrow $c$ such that $u=(c,c)$,

        \item [(3)] there is an $\Gamma$-maximal arrow $c$ and a path $\alpha\in \calB$ that neither begins nor ends with $c$ such that $u=(c,\alpha)$,

        \item [(4)] there is a cocomplete cycle $\delta$ starting with an arrow $c$ such that $u=(c,c\delta)$.
    \end{enumerate}

    The elements satisfying the  conditions $(1), (2), (3)$ and $(4)$ are linearly independent in $\bbK(\Gamma\| \calB)$. Write $Z$ for the space they span. Let $z$ be a non zero weight-homogeneous coboundary in $Z$. 
    It is clear that the elements satisfying conditions (1) or (3) are not coboundaries. Once again by \cite[Proposition 3.5]{CSSS24}, there are two cases:
    \begin{itemize}
        \item $z$ is a linear combination of elements of type 
        \[
        d^0(s(\epsilon),\epsilon)=(-1)^{|b||\epsilon|}(b,b\epsilon)-(a, a\rmrot(\epsilon))
        \]
        with $\epsilon$ a cocomplete cycle, $b$ its first arrow and $a$ the first arrow of $\rmrot(\epsilon)$, which is the last arrow of $\epsilon$. In this case, we know from Proposition~\ref{prop: graded HH0} that $(b,b\epsilon)\in\rmIm(d^0)$ if and only if $\Char(\bbK)\neq 2$ and $|\epsilon|$ is odd;

        \item $z$ is a linear combination of elements of type 
        \[
        d^0(e_i,e_i)=\sum\limits_{\substack{b\in Q_1\\ s(b)=i}}(b,b)-\sum\limits_{\substack{a\in Q_1\\ t(a)=i}}(a,a)
        \]
        with $i\in Q_0$. From the proof of \cite[Theorem 3.4]{CSS20} we know that $z$ is cohomologous to a unique linear combination of pairs $(c,c)$ with $c\in Q_1\setminus T$.
    \end{itemize}
    This completes the proof of the proposition.

\end{Proof}

Let $C=(c_r\cdots c_1)^{\frac{m}{r}}$ be a complete cycle of length $m$ and period $r$ in $(Q,R)$. We set 
\[
\llangle{C}_{gr}:=\sum_{i=0}^{r-1}(-1)^{i m+(|c_r|+\cdots+|c_{r-i+1}|)|C|}(\rmrot^i(C),s(\rmrot^i(C))).
\]
\begin{Prop}\label{prop:higher-graded-HH(gen)} (cf. \cite[Proposition 3.5]{CSSS24})
    Let $m\ge 2$. The graded vector space $\HH^{m,\sbl}(A)$ is freely spanned by the collection of the following elements:
    \begin{itemize}
        \item $\llangle{C}_{gr}\in \HH^{m,-|C|}(A)$ with $C\in \oC_m(\Gamma)$,

        \item $(bC,b)\in \HH^{m,-|C|}(A)$ with $C\in \oC_{m-1}(\Gamma)$ and $b$ the first arrow of $C$,

        \item $(\gamma,\alpha)\in \HH^{m,|\alpha|-|\gamma|}(A)$ with $\gamma$ a $\Gamma$-maximal element of $\Gamma_m$ and $\gamma$ and $\alpha$ neither beginning nor ending with the same arrow.
    \end{itemize}
\end{Prop}

\begin{Proof}
    Let $u$ be an irreducible $m$-cocycle which is not a coboundary. If $r=\mathrm{rank}(u)>1$, we can write $u=\sum\limits_{i=1}^rx_i(\gamma_i,\alpha_i)$ with $(\gamma_1,\alpha_1),\ldots,(\gamma_r,\alpha_r)$ all in $\Gamma_m\| \calB$ and pairwise different, and scalars $x_1,\ldots,x_r\in\bbK$ all non zero. As proved in \cite[Proposition 3.9]{CSSS24}, $u$ satisfies the following three conditions:
    \begin{itemize}
        \item for all $i\in\{1,\ldots, r \}$ the path $\gamma_i$ is not a power of a loop,

        \item all the paths $\alpha_1,\ldots,\alpha_r $ have length $0$,

        \item for each $i\in \{1,\ldots,r\}$, the path $\gamma_i$ is a complete cycle.
    \end{itemize}
    More precisely, there is a complete cycle $C=(c_r \cdots c_1)^{\frac{m}{r}}$ of period $r$ such that 
    \[\{\gamma_1,\ldots,\gamma_r \}=\{C,\rmrot^1(C),\ldots,\rmrot^{r-1}(C)\}.
    \]
    We assume again  that $u=\sum\limits_{i=0}^{r-1} x_i(\rmrot^i(C),s(\rmrot^i(C)))$; in order for $u$ to be a cocycle, it must satisfy 
    \[
    x_i=(-1)^{|c_{r-i+1}|\cdot|C|+m}x_{i-1},\qquad 1\le i\le r-1.
    \]
    Under this requirement, we have $d^m(\llangle{C}_{gr})=(-1)^{|c_1||C|}(1-(-1)^{m-|C|})\,x_0 (c_1C,c_1)$. Hence, $u$ is a cocycle if and only if $m-|C|$ is even or $\Char(\bbK)=2$. That is, the following condition holds:
    \begin{itemize}
        \item [(a)] there is a complete cycle $C$ of period $r\ge 2$ such that $m-|C|$ is even or $\Char(\bbK)=2$, and $u=\llangle{C}_{gr}$.
    \end{itemize}

    If $\mathrm{rank}(u)=1$, let $u=(\gamma,\alpha)\in \Gamma_m\| \calB$. In view of the definition of $d^m$, this tells us that 
    \begin{itemize}
        \item [(A)] either there are arrows $a$ and $b$ in $Q$ such that 
        \[
        (b\gamma,b\alpha)=(\gamma a,\alpha a)\in\Gamma_{m+1}\| \calB\quad \text{and}\quad (-1)^{(|\alpha|-|\gamma|)|b|}-(-1)^m=0,
        \]

        \item [(B)] or there is no arrow $b$ such that $(b\gamma,b\alpha)\in \Gamma_{m+1}\| \calB$ and there is no arrow $a$ such that $(\gamma a,\alpha a)\in \Gamma_{m+1}\| \calB$.
    \end{itemize}
    As obtained in the proof of \cite[Proposition 3.10]{CSSS24}, the first condition in (A) implies that $b$ is a loop  such that $\gamma=b^m$ and $\alpha=s(b)$, while the second condition tells us $m-|\gamma|$ is even or $\Char(\bbK)=2$. That is, the following condition holds:
    \begin{itemize}
        \item [(a')] there is a loop $b$ such that $m-m|b|$ is even or $\Char(\bbK)=2$, and $u=(b^m,s(b))$.
    \end{itemize}

    Let us suppose that we are in case (B). From the proof of \cite[Proposition 3.10]{CSSS24}, $u$ is one of the following elements:
    \begin{itemize}

        \item [(b)] $(bC,b)$ where $C$ is a complete cycle in $\Gamma_{m-1}$ with $m-|C|$ odd or $\Char(\bbK)=2$, and $b$ the first arrow of $C$;

        \item [(c)] $(\gamma,\alpha)$ where $\gamma$ is a $\Gamma$-maximal element of $\Gamma_m$ and $\gamma$ and $\alpha$ neither begin nor end with the same arrow.
    \end{itemize}

    The elements in (a), (b) and (c) are linearly independent in $\bbK(\Gamma\| \calB)$, let us write $Z$ for the subspace they span. By \cite[Proposition 3.11]{CSSS24} and the previous analysis, the subspace of coboundaries in $Z$ is spanned by the elements: 
    \begin{itemize}
        \item $(bC,b)$, with $C$ a complete circuit in $\Gamma_{m-1}$, $m-|C|$ even, $\Char(\bbK)\neq 2$, and $b$ the first arrow of $C$, and 

        \item $(-1)^{|b|\cdot|C|} (bC,b)-(-1)^{m-1}(a\rmrot(C),a)$, with $C$ a complete circuit in $\Gamma_{m-1}$, $m-|C|$ even, $\Char(\bbK)\neq 2$ and $b$ and $a$ the first arrows of $C$ and of $\rmrot(C)$, respectively. 
    \end{itemize}
    This proves the proposition.
    
\end{Proof}

By combining the three preceding propositions, we obtain the following theorem.

\begin{Thm}\label{basis-of-HH(gentle)-as-vs}
Let $A=\bbK Q/I$ be a graded gentle algebra, and let $T$ be a spanning tree for the graded quiver $Q$. The Hochschild cohomology $\HH^{*}(A)=\bigoplus\limits_{*=\bullet+\sbl}\HH^{\bullet, \sbl}(A)$ of $A$ is the graded vector space freely generated by the cohomology classes of the following homogeneous cocycles of the complex $\bbK (\Gamma\|  \calB)$:
    \begin{itemize}
        \item [$(\mathrm{G_I})$] the element 
        \[
        \mathds{1}:=\sum_{i\in Q_0} (e_i,e_i)\in \HH^{0,0}(A),
        \]

        \item [$(\mathrm{G_{II}})$] the pairs $(s(\alpha),\alpha)\in \HH^{0,|\alpha|}(A)$, with $\alpha$ a $\calB$-maximal path in  $(Q,I)$,

        \item [$(\mathrm{G_{III}})$] the sum 
        \[
        \llangle {\alpha}:=\sum_{i=0}^{r-1} (s(\rmrot^i(\alpha)),\rmrot^i(\alpha))\in \HH^{0,|\alpha|}(A),
        \]
        with $\alpha\in \oC(\calB)$and $r$ its period,

        \item [$(\mathrm{G_{IV}})$] the pairs $(c,c)\in \HH^{1,0}(A)$, with $c$ an arrow in the complement of the spanning tree $T$,

        \item [$(\mathrm{G_{V}})$] the pair $(c,c\alpha )\in \HH^{1,|\alpha|}(A)$, with $\alpha\in \oC(\calB)$ and $c$ the first arrow of $\alpha$,

        \item [$(\mathrm{G_{VI}})$] the pairs $(\gamma,\alpha)\in \HH^{l(\gamma),|\alpha|-|\gamma|}(A)$, with $\gamma$ a $\Gamma$-maximal element of $\Gamma$ and $\gamma$ and $\alpha$ neither beginning nor ending with the same arrow,

        \item [$(\mathrm{G_{VII}})$] the sums
        \[
       \llangle{C}_{gr}:=\sum_{i=0}^{r-1}(-1)^{i m+(|c_r|+\cdots+|c_{r-i+1}|)|C|}(\rmrot^i(C),s(\rmrot^i(C)))\in \HH^{m,-|C|}(A),
        \]
        where $C\in \oC_m(\Gamma)$ with period $r$ and $C=(c_r\cdots c_1)^{\frac{m}{r}}$,

        \item [$(\mathrm{G_{VIII}})$] the pairs $(bC,b)\in \HH^{m,-|C|}(A)$, with $C\in \oC_{m-1}(\Gamma)$ and $b$ the first arrow of $C$.
    \end{itemize}
\end{Thm}

\subsection{Gentle vs. skew-gentle}\label{section: gentle vs skew-gentle}\ 

\medskip

One can associate a graded gentle algebra $A'\coloneqq \bbK Q/\langle R\cup \{\varepsilon^2 |\ \varepsilon\in \Sp\} \rangle$ to each graded skew-gentle algebra $A:=\bbK Q/\langle R\cup \{\varepsilon^2-\varepsilon\ \mid \varepsilon\in \Sp\} \rangle$ as follows: 
\[A\coloneqq \bbK Q/\langle R\cup \{\varepsilon^2-\varepsilon\ |\ \varepsilon\in \Sp\} \rangle\xrightarrow{\text{replace $\varepsilon^2=\varepsilon $ by $\varepsilon^2=0$}}A':=\bbK Q/\langle R\cup \{\varepsilon^2\mid \varepsilon\in \Sp\} \rangle.\]
In fact, $A'$ is a degeneration of $A$ in the sense of \cite{GP95}, see also \cite{C19}, and we can also view $A$ as a deformation of $A'$, see \cite{BSW25} for details. Moreover, their Hochschild cohomologies are also related to each other.

Once we view $A'$ as a graded skew-gentle algebra without special loops, Proposition \ref{CS-cochain-complex} can be applied to $A'$. Note that in this case, the sets $\Gamma$ and $\calB$ associated to $A$, respectively coincide with the sets $\Gamma'$ and $\calB'$ associated to $A'$. We shall denote by $(\bbK(\Gamma\|\calB),d)$ and $(\bbK(\Gamma\|\calB),\delta)$ the complexes in Proposition \ref{CS-cochain-complex} for $A$ and $A'$ respectively. Recall that $\Sp^n$ denotes the set $\{\varepsilon^n\ |\ \varepsilon\in \Sp\}$ for each $n\geq 1$. Define $d^n_{sp}:=d^n|_{\bbK(\Sp^n\|\calB)}$ and $d^n_{or}:=d^n|_{\bbK((\Gamma_n\setminus\Sp^n)\|\calB)}$ as the restrictions of $d^n$ for $n\geq 1$, and similarly for $\delta^n$. 

\begin{Lem}\label{lem1-vs}
Let $(\bbK(\Gamma\|\calB),d)$ and $(\bbK(\Gamma\|\calB),\delta)$ be the complexes of Proposition \ref{CS-cochain-complex} for $A$ and $A'$ respectively. For each $n\geq 1$ we have 
\begin{itemize}
    \item $d^n=d^n_{sp}\oplus d^n_{or}$ and $\delta^n=\delta^n_{sp}\oplus \delta^n_{or}$.
    In particular, \[\begin{cases}
    \Ker(d^n) = \Ker(d^n_{sp}) \oplus \Ker(d^n_{or}) \\
    \mathrm{Im}(d^n) = \mathrm{Im}(d^n_{sp}) \oplus \mathrm{Im}(d^n_{or})
    \end{cases} \text{and} \quad \begin{cases}
    \Ker(\delta^n) = \Ker(\delta^n_{sp}) \oplus \Ker(\delta^n_{or}) \\
    \mathrm{Im}(\delta^n) = \mathrm{Im}(\delta^n_{sp}) \oplus \mathrm{Im}(\delta^n_{or})
    \end{cases}. \]
    \item $d^n_{or}=\delta^n_{or}$. Consequently, $\Ker(d^n_{or})=\Ker(\delta^n_{or})$ and $\mathrm{Im}(d^n_{or})=\mathrm{Im}(\delta^n_{or})$.
\end{itemize}
\end{Lem}
\begin{Proof}
We shall show that $d^n=d^n_{sp}\oplus d^n_{or}$, from which the first statement for $d^n$ will follow. Indeed, it suffices to prove that the summands of elements in $\mathrm{Im}(d^n_{sp})$ cannot be cancelled by the summands of elements in $\mathrm{Im}(d^n_{or})$.

By the formula of the differential in Proposition \ref{CS-cochain-complex}, for $(\varepsilon^n,\alpha)\in \Gamma_n\|\calB$ and $(q,\alpha')\in (\Gamma_n\setminus\Sp^n)\|\calB$, \[d^n_{sp}(\varepsilon^n,\alpha)=(\varepsilon^{n+1},\pi_A(\varepsilon\alpha))+(-1)^{n+1}(\varepsilon^{n+1},\pi_A(\alpha\varepsilon))-\delta_{(-1)^{n+1},1}\cdot (\varepsilon^{n+1},\alpha),\] \[d^n_{or}(q,\alpha')=\sum\limits_{\substack{b\in Q_1\\ bq\in \Gamma_{n+1}}} (-1)^{|b|(|\alpha'|-|q|)}(bq,\pi_A(b\alpha'))+(-1)^{n+1} \sum\limits_{\substack{a\in Q_1\\ q a\in \Gamma_{n+1}}} (q a,\pi_A(\alpha' a)), \] where $\delta_{i,j}$ denotes the Kronecker symbol. If there exists some $a,b\in Q_1$ such that $bq,qa\in \Gamma_{n+1}$, then such $bq$ and $qa$ cannot be $\varepsilon^{n+1}$ for $\varepsilon\in \Sp$ since $q\in \Gamma_n\setminus\Sp^n$. Therefore, the summands of $d^n_{sp}(\varepsilon^n,\alpha)$ cannot be cancelled by the summands of $d^n_{or}(q,\alpha')$. 

A similar proof holds for $\delta^n$ just by replacing $\pi_A$ by $\pi_{A'}$ in $d^n_{or}(q,\alpha')$ to get $\delta^n_{or}(q,\alpha')$, also note that \(\delta^n_{sp}(\varepsilon^n,\alpha)= (\varepsilon^{n+1},\pi_{A'}(\varepsilon\alpha))+(-1)^{n+1}(\varepsilon^{n+1},\pi_{A'}(\alpha\varepsilon))\). Thus, the first statement follows. For $q\in \Gamma_n\setminus \Sp^n$, if there exists some $a,b\in Q_1$ such that $bq,qa\in \Gamma_{n+1}$, then $b,a\notin \Sp$. Hence, $\pi_A(b\alpha')=\pi_{A'}(b\alpha')$ and $\pi_A(\alpha' a)=\pi_{A'}(\alpha' a)$. It follows that $d^n_{or}(q,\alpha')=\delta^n_{or}(q,\alpha')$ for any $(q,\alpha')\in (\Gamma_n\setminus \Sp^n)\| \calB$.
\end{Proof}

For any path $\alpha$ in $Q$ of positive length, denote the last and the first arrow of $\alpha$ by $\lfact{1}{\alpha}$ and $\rfact{2}{\alpha}$, respectively. Say $\alpha=a_l\cdots a_1\in Q_l$ of length $l\geq 1$, then $\lfact{1}{\alpha}=a_l$ and $\rfact{2}{\alpha}=a_1$. There are also two associated paths $\lfact{2}{\alpha}$ and $\rfact{1}{\alpha}$, possibly of length zero, such that 
\[\alpha=\lfact{1}{\alpha}\lfact{2}{\alpha}=\rfact{1}{\alpha}\rfact{2}{\alpha}.\]

\begin{Rmk}\label{compare-differentials}
The relation between the maps $\delta^n_{sp}, d^n_{sp}: \bbK(\Sp^n\|\calB)  \to \bbK(\Sp^{n+1}\|\calB)$ is as follows.
Given $(\varepsilon^n,\alpha)\in \bbK(\Sp^n\|\calB)$ with $n\geq 1$,
\begin{itemize}
    \item if $\rfact{2}{\alpha}= \varepsilon$ but $\lfact{1}{\alpha}\neq \varepsilon$, then \(d^n_{sp}(\varepsilon^n,\alpha)=\begin{cases}
        \delta^n_{sp}(\varepsilon^n,\alpha) & \text{if $n$ is odd,}\\
        \delta^n_{sp}(\varepsilon^n,\alpha)-(\varepsilon^{n+1},\alpha) & \text{if $n$ is even,}
    \end{cases}\)

    \item if $\rfact{2}{\alpha}\neq \varepsilon$ but $\lfact{1}{\alpha}= \varepsilon$, then \(d^n_{sp}(\varepsilon^n,\alpha)=\begin{cases}
        \delta^n_{sp}(\varepsilon^n,\alpha) & \text{if $n$ is odd,}\\
        -\delta^n_{sp}(\varepsilon^n,\alpha)+(\varepsilon^{n+1},\alpha) & \text{if $n$ is even,}
    \end{cases}\)

    \item if $\rfact{2}{\alpha}= \varepsilon$ and $\lfact{1}{\alpha}= \varepsilon$, then $\delta^n_{sp}(\varepsilon^n,\alpha)=0$ and \(d^n_{sp}(\varepsilon^n,\alpha)=\begin{cases}
        (\varepsilon^{n+1},\alpha) & \text{if $n$ is odd,}\\
        0 & \text{if $n$ is even,}
    \end{cases}\)

    \item in all other cases, we have \(d^n_{sp}(\varepsilon^n,\alpha)=\begin{cases}
        \delta^n_{sp}(\varepsilon^n,\alpha)-(\varepsilon^{n+1},\alpha) & \text{if $n$ is odd,}\\
        \delta^n_{sp}(\varepsilon^n,\alpha) & \text{if $n$ is even.}
    \end{cases}\)
\end{itemize}
\end{Rmk}

Write $V_{sp}^1:=\begin{cases}
        \langle(\varepsilon,s(\varepsilon)), (\varepsilon,\varepsilon)\mid \varepsilon\in \Sp\rangle & \text{if $\Char(\bbK)=2,$}\\ \langle (\varepsilon,\varepsilon)\mid \varepsilon\in \Sp\rangle & \text{if $\Char(\bbK)\neq2,$}
    \end{cases}$ and define $V_{sp}^n:=\Ker(\delta^n_{sp})/\mathrm{Im}(\delta^{n-1}_{sp})$ for $n\geq 2$. A direct calculation shows that for all $n\geq 1$, \[V_{sp}^n\cong \begin{cases}
        \langle (\varepsilon^n,s(\varepsilon)), (\varepsilon^n,\varepsilon)\mid \varepsilon\in\Sp\rangle & \text{if $\Char(\bbK)=2,$}\\ \langle (\varepsilon^n,s(\varepsilon)) \mid \varepsilon\in \Sp\rangle & \text{if $\Char(\bbK)\neq 2$ and $n$ is even,} \\ \langle (\varepsilon^n,\varepsilon) \mid \varepsilon\in \Sp\rangle & \text{if $\Char(\bbK)\neq 2$ and $n$ is odd.}
    \end{cases}\]

The proof of the next lemma is straightforward after Remark \ref{compare-differentials}.
\begin{Lem}\label{lem2-vs}
    There are $\bbK$-linear isomorphisms:
    \begin{itemize}
        \item for $n\geq 2$, $\Ker(d^n_{sp})/\mathrm{Im}(d^{n-1}_{sp})\cong 0$,

        \item for $n\geq 1$,  $\Ker(\delta^n_{sp})\cong \begin{cases}
            \Ker(d^n_{sp}) & \text{if $n$ is even,}\\
            \Ker(d^n_{sp})\oplus V_{sp}^n & \text{if $n$ is odd,}
        \end{cases}$

        \item  for $n\geq 1$,  $\mathrm{Im}(d^n_{sp})\cong \begin{cases}
            \mathrm{Im}(\delta^n_{sp}) & \text{if $n$ is even,}\\
            \mathrm{Im}(\delta^n_{sp})\oplus V_{sp}^{n+1} & \text{if $n$ is odd.}
        \end{cases}$
    \end{itemize}
\end{Lem}

In the next lemma we take care of the case $n=0$.

\begin{Lem}\label{lem3-vs}
$\mathrm{Im}(\delta^0)\cong\mathrm{Im}(d^0)$ and $\Ker(\delta^0)\cong\Ker(d^0)$.
\end{Lem}
\begin{Proof}
    By definition of $\delta^0$ and $d^0$ there is a bijection between $\mathrm{Im}(\delta^0)$ and $\mathrm{Im}(d^0)$ given by \begin{align*}
        \theta: \mathrm{Im}(\delta^0)&\to \mathrm{Im}(d^0)\\
        \delta^0(e,p)&\mapsto d^0(e,p)=\begin{cases}
            \delta^0(e,p)-(\varepsilon,p) & \text{if $(e,p)\in \Sp\|\calB_{\geq 2}$ with $\rfact{2}{p}=\varepsilon$ but $\lfact{1}{p}\neq\varepsilon$,}\\
            \delta^0(e,p)+(\varepsilon,p) & \text{if $(e,p)\in \Sp\|\calB_{\geq 2}$ with $\rfact{2}{p}\neq\varepsilon$ but $\lfact{1}{p}=\varepsilon$,}\\
            \delta^0(e,p) & \text{otherwise.}
        \end{cases}
    \end{align*} 
    Therefore, we obtain $\mathrm{Im}(\delta^0)\cong\mathrm{Im}(d^0)$ as $\bbK$-vector spaces. Also note that $\bbK(Q_0\|\calB)/\mathrm{Im}(\delta^0) \cong \Ker(\delta^0)$ and $\bbK(Q_0\|\calB)/\mathrm{Im}(d^0) \cong \Ker(d^0)$, thus $\Ker(\delta^0)\cong\Ker(d^0)$.
\end{Proof}

\smallskip

We are now ready to state the main theorem relating the graded Hochschild cohomology of the gentle algebra with the graded Hochschild cohomology of its skew-gentle  deformation. We already know from \cite{GP95}, that in the ungraded case, the dimensions of the Hochschild cohomology spaces of the deformation are smaller than or equal to the dimensions of the respective Hochschild cohomology spaces of the original algebra. In the next theorem we prove that in the case of  graded gentle algebras and their graded skew-gentle deformations, not only the same happens but we explicitly obtain generators of the spaces providing this difference of dimensions: the $\bbK$-vector spaces $V_{sp}^n$, which are completely determined by the special loops in the skew-gentle algebras.

\begin{Thm}\label{gentle-vs-sg}
    Let $A'$ be the graded gentle algebra associated to the graded skew-gentle algebra $A$. We have the following isomorphisms of $\bbK$-vector spaces:
    \begin{itemize}
        \item $\HH^{n,j}(A')\cong\begin{cases}
            \HH^{n,j}(A)\bigoplus V_{sp}^n & \text{if $n\geq 1$ and $j=0$,}\\ \HH^{n,j}(A) & \text{if $n=0$ or $j\neq 0$.}
        \end{cases}$
        \item $\HH^N(A')\cong\begin{cases}
            \HH^N(A) & \text{if $N\leq 0$,}\\ \HH^N(A)\bigoplus V_{sp}^N & \text{if $N>0$.}
        \end{cases}$
        \item $\HH^*(A')\cong \HH^*(A)\bigoplus(\bigoplus\limits_{N>0}V_{sp}^N)$.
    \end{itemize}
\end{Thm}

\begin{Proof}
    Note that the first statement implies the second, which in turn implies the third. Indeed, 
\[
\HH^0(A')=\bigoplus\limits_{\substack{n+j=0\\n\geq 0, j\in \bbZ}}\HH^{n,j}(A')=\bigoplus\limits_{n\geq 0}\HH^{n,-n}(A')=\HH^{0,0}(A')\bigoplus \big(\bigoplus\limits_{n\geq 1}\HH^{n,-n}(A')\big) \cong \]
\[  \HH^{0,0}(A)\bigoplus \big(\bigoplus\limits_{n\geq 1}\HH^{n,-n}(A)\big)=\bigoplus\limits_{\substack{n+j=0\\n\geq 0, j\in \bbZ}}\HH^{n,j}(A)=\HH^0(A).
\] 
    
    If $N\neq 0$, then 
    \begin{align*}
        \HH^N(A')&=\bigoplus\limits_{\substack{n+j=N\\n\geq 0, j\in \bbZ}}\HH^{n,j}(A')=\bigoplus\limits_{n\geq 0}\HH^{n,N-n}(A')=\HH^{0,N}(A')\bigoplus \big( \bigoplus\limits_{n\geq 1}\HH^{n,N-n}(A') \big)\\ & =
        \HH^{0,N}(A') \bigoplus \begin{cases}
            \HH^{N,0}(A')\bigoplus \big( \bigoplus\limits_{\substack{N\neq n\\n\geq 1}}\HH^{n,N-n}(A') \big) & \text{if $N>0$,}\\ \bigoplus\limits_{\substack{N\neq n\\n\geq 1}}\HH^{n,N-n}(A') & \text{if $N<0$}.
        \end{cases}\\
        & \cong \HH^{0,N}(A)\bigoplus \begin{cases}
            \bigoplus\limits_{n\geq 1}\HH^{n,N-n}(A)\bigoplus V_{sp}^N & \text{if $N>0$,} \\ \bigoplus\limits_{n\geq 1}\HH^{n,N-n}(A) & \text{if $N<0$.} 
        \end{cases}=\begin{cases}
            \HH^N(A)\bigoplus V_{sp}^N & \text{if $N>0$,}\\ \HH^N(A) & \text{if $N<0$.}
        \end{cases}
    \end{align*}
    
    Now we prove the first statement. Clearly, to compare $\HH^{n,j}(A')$ with $\HH^{n,j}(A)$,
    it is enough to compare $\Ker(\delta^n)$ with $\Ker(d^n)$, and $\mathrm{Im}(\delta^n)$ with $\mathrm{Im}(d^n)$. By Lemma \ref{lem1-vs} and Lemma \ref{lem2-vs}, we have 
    \begin{align*}
        \frac{\Ker(d^n)}{\mathrm{Im}(d^{n-1})} & \cong \frac{\Ker(d^n_{sp})}{\mathrm{Im}(d^{n-1}_{sp})}\oplus\frac{\Ker(d^n_{or})}{\mathrm{Im}(d^{n-1}_{or})}\cong \frac{\Ker(d^n_{or})}{\mathrm{Im}(d^{n-1}_{or})}\cong \frac{\Ker(\delta^n_{or})}{\mathrm{Im}(\delta^{n-1}_{or})}\quad \text{and}\\ \frac{\Ker(\delta^n)}{\mathrm{Im}(\delta^{n-1})} & \cong \frac{\Ker(\delta^n_{sp})}{\mathrm{Im}(\delta^{n-1}_{sp})}\oplus\frac{\Ker(\delta^n_{or})}{\mathrm{Im}(\delta^{n-1}_{or})}\cong V_{sp}^{n}\oplus\frac{\Ker(\delta^n_{or})}{\mathrm{Im}(\delta^{n-1}_{or})}\quad \text{for each $n\geq 2$},\\
        \frac{\Ker(\delta^1)}{\mathrm{Im}(\delta^0)} & \cong \frac{\Ker(d^1_{or})\oplus(\Ker(d^1_{sp})\oplus V_{sp}^1)}{\mathrm{Im}(\delta^0)}\cong \frac{\Ker(d^1)\oplus V_{sp}^1}{\mathrm{Im}(d^0)}\cong \frac{\Ker(d^1)}{\mathrm{Im}(d^0)}\oplus V_{sp}^1.
    \end{align*}
    Note that $\HH^{n,j}(A'):=(\frac{\Ker(\delta^n)}{\mathrm{Im}(\delta^{n-1})})^j$ and $\HH^{n,j}(A):=(\frac{\Ker(d^n)}{\mathrm{Im}(d^{n-1})})^j$ for all $n\geq 1$ and $j\in \bbZ$. It follows that $\HH^{n,j}(A')\cong \HH^{n,j}(A)\oplus V_{sp}^{n,j}$ for $n\geq 1$, where $V_{sp}^{n,j}$ be the subspace of $V_{sp}^n$ consisting of elements of internal degree $j$. But all elements in $V_{sp}^n$ are of internal degree $0$, so $\HH^{n,0}(A')\cong \HH^{n,0}(A)\oplus V_{sp}^n$ and $\HH^{n,j}(A')\cong \HH^{n,j}(A)$ for $n\geq 1$ and $j\neq 0$. The isomorphism  $\HH^{0,j}(A')\cong \HH^{0,j}(A)$ is a consequence of Lemma \ref{lem3-vs}. We are done.
    
\end{Proof}

\subsection{Hochschild cohomology spaces of graded skew-gentle algebras}\

\medskip

We set, like in previous sections, a graded skew-gentle triple $(Q,R,\Sp)$ and set $A=\bbK Q/I$ the corresponding graded skew-gentle algebra, where $I$ denotes the two-sided ideal of $\bbK Q$ generated by $R\cup \{\varepsilon^2-\varepsilon\ |\ \varepsilon\in \Sp\}$. Let $T$ be a spanning tree of $Q$. Throughout, we follow the notation in Subsection \ref{section: gentle vs skew-gentle} and associate a graded gentle algebra $A':=\bbK Q/\langle R\cup \{\varepsilon^2\mid \varepsilon\in \Sp\} \rangle$ to $A$. 

Now, we will view $A'$ as a graded skew-gentle algebra without special loops. In other words, the special loops in $\Sp$ for $A$ become  ordinary loops  in $A'$. We keep the notation $\Gamma$ 
defined in Section \ref{section: cs-projective resolution}. Note that in this case, $A$ and $A'$ have the same quiver $Q=(Q_0,Q_1)$ but different relation sets $S:=R\cup \{\varepsilon^2-\varepsilon\mid \varepsilon\in \Sp\}$ and $S':=R\cup \{\varepsilon^2\mid \varepsilon\in \Sp\}$, respectively. Hence, the set $\Gamma'$ for $A'$ consisting of those paths in $Q$ whose subpaths of length $2$ are in $S'$ is exactly the set $\Gamma$ for $A$, and the basis $\calB'$ of $A'$ is also the same as the basis $\calB$ of $A$. Let $I'$ be the ideal of $\bbK Q$ generated by $S'$, denote $\pi=\pi_A:\bbK Q\to A=\bbK Q/I$ and $\pi_{A'}:\bbK Q\to A'=\bbK Q/I'$ by the canonical projections for $A$ and $A'$ respectively.

Before stating the main results, we need one more definition. For any cycle $\alpha=a_m\cdots a_1\in Q_m$, we define the {\it special rotation} of $\alpha$, denoted by $\rmsrot(\alpha)$, as follows:
\[
\rmsrot(\alpha):=
\begin{cases}
    a_{m-1}\cdots a_1a_m=\rmrot(\alpha) & \text{if $a_m\notin \Sp$,} \\ a_{m-1}\cdots a_1(\varepsilon-s(\varepsilon))=\rmrot(\alpha)-a_{m-1}\cdots a_1 & \text{if $a_m=\varepsilon \in \Sp$.}
\end{cases}
\]
For $i\in \bbZ_+$, define $\rmsrot^i(\alpha):=\rmsrot(\rmsrot^{i-1}(\alpha))$ inductively with the convention $\rmsrot^0(\alpha)=\alpha$. We can also extend this definition to the  elements  in the basis $\calB$ of $A$ by linearity. 

\begin{Rmk}\label{property-of-srot}
   Let $A$ be a (graded) skew-gentle algebra, and let $\alpha$ be a cocomplete cycle in the basis $\calB$. The following facts easily follow by induction:
   \[ \rmsrot^i(\alpha)=\rmsrot(\rmrot^{i-1}(\alpha)), \hspace{3ex} \rmsrot^i(\alpha^j)=(\rmsrot^i(\alpha))^j, \hbox{ for any } i,j\geq 1. 
\]
\end{Rmk}

Given a cocomplete cycle $\alpha\in\oC(\calB)$ with period $r$, for each $0\leq i\leq r-1$ write  \[\rmrot^i(\alpha)=\rfact{1}{(\rmrot^i(\alpha))}\rfact{2}{(\rmrot^i(\alpha))}=\lfact{1}{(\rmrot^i(\alpha))}\lfact{2}{(\rmrot^i(\alpha))}.\] Recall that $\llangle{\alpha}:=\sum\limits_{i=0}^{r-1} (s(\rmrot^i(\alpha)),\rmrot^i(\alpha))$. Similarly, denote by \[\llangle {\alpha}_s:=\sum_{i=1}^{r} (s(\rmsrot^i(\alpha)),\rmsrot^i(\alpha)).\] 
Note that if $\alpha$ is a cocomplete cycle starting from a special loop, i.e., $\alpha'_\rt \in \Sp$, then $\sum\limits_{i=1}^{r} (s(\rmsrot^i(\alpha)),\rmsrot^i(\alpha)) \neq \sum\limits_{i=0}^{r-1} (s(\rmsrot^i(\alpha)),\rmsrot^i(\alpha))$. In fact, we have \[\llangle{\alpha}_s=\sum\limits_{i=1}^{r} (s(\rmsrot^i(\alpha)),\rmsrot^i(\alpha))=\sum\limits_{i=0}^{r-1} (s(\rmsrot^i(\alpha)),\rmsrot^i(\alpha))-(s(\alpha),\alpha_\rt).\]

The following lemma is helpful to state the main result in this subsection.

\begin{Lem}\label{lem-rot-vs-srot}
    Let $A$ be a (graded) skew-gentle algebra, and $\alpha\in\oC(\calB)$ with period $r$. Then 
    \begin{enumerate}
       \item[$(1)$] $\llangle{\alpha}-\llangle{\alpha}_s=\sum\limits_{\substack{i=1\\ \rfact{2}{(\rmrot^i(\alpha))}\in\Sp}}^{r} \Big(s\big(\rmrot^i(\alpha)\big),\rfact{1}{\big(\rmrot^i(\alpha)\big)}\Big)=\sum\limits_{\substack{i=1\\ \lfact{1}{(\rmrot^i(\alpha))}\in\Sp}}^{r} \Big(s\big(\rmrot^i(\alpha)\big),\lfact{2}{\big(\rmrot^i(\alpha)\big)}\Big)$,

       \item[$(2)$] $d^0(\llangle{\alpha})=-\delta^0(\llangle{\alpha}_s)$,

       \item[$(3)$] $d^0(\llangle{\alpha}_s)=0$.
   \end{enumerate}
\end{Lem}
\begin{Proof}
    The first equality is straightforward, the second one follows from the fact that 
    $\lfact{1}{(\rmrot^i(\alpha))}=\rfact{2}{(\rmrot^{i+1}(\alpha))}$ and $\lfact{2}{(\rmrot^i(\alpha))}=\rfact{1}{(\rmrot^{i+1}(\alpha))}$ for each $1\leq i\leq r$. The second item can be derived from the formula in the proof of Lemma \ref{lem3-vs} and the third statement of Proposition \ref{prop: graded HH0}. 
    More precisely, 
    \begin{align*}
       d^0(\llangle{\alpha}) & = \sum_{i=1}^{r} d^0\big(s(\rmrot^i(\alpha)),\rmrot^i(\alpha)\big)\\
        & =\sum_{i=1}^{r}\begin{cases}
           \delta^0\big(s(\rmrot^i(\alpha)),\rmrot^i(\alpha)\big)+(\lfact{1}{(\rmrot^i(\alpha))},\rmrot^i(\alpha)) & \text{if $\lfact{1}{(\rmrot^i(\alpha))}\in \Sp$,}\\
            \delta^0\big(s(\rmrot^i(\alpha)),\rmrot^i(\alpha)\big)-(\rfact{2}{(\rmrot^i(\alpha))},\rmrot^i(\alpha)) & \text{if $\rfact{2}{(\rmrot^i(\alpha))}\in \Sp$,}\\
            \delta^0\big(s(\rmrot^i(\alpha)),\rmrot^i(\alpha)\big) & \text{otherwise.}
        \end{cases}\\ &=\delta^0(\llangle{\alpha})+ \sum_{\substack{i=1\\\lfact{1}{(\rmrot^i(\alpha))}\in\Sp}}^{r}(\lfact{1}{(\rmrot^i(\alpha))},\rmrot^i(\alpha))-\sum_{\substack{i=1\\\rfact{2}{(\rmrot^i(\alpha))}\in\Sp}}^{r}(\rfact{2}{(\rmrot^i(\alpha))},\rmrot^i(\alpha))\\ &=\sum_{\substack{i=1\\\lfact{1}{(\rmrot^i(\alpha))}\in\Sp}}^{r}(\lfact{1}{(\rmrot^i(\alpha))},\rmrot^i(\alpha))-\sum_{\substack{i=1\\\rfact{2}{(\rmrot^i(\alpha))}\in\Sp}}^{r}(\rfact{2}{(\rmrot^i(\alpha))},\rmrot^i(\alpha)).
    \end{align*} 
    On the other hand, combining the first statement with $\delta^0(\llangle{\alpha})=0$, we obtain that \[-\delta^0(\llangle{\alpha}_s)=\delta^0(\llangle{\alpha}-\llangle{\alpha}_s)=\sum\limits_{\substack{i=1\\ \lfact{1}{(\rmrot^i(\alpha))}\in\Sp}}^{r} \delta^0\Big(s\big(\rmrot^i(\alpha)\big),\lfact{2}{\big(\rmrot^i(\alpha)\big)}\Big).\] If there exists some $i$ with  $1\leq i\leq r$ such that $\lfact{1}{(\rmrot^i(\alpha))}\in \Sp$, then $\rmrot^i(\alpha)\in \oC(\calB)$ implies that $\lfact{2}{(\rmrot^i(\alpha))}$ neither begins nor ends with a special loop. 
    In fact, $\lfact{2}{(\rmrot^i(\alpha))}$ starts from $\lfact{1}{(\rmrot^{i-1}(\alpha))}$ and ends at $\lfact{1}{(\rmrot^{i+1}(\alpha))}$ with $\lfact{1}{(\rmrot^{i-1}(\alpha))}\cdot\lfact{1}{(\rmrot^{i+1}(\alpha))}\in R$. As a consequence, the definition of $\delta^0$ implies that $\delta^0\Big(s\big(\rmrot^i(\alpha)\big),\lfact{2}{\big(\rmrot^i(\alpha)\big)}\Big)=\big(\lfact{1}{(\rmrot^i(\alpha))},\rmrot^i(\alpha)\big)-\big(\rfact{2}{(\rmrot^{i+1}(\alpha))},\rmrot^{i+1}(\alpha)\big)$. 
    Therefore, we obtain that $d^0(\llangle{\alpha})=-\delta^0(\llangle{\alpha}_s)$, which further implies that $d^0(\llangle{\alpha}_s)=0$ since  $d^0(\llangle{\alpha})=-\delta^0(\llangle{\alpha}_s)=\delta^0(\llangle{\alpha}-\llangle{\alpha}_s)=d^0(\llangle{\alpha}-\llangle{\alpha}_s)$.
\end{Proof}

Roughly speaking, combining Theorem \ref{basis-of-HH(gentle)-as-vs} with Theorem \ref{gentle-vs-sg}, we can recover the generators of $\HH^*(A)$ from $\HH^*(A')$ by killing the elements in $V_{sp}^n$ given by special loops of $A$. However, what we want is to get the specific generators of $\HH^*(A)$.
To achieve this goal, we need to know exactly how the generators of $\HH^*(A')$ transfer to those of $\HH^*(A)$. Indeed, bringing into account the equalities $\HH^{\bullet,\sbl}(A)=(\frac{\Ker(d^\bullet)}{\mathrm{Im(d^{\bullet-1})}})^{\sbl}$ and $\HH^{\bullet,\sbl}(A')=(\frac{\Ker(\delta^\bullet)}{\mathrm{Im(\delta^{\bullet-1})}})^{\sbl}$ in mind, we can further reduce this problem to describing  
$\frac{\Ker(d^\bullet)}{\mathrm{Im(d^{\bullet-1})}}$ and $\frac{\Ker(\delta^\bullet)}{\mathrm{Im(\delta^{\bullet-1})}}$ since the internal degree of a specific element in $\HH^{\bullet,\sbl}$ is clear.

Let us start from the lower degree Hochschild cohomology  spaces $\HH^{\bullet,\sbl}$.
\begin{Prop}\label{0th-graded-HH(sg)}
The graded vector space $\HH^{0,\sbl}(A)$ is freely generated by the collection of the following elements of $\bbK (\Gamma_0\| \calB)$:
    \begin{itemize}
        \item the element 
        \[
        \mathds{1}:=\sum_{i\in Q_0} (e_i,e_i)\in \HH^{0,0}(A),
        \]

        \item the pairs $(s(\alpha),\alpha)\in \HH^{0,|\alpha|}(A)$ with $\alpha$ a $\calB$-maximal path in  $(Q,I)$, and
        
        \item the sums 
        \[
        \llangle {\alpha}_s:=\sum_{i=1}^{r} (s(\rmsrot^i(\alpha)),\rmsrot^i(\alpha))\in \HH^{0,|\alpha|}(A)
        \]
        with $\alpha\in \oC(\calB)$ and $r$ its period.
        
    \end{itemize}    
\end{Prop}
\begin{Proof}
    Combining Proposition \ref{prop: graded HH0} with the first  statement of Theorem \ref{gentle-vs-sg}, it suffices to show that the elements listed above belong to $\HH^{0,\sbl}(A)$. Given an element in $\HH^{\bullet,\sbl}(A)$, we know its internal degree. Thus, we can reduce further to prove that these elements belong to $\HH^{0,\sbl}(A)=\Ker(d^0)$. 
    
    With the formula of $d^0$ and $\delta^0$ in the proof of Lemma \ref{lem3-vs} in mind, it is clear that $d^0(e_i,e_i)=\delta^0(e_i,e_i)$ for all $i\in Q_0$. Note that a path $\alpha$ is $\calB'$-maximal in $(Q,I')$ if and only if $\alpha$ is $\calB$-maximal in $(Q,I)$, and if it is a cycle of length at least $2$ then $\alpha$ cannot start or end at some $\varepsilon\in\Sp$. Thus, we also have $d^0(s(\alpha),\alpha)=\delta^0(s(\alpha),\alpha)=0$. Consequently, the first and second classes of elements are exactly in $\HH^{0,\sbl}(A)$. For any $\alpha\in\oC(\calB)$, $\llangle{\alpha}_s\in \Ker(d^0)$ the result follows immediately from Lemma \ref{lem-rot-vs-srot}. 
\end{Proof}

Next we turn our attention to the first Hochschild cohomology $\HH^{1,\sbl}(A)$.

\begin{Prop}\label{1th-graded-HH(sg)}
    The graded vector space $\HH^{1,\sbl}(A)$ is freely generated by:

    \begin{itemize}
        \item the pairs $(c,c)\in \HH^{1,0}(A)$ with $c\in Q_1\setminus (T\cup \Sp)$,

        \item the pairs $(c,\alpha)\in \HH^{1,|\alpha|-|c|}(A)$ with $c$ an $\Gamma$-maximal arrow and $\alpha$ a path that neither begins nor ends with $c$,

        \item the pair $(c,c\alpha )\in \HH^{1,|\alpha|}(A)$ with $\alpha\in \oC(\calB)$ and $c\notin \Sp$ the first arrow of $\alpha$,

        \item the pair $(b,s(b))\in \HH^{1,-|b|}(A)$ with $b$ a loop in $Q$ such that $b^2\in R$, $|b|$ odd or $\Char(\bbK)=2$.
    \end{itemize}
\end{Prop}
\begin{Proof}
    By Theorem \ref{gentle-vs-sg}, we know that $\HH^{1,\sbl}(A')\cong \HH^{1,\sbl}(A)\oplus V_{sp}^1$. Observe that there is a canonical inclusion from $V_{sp}^1$ to $\HH^{1,\sbl}(A')$. Indeed, it is clear that $V_{sp}^1\subseteq \Ker(\delta_{sp}^1)\subseteq \Ker(\delta^1)$. Moreover, $\mathrm{Im}(\delta^0)\subseteq \bbK(Q_1\|\calB_{\geq 1})$ implies that $(\varepsilon,s(\varepsilon))\notin \mathrm{Im}(\delta^0)$ for any $\varepsilon\in \Sp$. In addition, $(\varepsilon,\varepsilon)\notin \mathrm{Im}(\delta^0)$ can be deduced from the fact that $\varepsilon$ is a loop in $Q_1$. Then the generators in $\HH^{1,\sbl}(A')$ excluding those in $V_{sp}^1$ would be in $1-1$-correspondence with the generators in $\HH^{1,\sbl}(A)$. 
    
    We claim that it is enough to show the four types of elements listed above lie in $\HH^{1,\sbl}(A)=\Ker(d^1)/\mathrm{Im}(d^0)$. 
    Indeed, those elements giving rise to $V_{sp}^1$ can only appear in the first and the fourth classes of elements in Proposition \ref{prop:1st-graded-HH(gen)}. Moreover, for any $\alpha\in \oC(\calB)$ having $c$ as the first arrow,    if $c\in\Sp$ then $(c,c\alpha)\in \HH^{1,\sbl}(A')$ but it no longer belongs to $\HH^{1,\sbl}(A)$. This follows from the fact that in this case $(c,c\alpha)\in \Sp\|\calB_{\geq 2}$, where $c\alpha$ is a path such that its starting arrow and its ending arrow are both equal to the same special loop $c\in\Sp$. Then, the third statement of Remark \ref{compare-differentials} implies that $d^1(c,c\alpha)=d^1_{sp}(c,c\alpha)=(c^2,c\alpha)\neq 0$. In other words, $(c,c\alpha)$ is not a cocycle for $A$.

    Let $(a,\alpha)$ be in any of the above four classes of elements, then $a\in Q_1\setminus\Sp$. Clearly, $(a,\alpha)\in\Ker(d^1)$ since $d^1(a,\alpha)=d_{or}^1(a,\alpha)=\delta_{or}^1(a,\alpha)=0$ by Lemma \ref{lem1-vs} and Proposition \ref{prop:1st-graded-HH(gen)}. It remains to show $(a,\alpha)$ is not a coboundary for $A$. Actually, this holds for the second and the fourth classes of elements since the definition of $d^0$ yields to $\mathrm{Im}(d^0)\subseteq\langle (a,a),(a,a\alpha),(a,\alpha a)\mid a\in Q_1,\alpha\in \calB_{\geq 1} \rangle$.

    To show $(a,\alpha)\notin \mathrm{Im}(d^0)$ for the remaining two classes of elements, we proceed as follows. If $(c,c)\in \mathrm{Im}(d^0)$ for any $c\in Q_1\setminus(T\cup\Sp)$, then $(c,c)\in \mathrm{Im}(d^0|_{\bbK(Q_0\|Q_0)})=\mathrm{Im}(\delta^0|_{\bbK(Q_0\|Q_0)})$ since $d^0(e,e)=\delta^0(e,e)$ for any $e\in Q_0$ by Lemma \ref{lem3-vs}. Thus, we have $(c,c)\in \mathrm{Im}(\delta^0)$, which contradicts Proposition \ref{prop:1st-graded-HH(gen)}. For the third class of elements $(c,c\alpha)$ with $\alpha\in \oC(\calB)$, we claim that $(c,c\alpha)\in \mathrm{Im}(d^0)$ implies $(c,c\alpha)\in \mathrm{Im}(\delta^0)$, which is impossible by Proposition \ref{prop:1st-graded-HH(gen)}. In fact, if $(c,c\alpha)\in \mathrm{Im}(d^0)$, then $(c,c\alpha)\in \mathrm{Im}(d^0|_{\bbK(Q_0\|\calB_{\geq 1})})$. Suppose that $(c,c\alpha)=\sum \lambda_{\beta}\cdot d^0(s(\beta),\beta)$,
    where $\beta$ can only be $\rmrot^i(\alpha)$ for $0\leq i\leq r-1$. Note that all such $\lambda_{\beta}\in\bbK$ should be equal. In other words, we can write $(c,c\alpha)=\lambda\cdot\sum\limits_{i=0}^{r-1} d^0\Big(s\big(\rmrot^i(\alpha)\big),\rmrot^i(\alpha)\Big)=\lambda\cdot d^0(\llangle{\alpha})$ with $\lambda\in\bbK$. As a consequence, $(c,c\alpha)=-\lambda\cdot\delta^0(\llangle{\alpha}_s)\in \mathrm{Im}(\delta^0)$ by Lemma \ref{lem-rot-vs-srot}. Therefore,  no $(a,\alpha)$ listed in this proposition can be coboundary.
     \end{Proof}

Now we describe the higher Hochschild cohomology spaces.
\begin{Prop}\label{higher-graded-HH(sg)}
    Let $m\ge 2$. The graded vector space $\HH^{m,\sbl}(A)$ is freely spanned by the  following elements:
    \begin{itemize}
        \item $\llangle{C}_{gr}\in \HH^{m,-|C|}(A)$ with $C\in \oC_m(\Gamma)\setminus \Sp^m$,

        \item $(bC,b)\in \HH^{m,-|C|}(A)$ with $C\in \oC_{m-1}(\Gamma)\setminus \Sp^{m-1}$ and $b$ the first arrow of $C$,

        \item $(\gamma,\alpha)\in \HH^{m,|\alpha|-|\gamma|}(A)$ with $\gamma$ an $\Gamma$-maximal element of $\Gamma_m$ and $\gamma$ and $\alpha$ neither beginning nor ending with the same arrow.
    \end{itemize}
\end{Prop}
\begin{Proof}
    Let $m\geq 2$. It follows from Theorem \ref{gentle-vs-sg} that $\HH^{m,\sbl}(A')\cong \HH^{m,\sbl}(A)\oplus V_{sp}^m$. Indeed, once we apply Proposition \ref{prop:higher-graded-HH(gen)} to $\HH^{m,\sbl}(A')$,  the definition $V_{sp}^m:=\Ker(\delta_{sp}^n)/\mathrm{Im}(\delta_{sp}^{n-1})$ shows that there is a canonical inclusion of $V_{sp}^m$ into $\HH^{m,\sbl}(A')$. Clearly, $V_{sp}^m$ can come only from the first two classes of elements in Proposition \ref{prop:higher-graded-HH(gen)}. Also, observe that the set $\Gamma'$ corresponding to  $A'$ is the same as $\Gamma$ of $A$, hence $C\in \oC(\Gamma')$ if and only if $C\in \oC(\Gamma)$ and $\gamma$ is $\Gamma'$-maximal if and only if $\gamma$ is $\Gamma$-maximal. It suffices to show that the elements listed above are exactly in $\HH^{m,\sbl}(A)=\Ker(d^m)/\rmIm(d^{m-1})$. 

    Let $x$ be any of the three elements listed in the statement. It is easy to see that $x\in \bbK((\Gamma_m\setminus\Sp^m)\|\calB)$. Combining the definition of $d^m$ with  Lemma \ref{lem1-vs}, we have $d^m(x)=d^m_{or}(x)=\delta^m_{or}(x)=\delta^m(x)$, which is equal to zero by Proposition \ref{prop:higher-graded-HH(gen)}. Thus, $x\in \Ker(d^m)$. Note that $\rmIm(d^{m-1}_{sp})\subseteq \bbK(\Sp^m\|\calB)$, therefore, if $x\in\rmIm(d^{m-1})$ then $x\in\rmIm(d^{m-1}_{or})$. Moreover, Lemma \ref{lem1-vs} shows that $x\in\rmIm(d^{m-1}_{or})=\rmIm(\delta^{m-1}_{or})\subseteq \rmIm(\delta^{m-1})$, contradicting  Proposition \ref{prop:higher-graded-HH(gen)}. Therefore, we have $x\notin\rmIm(d^{m-1})$. 
\end{Proof}

We are now able to prove the main theorem of this subsection.
\begin{Thm}\label{Thm: coho-basis-sg}
    Let $A$ be the graded skew-gentle algebra corresponding to a graded skew-gentle triple $(Q,R,\Sp)$. The graded Hochschild cohomology $\HH^*(A)=\bigoplus\limits_{*=\bullet+\sbl}\HH^{\bullet, \sbl}(A)$ of $A$ is the graded vector space freely generated by the cohomology classes of the following homogeneous cocycles of the complex $\bbK (\Gamma\| \calB)$:
    \begin{itemize}
        \item [$(\mathrm{H_I})$] the element,
        \[
        \mathds{1}:=\sum_{i\in Q_0} (e_i,e_i)\in \HH^{0,0}(A),
        \]

        \item [$(\mathrm{H_{II}})$] the pairs $(s(\alpha),\alpha)\in \HH^{0,|\alpha|}(A)$, with $\alpha$ a $\calB$-maximal path in  $(Q,I)$,

        \item [$(\mathrm{H_{III}})$] the sums 
        \[
        \llangle {\alpha}_s:=\sum_{i=1}^{r} (s(\rmsrot^i(\alpha)),\rmsrot^i(\alpha))\in \HH^{0,|\alpha|}(A),
        \]
        with $\alpha\in \oC(\calB)$ and $r$ its period,

        \item [$(\mathrm{H_{IV}})$] the pairs $(c,c)\in \HH^{1,0}(A)$, with  $c\in Q_1\setminus(T\cup \Sp)$,

        \item [$(\mathrm{H_{V}})$] the pairs $(c, c\alpha)\in \HH^{1,|\alpha|}(A)$, with $\alpha\in \oC(\calB)$ and $c$ the first arrow of $\alpha$ and $c\notin \Sp$,

        \item [$(\mathrm{H_{VI}})$] the pairs $(\gamma,\alpha)\in \HH^{l(\gamma),|\alpha|-|\gamma|}(A)$, with $\gamma$ a $\Gamma$-maximal element of $\Gamma$ and $\gamma$ and $\alpha$ neither beginning nor ending with the same arrow,

        \item [$(\mathrm{H_{VII}})$]  the sums
        \[ 
        \llangle{C}_{gr}:=\sum_{i=0}^{r-1}(-1)^{im+(|c_{r-i+1}|+\cdots+|c_r|)|C|}(\rmrot^i(C),s(\rmrot^i(C)))\in \HH^{m,-|C|}(A),
        \]
        where $C\in \oC_m(\Gamma)\setminus\Sp^m$ with period $r$ and $C=(c_r\cdots c_1)^{\frac{m}{r}}$,

        \item [$(\mathrm{H_{VIII}})$] the pairs $(bC,b)\in \HH^{m,-|C|}(A)$, with $C\in \oC_{m-1}(\Gamma)\setminus\Sp^{m-1}$ and $b$ the first arrow of $C$.
    \end{itemize}
\end{Thm}

\begin{Proof}
The proof follows from Propositions \ref{0th-graded-HH(sg)}, \ref{1th-graded-HH(sg)} and \ref{higher-graded-HH(sg)}.
\end{Proof}

\subsection{The cup product}\

\medskip

Let $(Q,R,\Sp)$ be a graded skew-gentle triple and let $A$ be the the graded skew-gentle algebra  associated to this data. 

The cup product from Section~\ref{section: TT calculus of dg algebra} is defined in terms of the complex $C^*(A)$. Using the comparison morphisms $F$ and $G$ - see Section~\ref{section: comparison-map}, we can describe the cup product on the complex $\bbK(\Gamma\| \calB)$. 

\begin{Lem}

     Let $A$ be the graded skew-gentle algebra associated to a graded skew-gentle triple $(Q,R,\Sp)$. If $m, n \geq 0$, given  $(\gamma,\alpha)\in \bbK(\Gamma_m\|  \calB)$ and $(\gamma',\alpha')\in \bbK(\Gamma_n\|  \calB)$, then
     \begin{align*}
       \smile\ :\ \bbK(\Gamma_m\| \calB) \times \bbK(\Gamma_n\| \calB)
        & \longrightarrow \bbK(\Gamma_{m+n}\| \calB) \\
        \big(\ (\gamma,\alpha)\;,\; (\gamma',\alpha')\ \big)\; & \longmapsto 
        \begin{cases}
            (-1)^{(|\alpha'|-|\gamma'|)(m+|\gamma|)+mn}\ (\gamma\gamma',\pi(\alpha\alpha')) & \text{if } \gamma\gamma'\in \Gamma_{m+n}, \\  
            0 & \text{otherwise.}
        \end{cases}
\end{align*}

\end{Lem}

The proof is direct. 

Next we recall the following lemma from \cite{CSSS24}.

\begin{Lem}\label{lemma: one cycle quiver}( \cite[Lemma 5.1]{CSSS24}). 
    Let $\gamma$ be a $\Gamma$-complete path in $\Gamma$ and $\alpha$ a $\calB$-maximal cycle in $\calB$, and suppose that $s(\gamma)=s(\alpha)$. If either $\gamma$ or $\alpha$ has length $1$, then the quiver has exactly one vertex and one loop.
\end{Lem}

Let us suppose that the quiver is not a loop and go through the entries of Table~1. We will deal with that exceptional case later.
\begin{table}[H]
\centering
{\renewcommand{\arraystretch}{1.5}
\begin{tabular}{|c|c|c|c|c|c|c|c|}
\hline 
$\smile$ & $(s(\alpha),\alpha)$ & $\llangle{\alpha}_s$ & $(c,c\alpha)$ & $(c,c)$ & $(\gamma,\alpha)$ & $\llangle{C}_{gr}$ & $(b C,b)$ \\ \hline
 $(s(\alpha),\alpha)$ &\cellcolor{pink!30} $0_{\mathrm{max}}$ & \cellcolor{pink!30} $0_{\mathrm{max}}$ & \cellcolor{pink!30} $0_{\mathrm{max}}$ & \cellcolor{pink!30} $0_{\mathrm{max}}$ & $0$ (v) &  $0$ (vi.a)&\cellcolor{pink!30} $0_{\mathrm{max}}$\\ \hline
 $\llangle{\alpha}_s$ & \cellcolor{gray!30}  &\cellcolor{green!10} (iv.b)  &\cellcolor{green!10} (iv.b) &\cellcolor{green!10} (iv.a) & \cellcolor{pink!30} $0_{\mathrm{max}}$ & $0$ (vi.b)& $0$ (vi.b) \\ \hline
$(c,c\alpha)$ & \cellcolor{gray!30} & \cellcolor{gray!30} & \cellcolor{blue!15} $0_{\mathrm{glue}}$ & \cellcolor{blue!15} $0_{\mathrm{glue}}$ & \cellcolor{pink!30} $0_{\mathrm{max}}$ &$0$ (vi.c) & \cellcolor{blue!15} $0_{\mathrm{glue}}$\\ \hline
  $(c,c)$ & \cellcolor{gray!30} & \cellcolor{gray!30} &\cellcolor{gray!30} & \cellcolor{blue!15} $0_{\mathrm{glue}}$ & \cellcolor{pink!30} $0_{\mathrm{max}}$ &\cellcolor{green!10} (vi.d) &\cellcolor{blue!15} $0_{\mathrm{glue}}$  \\ \hline
 $(\gamma,\alpha)$ & \cellcolor{gray!30} & \cellcolor{gray!30} &\cellcolor{gray!30} & \cellcolor{gray!30} & \cellcolor{pink!30} $0_{\mathrm{max}}$ & \cellcolor{pink!30} $0_{\mathrm{max}}$& \cellcolor{pink!30} $0_{\mathrm{max}}$\\ \hline
$\llangle{C}_{gr}$ &\cellcolor{gray!30} &\cellcolor{gray!30} &\cellcolor{gray!30} &\cellcolor{gray!30} &\cellcolor{gray!30} &\cellcolor{green!10} (vi.e) &\cellcolor{green!10} (vi.e) \\ \hline
$(bC,b)$ &\cellcolor{gray!30} &\cellcolor{gray!30} &\cellcolor{gray!30} &\cellcolor{gray!30} &\cellcolor{gray!30} &\cellcolor{gray!30} & \cellcolor{blue!15} $0_{\mathrm{glue}}$ \\ \hline
\end{tabular}}
\caption{The cup product of elements of our basis of $\HH^*(A)$, for $Q$ not a quiver with one vertex and one loop.}
\label{table: cup product}
\end{table}

The cup product essentially coincides with the gentle case, see \cite[Table~1]{CSSS24}. 
Besides the proof of (vi.b), the remaining proofs follow along the same lines. For the convenience of the reader, we state the results here, but only prove (vi.b). The numbering here starts from (iv), which may seem somewhat unusual. This choice is made in order to keep a complete correspondence with the gentle case.
\begin{itemize}
    \item [(iv.a)] $\llangle {\alpha}_s \smile (c,c)=
    \begin{cases}
        (d,d\alpha) &\text{if the cycle $\alpha$ goes through the arrow $c$, where $d$ is the first arrow of $\alpha$,} \\
        0 & \text{otherwise,}
    \end{cases}$    \item [(iv.b)] $\llangle{\alpha}_s\smile \llangle{\delta}_s=
        \begin{cases}
        \llangle{\alpha\delta}_s     &\text{if $\alpha$ and $\delta$ are powers of same primitive cycle,}  \\
         0    & \text{otherwise,}
        \end{cases}$ \\
        $\llangle{\alpha}_s\smile (c,c\delta)=
        \begin{cases}
        (c,c\alpha\delta)     &\text{if $\alpha$ and $\delta$ are powers of same primitive cycle,}   \\
         0    & \text{otherwise,}
        \end{cases}$
    \item [(vi.d)] $(c,c)\smile \llangle{C}_{gr}=
    \begin{cases}
        (-1)^{|d|\cdot |C|}(d C,d) &\text{if $C$ goes through the arrow $c$, where $d$ is the first arrow of $C$,}\\
        0 &\text{otherwise,}
    \end{cases}$
    \item [(vi.e)] $\llangle{C}_{gr}\smile \llangle{D}_{gr}=
        \begin{cases}
        (-1)^{l(C)\cdot l(D)}\llangle{CD}_{gr}     &\text{if $C$ and $D$ are powers of same primitive cycle,}  \\
         0    & \text{otherwise,}
        \end{cases}$\\
        $(bC,b)\smile \llangle{D}_{gr}=
        \begin{cases}
        (-1)^{(|b|+l(C))\cdot l(D)} (bCD,b)     &\text{if $C$ and $D$ are powers of same primitive cycle,}   \\
         0,    & \text{otherwise.}
        \end{cases}$
\end{itemize}


We start by observing that for $\alpha\in \C(\calB)$, $\rmsrot^i(\alpha)=\rmrot^i(\alpha)$ if the first arrow of $\rmrot^i(\alpha)$ is not a special loop. 

Now we prove (vi.b).

\begin{itemize}
    \item [(vi.b)] 
    Let $C=(C_r\cdots C_1)^{\frac{m}{r}}\in \oC(\Gamma)\setminus \Sp^m$ with length $m$ and period $r$.
     Let $\alpha\in\oC(\calB)$ and $r'$ be its period. We define two sets 
        \[
        I:=\{s(\rmrot^i(C))\ |\ i=0,\ldots, r-1\},\quad  J:= \{s(\rmsrot^j(\alpha))\ |\ j=1,\ldots,r'\},
        \]
        and $K=I\cap J$. If $K=\emptyset$, i.e., $s(\rmrot^i(C))\neq s(\rmsrot^j(\alpha))$ for all $i$ and $j$, then we clearly have $\llangle{\alpha}_s\smile \llangle{C}_{gr}=0$. Let us suppose that, on the contrary, $K\neq\emptyset$, then for any $e\in K$, there are integers $i_e\in\{ 0,\ldots, r-1 \}$ and $j_e\in\{ 1,\ldots,r' \}$ such that $s(\rmrot^i(C))= s(\rmsrot^j(\alpha))=e$, let $I_e$ and $J_e$ be the sets whose elements are such $i_e$ and $j_e$, respectively. We can write $\llangle{\alpha}_s \smile \llangle{C}_{gr}= \sum\limits_{e\in K} x_e$ with 
        \[
        x_e:= \sum_{i\in I_e,j\in J_e} (-1)^{im+(|C_{r-i+1}|+\cdots |C_r|)|C|}(\rmrot^i(C),\rmsrot^j(\alpha)).
        \]
        The set $K$ can be written as a disjoint union of connected components $K=K_1\cup \cdots \cup K_n$.  
        For any connected component $K_t$ containing no special point, the definition of skew-gentle algebra implies that $K_t$ consists of exactly two points $\{ e_1,e_2\}$, and the arrow $a:e_1\to e_2$ is the same arrow $\alpha$ which must be contained in $C$. Then 
       \[
        \sum\limits_{e\in K_t}x_e=x_{e_1}+x_{e_2}=\pm \big( (\rmrot^{i_{e_2}}(C),\rmsrot^{j_{e_2}}(\alpha))+(-1)^{m+|a|\cdot |C|} (\rmrot^{i_{e_2+1}}(C),\rmsrot^{j_{e_2+1}}(\alpha)) \big)
        \]
        is a coboundary, so $\sum\limits_{e\in K_t}x_e\equiv 0$. If $K_t$ contains special points, assume that the cardinality of $K_t$ is $s$, i.e., $K_t=\{e_1,\ldots,e_s\}$. Then $s\ge 3$, and among these $s$ points, $s-2$ are special points, while the  remaining two are not, that is, $Q$ contains the following subquiver
        
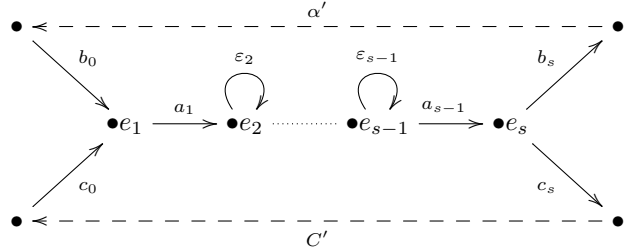
\[
        \xymatrix{ \bullet\ar[dr]^{b_0}&&&&& \bullet \ar@{-->}[lllll]_{\alpha'} \\
        & \bullet e_1\ar[r]^{a_1} & \bullet e_2 \ar@(ul,ur)[]^{\varepsilon_2}\ar@{.}[r]^{}&\bullet e_{s-1}\ar[r]^{a_{s-1}}\ar@(ul,ur)[]^{\varepsilon_{s-1}} &\bullet e_s\ar[ur]^{b_s}\ar[dr]_{c_s}\\
        \bullet\ar[ur]_{c_0}&&&&& \bullet\ar@{-->}[lllll]^{C'}}
        \]
        with $c_sa_{s-1},\, a_{s-1} a_{s-2},\,\ldots,\,  a_2 a_1,\, a_1 c_0\in R$  and $\varepsilon_2,\ldots, \varepsilon_{s-1}$ special loops. Here, $a_1,\ldots, a_{s-1}$ denote the arrows through which both $\alpha$ and $C$ pass, and $\alpha'$, $C'$ are the paths with length greater than or equal to $0$ such that 
       \[
        \rmrot^{i_{e_1}}(C)=c_0 C' c_s a_{s-1}\cdots a_2 a_1, \text{ and } \rmsrot^{j_{e_1}}(\alpha)=b_0 \alpha' b_s a_{s-1} \varepsilon_{s-1}\cdots a_2\varepsilon_2 a_1.
        \]
        Then 
        \begin{align*}
            \sum_{e\in K_t}x_e 
            = &x_{e_1}+\cdots x_{e_s}\\
            = &\sum_{i=2}^{s-1} \pm \big( (a_i\cdots a_1 c_0 C'c_s a_{s-1}\cdots a_{i+1},a_i\varepsilon_i\cdots a_2\varepsilon_2a_1b_0\alpha'b_s  a_{i+1}\varepsilon_{i+1})\\
            &\quad \quad \quad+(-1)^{m+|a_i|\cdot |C|}(a_{i-1}\cdots a_1 c_0 C'c_s a_{s-1}\cdots a_i,\varepsilon_i\cdots a_2\varepsilon_2a_1b_0\alpha'b_s  a_{i+1}\varepsilon_{i+1}a_i)\big)\\
            &\pm \big( (a_1c_0 C'c_s a_{s-1}\cdots a_2, a_1 b_0\alpha' b_s a_{s-1}\varepsilon_{s-1}\cdots a_2 \varepsilon_2)\\
            &\quad \quad +(-1)^{m+|a_1|\cdot |C|} (c_0 C'c_s a_{s-1}\cdots a_1,  b_0\alpha' b_s a_{s-1}\varepsilon_{s-1}\cdots a_2 \varepsilon_2 a_1)\big)\\
            &+\sum_{j=1}^{s-2}\pm(a_j\cdots a_1 c_0 C' c_s a_{s-1}\cdots a_{j+1},a_j\varepsilon_j\cdots a_2\varepsilon_2a_1 b_0\alpha'a_{s-1}\varepsilon_{s-1}\cdots a_{j+1})\\
            \equiv &0,
        \end{align*}
        so that the classes $\llangle{\alpha}_s\smile \llangle{C}_{gr}$ are zero in cohomology.

        Now we calculate $\llangle{\alpha}_s\smile (bC,b)$ with $b$ the first arrow in $C$. If $I_{t(b)}=\emptyset$, then we have $\llangle{\alpha}_s\smile (bC,b)=0$. We assume that $I_{t(b)}\neq \emptyset$, then
        \[
        \llangle{\alpha}_s \smile (bC,b)=\sum_{j\in J_{t(b)}} (bC,\pi(\rmsrot^j(\alpha)b))=\sum_{j\in J_{t(b)}}(\rmrot^{r-1}(C)b,\pi(\rmsrot^j(\alpha)b)).
        \]
        If $t(b)$ is not a special vertex, then $I_{t(b)}=\{j_{t(b)}\}$. If the first arrow of $\rmsrot^{j_{t(b)}}(\alpha)$ coincides with the first arrow of $\rmrot^{r-1}(C)$, then $\pi(\rmsrot^j(\alpha)b)=0$, that is,
        \[
        \llangle{\alpha}_s \smile (bC,b)=(\rmrot^{r-1}(C)b,\pi(\rmsrot^j(\alpha)b))=0.
        \]
        On the other hand, if the first arrow of $\rmsrot^{j_{t(b)}}(\alpha)$ is different from that of $\rmrot^{r-1}(C)$, then $b$ is the last arrow of $\rmsrot^{j_{t(b)}}(\alpha)$. In this case, $(\rmrot^{r-1}(C)b,\rmsrot^j(\alpha)b)$ is a coboundary, and hence $\llangle{\alpha}_s \smile (bC,b)\equiv 0$.

        If $t(b)$ is a special vertex, let $\varepsilon$ denote the special loop at $t(b)$, then we can set $J_{t(b)}=\{j_{t(b)},j_{t(b)}+1\}$. Suppose that $\alpha'$ is the cycle such that $\rmsrot^{j_{t(b)}}(\alpha)=\varepsilon \alpha'$. In this case, we have $\rmsrot^{j_{t(b)}+1}(\alpha)=\alpha'\varepsilon-\alpha'$, and the first arrow of $\alpha'$ coincides with the first arrow of $\rmrot^{r-1}(C)$, while the last arrow of $\alpha'$ is $b$. Therefore,
        \[
         \llangle{\alpha}_s \smile (bC,b)=(\rmrot^{r-1}(C)b,\pi(\varepsilon \alpha' b+\alpha' \varepsilon b-\alpha' b))=(\rmrot^{r-1}(C)b,\alpha' \varepsilon b),
        \]
        which is a coboundary.
\end{itemize}

Now we consider the cases previously excluded.

\begin{Rmk}\label{remark: cup product of one cycle quiver}
    Let us suppose that the quiver $Q$ has exactly one vertex and one loop $a$, so that the spanning tree $T$ is empty 
    Depending on $\Char(\bbK)$, the degree $|a|$ and the relations, we have five cases:
    \begin{itemize}
         \item If $a\in \Sp$, then:
        \[
        \HH^k(A)=
        \left\{
        \begin{array}{ll}
        \left\langle \mathds{1},(s(a),a) \right\rangle     & \text{if } k=0, \\
        0     & \text{if } k\ge 1. 
        \end{array}
        \right.
        \]
        Also, $(s(a),a)\smile (s(a),a)=(s(a),a)$.
     \end{itemize}

     Now, we assume that $a\notin \Sp$ and consider the following four cases.
     \begin{itemize}
        \item If $a^2=0$, $\Char(\bbK)= 2$ or $|a|$ is odd, then $\oC(\Gamma)=\{a^l\ |\ l\ge 1\}$, the cohomology vector space is as follows:
        \[
        \HH^k(A)=
        \left\{
        \begin{array}{ll}
        \langle\mathds{1},(s(a),a) \rangle     & \text{if } k=0, \\
        \langle(a^k,s(a)),(a^k,a)\rangle     & \text{if } k\ge 1. 
        \end{array}
        \right.
        \]
        For all $m,n\ge 1$, we have 
        \[
        \begin{array}{lll}
        &(s(a),a)\smile (s(a),a)=0,  &\quad (s(a),a)\smile (a^n,s(a))=(a^n,a) \\
        &(s(a),a)\smile (a^n,a)=0, &  \quad (a^m,s(a))\smile (a^n,s(a))=(-1)^{mn}(a^{m+n},s(a)) \\
        &(a^m,s(a))\smile (a^n,a)=(-1)^{mn}(a^{m+n},a), &\quad (a^m,a)\smile (a^n,a)=0.
        \end{array}
        \]

        \item If $a^2=0$, $\Char(\bbK)\neq 2$ and $|a|$ is even, we get:
        \[
        \HH^k(A)=
        \left\{
        \begin{array}{ll}
        \langle \mathds{1},(s(a),a) \rangle     & \text{if } k=0, \\
        \langle(a^k,s(a))\rangle    & \text{if $k\ge 1$ is even},\\
        \langle(a^k,a)\rangle     & \text{if $k\ge 1$ is odd}. 
        \end{array}
        \right.
        \]
         For all $m,n\ge 1$, we have 
        \[
        \begin{array}{lll}
        &(s(a),a)\smile (s(a),a)=0,  &\quad (s(a),a)\smile (a^{2n},s(a))=0 \\
        &(s(a),a)\smile (a^n,a)=0, & \quad (a^{2m},s(a))\smile (a^{2n},s(a))=(a^{2(m+n)},s(a)) \\
        &(a^{2m},s(a))\smile (a^{2n-1},a)=(a^{2(m+n)-1},a), &\quad (a^{2m-1},a)\smile (a^{2n-1},a)=0.
        \end{array}
        \]

        \item If $a^2\neq 0$, $\Char(\bbK)\neq 2$ and $|a|$ is odd, then:
        \[
        \HH^k(A)=
        \left\{
        \begin{array}{ll}
        \langle \mathds{1},(s(a),a^{2m})\mid m\ge 1 \rangle    & \text{if } k=0, \\
        \langle(a,s(a)), (a,a^{2m-1})\mid m\ge 1\rangle     & \text{if } k= 1,\\
        0 &\text{if } k>1.
        \end{array}
        \right.
        \]
        It is clear that we have 
        \[
        \begin{array}{lll}
        &(s(a),a^{2m})\smile (a,s(a))=0,  &\quad (s(a),a^{2m})\smile (s(a),a^{2n})=(s(a),a^{2(m+n)}) \\
        &(s(a),a^{2m})\smile (a,a^{2n-1})=(a,a^{2(m+n)-1}), &  \quad (a,s(a))\smile (a,s(a))=0 \\
        &(a,s(a))\smile (a,a^{2n-1})=0, &\quad (a,a^{2m-1})\smile (a,a^{2n-1})=0,
        \end{array}
        \]
        for all $m,n\ge 1$.

         \item If $a^2\neq 0$, $\Char(\bbK)= 2$ or $|a|$ is even, then the cohomology vector spaces are as follows:
         \[
        \HH^k(A)=
        \left\{
        \begin{array}{ll}
\langle(s(a),a^m)\mid\ m\ge 0 \rangle     & \text{if } k=0, \\
        \langle(a,a^m)\mid m\ge 0\rangle     & \text{if } k= 1,\\
        0 &\text{if } k>1.
        \end{array}
        \right.
        \]
        Clearly,
        \[
        (s(a),a^m)\smile (s(a),a^n)=(s(a),a^{m+n}),\quad (s(a),a^m)\smile (a,a^n)=(a,a^{m+n}), 
        \]
        \[
        \text{and } \ (a,a^m)\smile (a,a^n)=0 \text{ for all } m,n\ge 0.
        \]
    \end{itemize}
\end{Rmk}

\subsection{Hochschild cohomology algebra}\

\medskip

In this section, we exhibit a presentation of the graded algebra $\HH^*(A)$. We let $\oC^{\basic}(\Gamma)$ be the set of those elements of $\oC(\Gamma)$ that are not proper powers of another element of $\oC(\Gamma)$. The elements of $\oC^{\basic}(\Gamma)$ are the primitive elements of $\oC(\Gamma)$ when $\Char(\bbK)=2$. But when $\Char(\bbK)\neq 2$, the elements of $\oC^{\basic}(\Gamma)$ are the primitive elements  $C\in \oC(\Gamma)$ such that $l(C)+|C|$ is even together with the squares of the primitive $\Gamma$-complete cycles $D$ with $l(D)+|D|$ odd. Let $\Sp^{\basic}$ be the subset of $\oC^{\basic}(\Gamma)$ whose elements  correspond to the special loops in $\Sp$. That is, if $\Char(\bbK)=2$ then $\Sp^{\basic}=\Sp$, otherwise $\Sp^{\basic}=\Sp^2$. Similarly, denote by $\oC^{\basic}(\calB)$ the set of elements of $\oC(\calB)$ that are not proper powers of the elements of $\oC(\calB)$.

As a generalization of the results obtained in \cite[Proposition 5.5]{CSSS24} for the gentle case, we have the following result.

\begin{Thm}\label{Thm: basis-as-alg}
The set $\mathscr{G}$ of cohomology classes of the following cocycles of $\bbK(\Gamma\| \calB)$ is a generating set for the graded algebra $\HH^*(A)$:
    \begin{cenum}
        \item the pairs $(s(\alpha),\alpha)\in \HH^{0,|\alpha|}(A)$ with $\alpha$ a $\calB$-maximal path in $(Q,I)$,
        
        \item the sums $\llangle {\alpha}_s\in \HH^{0,|\alpha|}(A)$ with $\alpha\in \oC^{\basic}(\calB)$,

        \item the pairs $(c,c)\in \HH^{1,0}(A)$ with $c$ an arrow in the complement of the union of the spanning tree $T$ and the set of  special loops $\Sp$,

        \item the pairs $(\gamma,\alpha)\in \HH^{l(\gamma),|\alpha|-|\gamma|}(A)$ with $\gamma$ a $\Gamma$-maximal element of $\Gamma$ and $\gamma$ and $\alpha$ neither beginning nor ending with the same arrow,

        \item the sums $\llangle{C}_{gr}\in \HH^{l(C),-|C|}(A)$ with $C\in \oC^{\basic}(\Gamma)\setminus\Sp^{\basic}$.
    \end{cenum} 
\end{Thm}

It follows from  the cup product table (see Table \ref{table: cup product} that the set $\mathscr{G}$ generates $\HH^*(A)$ as a graded algebra. From now on, we will always use the symbol \cnum{$i$}  for the $i$-th class of generators when enumerating the generators of $\scrG$ for all $1\le i\le 5$. 

\begin{Rmk}\label{Remark: quiver that except}(see \cite[Remark 5.6]{CSSS24}).
    The set $\mathscr{G}$ of Theorem~\ref{Thm: basis-as-alg} generates $\HH^*(A)$ minimally except when \emph{the quiver $Q$ is given by one vertex and one loop $a$, which is not a special loop, and either $\Char(\bbK)=2$, or ``$a^2\in R$ and $|a|$ is odd", or ``$a^2\not\in R$ and $|a|$ is even"}, in which cases the algebra is gentle.  In these three cases,  we have $(s(a),a)\smile (a,s(a))=(a,a)$, so that the generator $(a,a)$ listed in Theorem~\ref{Thm: basis-as-alg} is not really needed.
\end{Rmk}

Now, we will describe a set of relations that present the cohomology algebra. For this, note that each generator in $\mathscr{G}$ has a total degree given by the sum of its external and internal degrees. The total degree is important to determine the graded commutativity relations implicit in the following presentation. 

\begin{Thm}\label{Thm: relation of HH}
    Let $A$ be the graded skew-gentle algebra associated to a graded skew-gentle triple $(Q,R,\Sp)$. Suppose that the quiver $Q$ is not like those excluded by the hypotheses of Remark~\ref{Remark: quiver that except}. Suppose furthermore that the quiver is not the one with one vertex and one special loop. The graded cohomology algebra $\HH^*(A)$ is the quotient of the free graded commutative algebra generated by the set $\scrG$ of Theorem \ref{Thm: basis-as-alg} by the ideal generated by a set $\scrR$ consisting of the following elements:
    \begin{itemize}
        \item $u\smile v$, one for each choice of $u$ and $v$ in $\scrG$ except those in which
        \begin{itemize}[label=\textendash]
        \item $u=v=\llangle{C}_{gr}$ for some $C\in \oC^{\basic}(\Gamma)\setminus\Sp^{\basic}$,
        
        \item $u=v=\llangle{\alpha}_s$ for some $\alpha\in \oC^{\basic}(\calB)$,
        
        \item $u$ and $v$ are, in some order, $\llangle{C}_{gr}$ and $(c,c)$ for some $C\in \oC^{\basic}(\Gamma)\setminus\Sp^{\basic}$ and $c$ an arrow in $Q_1\setminus(T\cup \Sp)$ through which $C$ passes,

        \item $u$ and $v$ are, in some order, $\llangle{\alpha}_s$ and $(c,c)$ for some $\alpha\in \oC^{\basic}(\calB)$ and $c$ an arrow in $Q_1\setminus(T\cup \Sp)$ through which $\alpha$ passes,
        \end{itemize}

        \item $(c,c)\smile \llangle{C}_{gr}-(d,d)\smile \llangle{C}_{gr}$, one for each choice of two distinct arrows $c$ and $d$ in $Q_1\setminus(T\cup \Sp)$ and a complete cycle $C$ in $\oC^{\basic}(\Gamma)\setminus\Sp^{\basic}$ that passes both
        through $c$ and through $d$,

        \item $(c,c)\smile \llangle{\alpha}_s-(d,d)\smile \llangle{\alpha}_s$, one for each choice of two distinct arrows $c$ and $d$ in $Q_1\setminus(T\cup \Sp)$ and a cocomplete cycle $\alpha$ in $\oC^{\basic}(\calB)$ that passes both
        through $c$ and through $d$.
    \end{itemize}
    \end{Thm}

\begin{Proof}
Indeed, by Remark \ref{Remark: quiver that except} and the first case of Remark \ref{remark: cup product of one cycle quiver}, it is straightforward to check that the cases we excluded in the statement of Theorem \ref{Thm: relation of HH} do not satisfy Theorem \ref{Thm: relation of HH}. On the other hand, 
when $\Char(\bbK)\neq 2$, the quiver $Q$ consists only of one vertex and one loop $a\notin \Sp$ and either ``$a^2\notin R$ and $|a|$ is odd" or ``$a^2\in R$ and $|a|$ is even", the computation in Remark \ref{remark: cup product of one cycle quiver} shows that Theorem \ref{Thm: relation of HH} holds. We will therefore assume that the quiver is not given by a vertex and a loop.

    Consider the free graded commutative algebra  $\mathscr{H}$ generated by $\mathscr{G}$, and the canonical surjective morphism $\pi: \mathscr{H}\to \HH^*(A)$ with kernel $I$. It was proved in \cite[Theorem 5.8]{CSSS24} using a case by case argument that $I$ is generated by the relations listed in the theorem, for more details, see \cite[Theorem 5.8]{CSSS24}.
\end{Proof}
\begin{Rmk}
    In Theorem~\ref{Thm: relation of HH} there are three cases excluded by the hypotheses on the presentation. In all these cases, the quiver has one vertex and one loop $a$, which is not a special loop,
    \begin{itemize}
        \item [(1)] suppose $a\in \Sp$, in this situation the cohomology algebra $\HH^*(A)$ is freely generated as a graded commutative algebra by the classes of the elements $(s(a),a)$ of degree $0$, subject to the relation $(s(a),a)\smile (s(a),a)=(s(a),a)$,
    
        \item [(2)] if $a^2\in R$, $\Char(\bbK)=2$ or $|a|$ is odd, then  the cohomology algebra $\HH^*(A)$ is freely generated as a graded commutative algebra by the classes of the elements $(s(a),a)$ and $(a,s(a))$ of degree $|a|$ and $1-|a|$, subject to the relation $(s(a),a)\smile (s(a),a)=0$,

        \item [(3)] when $a^2\not \in R$,  and either $\Char(\bbK)\neq 2$ and $|a|$ is even or $\Char(\bbK)=2$,  the cohomology algebra $\HH^*(A)$ is freely generated as a graded commutative algebra by the classes of the elements $(s(a),a)$ and $(a,s(a))$ of degree $|a|$ and $1-|a|$, subject to the relation $(a,s(a))\smile (a,s(a))=0$.
    \end{itemize}
\end{Rmk}

\subsection{The Gerstenhaber bracket}\

\smallskip

In this section, we employ the comparison morphisms to carry the Gerstenhaber bracket introduced on the Hochschild cochain complex induced by the bar resolution in Section~\ref{section: TT calculus of dg algebra} over to the Hochschild cochain complex induced by the CS-projective resolution.

Given two elements $(\gamma,\alpha)\in \Gamma_m\| \calB$ and $(\eta,\beta)\in \Gamma_n\| \calB$ with $m\ge 1$ and $n\ge 0$, the \emph{Gerstenhaber bracket} $[(\gamma,\alpha),(\eta,\beta)]\in \bbK(\Gamma_{m+n-1}\| \calB)$ is defined by
\[
[(\gamma,\alpha),(\eta,\beta)]=(\gamma,\alpha)\circ(\eta,\beta)-(-1)^{(|\alpha|-|\gamma|+m-1)(|\beta|-|\eta|+n-1)}(\eta,\beta)\circ (\gamma,\alpha),
\]
where $(\gamma,\alpha)\circ(\eta,\beta)=\sum\limits_{i=1}^m (\gamma,\alpha)\circ_i(\eta,\beta)$, and the operation $\circ_i$ for $1\le i\le m$ on $\bbK(\Gamma\| \calB)$ can be performed using the following diagram

\smallskip

\begin{eqnarray*}
\xymatrix{\bbK(\Gamma_m\| \calB)\ar[d]^{\simeq} \ar@{}[r]|-*+{\displaystyle \times} & \bbK (\Gamma_n\| \calB) \ar[d]^{\simeq} \ar[rr]^{\circ_i} && \bbK (\Gamma_{m+n-1}\| \calB)\\
\rmHom_{A^e}(A\ot\bbK \Gamma_m\ot A,A)\ar[d]^{\circ G^m} \ar@{}[r]|-*+{\displaystyle \times}& \rmHom_{A^e}(A\ot\bbK \Gamma_n\ot A,A)\ar[d]^{\circ G^n} &&\rmHom_{A^e}(A\ot\bbK \Gamma_{m+n-1}\ot A,A) \ar[u]_{\simeq}
\\
\rmHom_{A^e}(A\ot (s\oA)^m\ot A,A)\ar[d]^{\simeq} \ar@{}[r]|-*+{\displaystyle \times}& \rmHom_{A^e}(A\ot (s\oA)^n\ot A,A)\ar[d]^{\simeq} && \rmHom_{A^e}(A\ot (s\oA)^{m+n-1}\ot A,A)\ar[u]_{\circ F^{m+n-1}}\\
\rmHom_{E^e}((s\oA)^m,A)\ar@{}[r]|-*+{\displaystyle \times} & \rmHom_{E^e}( (s\oA)^n,A)\ar[rr]^{\circ_i} &&  \rmHom_{E^e}( (s\oA)^{m+n-1},A)\ar[u]_{\simeq}.}
\end{eqnarray*}

\medskip

That is, $f_{(\gamma,\alpha)\circ_i (\eta,\beta)}=\big( (f_{(\gamma,\alpha)}\circ G^m)\circ_i(f_{(\eta,\beta)}\circ G^n) \big)\circ F^{m+n-1}$. 

Assume $1\le i \le m$ and $x=x_{m+n-1}\cdots x_1\in \Gamma_{m+n-1}$. If $x_{i+n-1}\cdots x_i=\eta$, then:
\begin{align}\label{equation: gerstenhaber bracket 1}
    f_{(\gamma,\alpha)\circ_i (\eta,\beta)}(1\ot x\ot 1)=(-1)^{\sigma_i} f_{(\gamma,\alpha)}\circ G^m (1\ot s x_{m+n-1}\ot \cdots s x_{i+n}\ot s\beta\ot s x_{i-1}\ot\cdots\ot sx_1\ot 1),
\end{align}
where 
\small
\[\sigma_i=\sum\limits_{j=2}^{m+n-1}(|x_j|+\cdots+|x_{m+n-1}|)+(|x_{i+n}|+\cdots+|x_{m+n-1}|+m-i)(|\beta|-|\eta|+n+1)+\sum\limits_{j=i+1}^{i+n-1}(|x_j|+\cdots+|x_{i+n-1}|); \]
\normalsize
while if $x_{i+n-1}\cdots x_i\neq\eta,$
then 
$f_{(\gamma,\alpha)\circ_i (\eta,\beta)}(1\ot x\ot 1)=0.$

\begin{Fact}
    Given two elements $(\gamma,\alpha)\in \Gamma_m\| \calB$ and $(\eta,\beta)\in \Gamma_n\| \calB$ with $m\ge 1$ and $n\ge 0$, if $l(\beta)=0$, then $(\gamma,\alpha)\circ (\eta,\beta)=0$. Indeed, this 
follows directly from Equation~(\ref{equation: gerstenhaber bracket 1}).
\end{Fact}

\begin{Not}
    Let $\gamma=\gamma_m\cdots \gamma_1\in Q_m$ with $m\ge 1$, and let $\beta$ be a path in $Q$ and $a\in Q_1$, for any $1\le j\le m$ we denote 
    \[
    \gamma\vee_j^a \beta:=
    \left\{
    \begin{array}{ll}
    \gamma_m \cdots \gamma_{j+1} \beta \gamma_{j-1} \cdots \gamma_1 
    & \text{if } \gamma_j=a,\\
    0   
    & \text{otherwise}.
    \end{array}
    \right.
    \]
 \end{Not}

\begin{Prop}\label{prop: operation circ i}
    Let  $(\gamma,\alpha)\in \Gamma_m\| \calB$ and $(\eta,\beta)\in \Gamma_n\| \calB$ with $m\ge 1$, $n\ge 0$ and $l(\beta)\ge 1$. Assume $\gamma=\gamma_m\cdots \gamma_1$ and $\beta=\beta_{l(\beta)}\cdots \beta_1$. Then $(\gamma,\alpha)\circ_i (\eta,\beta)$ is  as follows:
    
\begin{itemize}
    \item [$(a)$] if $m=1$, then $(\gamma,\alpha)\circ_1 (\eta,\beta)=\sum\limits_{j=1}^{l(\beta)} (-1)^{(|\alpha|-|\gamma|)|\beta_{l(\beta)}\cdots \beta_{j+1}|} (\eta,\pi(\beta\vee_j^{\gamma} \alpha))$,

    \item [$(b)$] if $m>1$, then 
     \begin{align*}
     &(\gamma,\alpha)\circ_1 (\eta,\beta)
     =\left\{
     \begin{array}{ll}
     (-1)^{(|\gamma|-|\gamma_1|+m-1)|\beta|+m-1}(\gamma\vee_1^{\beta_{l(\beta)}} \eta, \pi(\beta\vee_{l(\beta)}^{\gamma_1}\alpha))  & \text{if } n=0,  \\
     (-1)^{(|\gamma|-|\gamma_1|+m-1)(|\beta|-|\eta|)+(m-1)(n-1)}(\gamma\vee_1^{\eta_n} \eta, \pi(\beta\vee_{l(\beta)}^{\eta_n}\alpha))  & \text{if } n\ge 1,
     \end{array}
    \right.
     \end{align*}
    and 
     \begin{align*}
     &(\gamma,\alpha)\circ_m (\eta,\beta)
    =\left\{
    \begin{array}{ll}
    (-1)^{(m-1)|\beta|+(|\alpha|-|\gamma|)(|\beta|-|\beta_1|)} (\gamma\vee_m^{\beta_1}\eta, \pi( \beta\vee_1^{\gamma_m}\alpha))  & \text{if } n=0,  \\
     (-1)^{(m-1)(|\beta|+|\eta|)+(|\alpha|-|\gamma|)(|\beta|-|\beta_1|)} (\gamma\vee_m^{\eta_1}\eta, \pi( \beta\vee_1^{\eta_1}\alpha))  & \text{if } n\ge 1,
     \end{array}
     \right.
     \end{align*}
    
    \item [$(c)$] if $m>2$ and $2\le i\le m-1$, then
     \begin{align*}
     &(\gamma,\alpha)\circ_i (\eta,\beta)\\
     &=\left\{
     \begin{array}{ll}
     (-1)^{(i+1)|\beta|+m-i}( \beta^{m-1},\alpha )  & \text{if $n=0$ and $\gamma=\beta^m$,}  \\
      (\gamma,\alpha) &\text{if $n=1$ and $\eta=\beta=\gamma_i$,} \\
      (-1)^{(|\gamma_m\cdots\gamma_{i+1}|+m-1)|\delta|+(m-i)(n-1)} ( \gamma_m\cdots\gamma_{i+1}\delta\gamma_{i-1}\cdots\gamma_1,\alpha ) & \text{if $n>1$, $\beta=\gamma_i$, $\eta=\beta\delta$ with $\delta$}\\
      & \text{a complete cycle starts with $\beta$,}\\
      0 & \text{otherwise}.
     \end{array}
     \right.
     \end{align*}
     \end{itemize}
\end{Prop}

\begin{Proof}
    \begin{itemize}
        \item [(a)] Let $m=1$. By Equation~(\ref{equation: gerstenhaber bracket 1}), we have      
        \begin{align*}
            f_{(\gamma,\alpha)\circ_1 (\eta,\beta)}(1\ot x\ot 1)
            &=\delta_{x,\eta} f_{\gamma,\alpha}\circ G^1(1\ot s\beta\ot 1)\\
            &=\delta_{x,\eta} \sum_{j=1}^{l(\beta)} f_{\gamma,\alpha}(\beta_{l(\beta)} \cdots \beta_{j+1} \ot \beta_j \ot \beta_{j-1}\cdots\beta_1)\\
            &= \delta_{x,\eta} \sum_{j=1}^{l(\beta)} (-1)^{(|\alpha|-|\gamma|)|\beta_{l(\beta)}\cdots \beta_{j+1}|} \delta_{\beta_j,\gamma} \pi(\beta_{l(\beta)} \cdots \beta_{j+1} \alpha \beta_{j-1}\cdots\beta_1),
        \end{align*}
      which implies $(\gamma,\alpha)\circ_1 (\eta,\beta)=\sum\limits_{j=1}^{l(\beta)} (-1)^{(|\alpha|-|\gamma|)|\beta_{l(\beta)}\cdots \beta_{j+1}|} (\eta,\pi(\beta\vee_j^{\gamma} \alpha))$.

        \item [(b)] Assume $m>1$. By Equation~(\ref{equation: gerstenhaber bracket 1}), if $n=0$, we have 
        \begin{align*}
            f_{(\gamma,\alpha)\circ_1 (\eta,\beta)}(1\ot x\ot 1)
            &=(-1)^{\sigma_1} f_{\gamma,\alpha}\circ G^m(1\ot sx_{m-1}\ot \cdots sx_1 \ot s\beta\ot 1)\\
            &=(-1)^{(|\gamma|-|\gamma_1|+m-1)|\beta|+m-1} \pi(\alpha \beta_{l(\beta)-1}\cdots \beta_1),
        \end{align*}
        if $\eta =s(x_1)$ and $\gamma=x_{m-1}\cdots x_1\beta_{l(\beta)}$.
        Otherwise $f_{(\gamma,\alpha)\circ_1 (\eta,\beta)}(1\ot x\ot 1)=0$. So,
        \[
        (\gamma,\alpha)\circ_1 (\eta,\beta)=(-1)^{(|\gamma|-|\gamma_1|+m-1)|\beta|+m-1} (\gamma\vee_1^{\beta_{l(\beta)}} \eta, \pi(\alpha \beta_{l(\beta)-1}\cdots \beta_1)).
        \]

        If $n\ge 1$, then
        \begin{align*}
            f_{(\gamma,\alpha)\circ_1 (\eta,\beta)}(1\ot x\ot 1)
            &=(-1)^{\sigma_1} \delta_{x_n\cdots x_1,\eta} f_{\gamma,\alpha}\circ G^m(1\ot sx_{m+n-1}\ot \cdots sx_{n+1} \ot s\beta\ot 1)\\
            &=(-1)^{(|\gamma|-|\gamma_1|+m-1)(|\beta|-|\eta|)+(m-1)(n-1)} \pi(\alpha \beta_{l(\beta)-1}\cdots \beta_1) ,
        \end{align*}
        if $\eta =x_n\cdots x_1$ and $\gamma=x_{m+n-1}\cdots x_{n+1}\beta_{l(\beta)}$. Otherwise $f_{(\gamma,\alpha)\circ_1 (\eta,\beta)}(1\ot x\ot 1)=0$ and
        \[
        (\gamma,\alpha)\circ_1 (\eta,\beta)=(-1)^{(|\gamma|-|\gamma_1|+m-1)(|\beta|-|\eta|)+(m-1)(n-1)} (\gamma\vee_1^{\beta_{l(\beta)}} \eta, \pi(\alpha \beta_{l(\beta)-1}\cdots \beta_1)),
        \]
        if $\gamma\vee_1^{\beta_{l(\beta)}} \eta\in\Gamma_{m+n-1}$, and zero if not.
       
\smallskip 

       The computation for  $(\gamma,\alpha)\circ_m (\eta,\beta)$ is analogous.
    
\smallskip

    \item [(c)] Assume $m>2$ and $2\le i\le m-1$.  
By Equation~(\ref{equation: gerstenhaber bracket 1}) and the definition of $G^m$, it follows that 
if $l(\beta)\ge 2$, then for any $x\in \Gamma_{m+n-1}$
\[
f_{(\gamma,\alpha)\circ_i (\eta,\beta)}(1\ot x\ot 1)=0,
\]
that is, $(\gamma,\alpha)\circ_i (\eta,\beta)=0$.  

Now suppose $\beta\in Q_1$. If $n=0$, then $\beta$ is a loop.  
In this case, we have
\[
\begin{aligned}
    f_{(\gamma,\alpha)\circ_i (\eta,\beta)}(1\ot x\ot 1)
    &=(-1)^{\sigma_i}\, f_{(\gamma,\alpha)}\circ G^m 
    (1\ot s x_{m-1}\ot \cdots \ot s x_i\ot s\beta\ot s x_{i-1}\ot\cdots\ot s x_1\ot 1)\\
    &=(-1)^{(|x_{m-1}\cdots x_i|+m-1)|\beta|+m-i}
    f_{(\gamma,\alpha)} (1\ot x_{m-1}\cdots x_i\beta x_{i-1}\cdots x_1\ot 1),
\end{aligned}
\]
if $\eta=s(x_i)$ and $x_{m-1}\cdots x_i\beta x_{i-1}\cdots x_1\in \Gamma_m$,  
as a consequence, $\beta^2\in S$ and $x=\beta^{m-1}$.  
Hence, when $n=0$, we obtain
\[
(\gamma,\alpha)\circ_i (\eta,\beta)
=(-1)^{(i+1)|\beta|+m-i}\,(\beta^m,\alpha),
\]
provided that $\gamma=\beta^m$.

When $n\ge 1$, we have
\begin{align*}
        &f_{(\gamma,\alpha)\circ_i (\eta,\beta)}(1\ot x\ot 1)\\
        &=(-1)^{\sigma_i} f_{(\gamma,\alpha)}\circ G^m (1\ot s x_{m+n-1}\ot \cdots s x_{i+n}\ot s\beta\ot s x_{i-1}\ot\cdots\ot sx_1\ot 1)\\
        &=(-1)^{(|x_{m+n-1}\cdots x_{i+n}|+m-1)(|\beta|-|\eta|)+(m-i)(n-1)} f_{(\gamma,\alpha)} (1\ot x_{m+n-1}\cdots x_{i+n}\beta x_{i-1}\cdots x_1\ot 1).
    \end{align*}
This holds whenever $\eta=x_{i+n-1}\cdots x_i$ and 
$x_{m+n-1}\cdots x_{i+n}\beta x_{i-1}\cdots x_1\in \Gamma_m$.  
This condition implies that $\beta=\eta=x_i$ if $n=1$,  
and that $\beta=x_{i+n-1}=x_i$ and $\eta=\beta \delta$ with $\delta=x_{i+n-2}\cdots x_i$ 
a complete cycle starting from $\beta$ when $n>1$.  

Therefore, we further obtain
\[
f_{(\gamma,\alpha)\circ_i (\eta,\beta)}(1\ot x\ot 1)
=(-1)^{(|\gamma_m\cdots\gamma_{i+1}|+m-1)(|\beta|-|\eta|)+(m-i)(n-1)}\,\alpha,
\]
if $\gamma=x_{m+n-1}\cdots x_{i+n}\beta x_{i-1}\cdots x_1$;  
otherwise, $f_{(\gamma,\alpha)\circ_i (\eta,\beta)}(1\ot x\ot 1)=0$.

As a conclusion, when $n=1$, we have
\[
(\gamma,\alpha)\circ_i (\eta,\beta)
=(\gamma,\alpha),
\quad\text{if } (\eta,\beta)=(\gamma_i,\gamma_i).
\]
When $n>1$, we obtain
\[
(\gamma,\alpha)\circ_i (\eta,\beta)
=(-1)^{(|\gamma_m\cdots\gamma_{i+1}|+m-1)|\delta|+(m-i)(n-1)} 
(\gamma_m\cdots\gamma_{i+1}\delta\gamma_{i-1}\cdots\gamma_1,\alpha),
\]
if $(\eta,\beta)=(\gamma_i\delta,\gamma_i)$, where $\delta$ is a complete cycle starting from $\gamma_i$.  
In all other cases, $(\gamma,\alpha)\circ_i (\eta,\beta)=0$.
\end{itemize}
\end{Proof}

\begin{Not}
If $c$ is an arrow in $Q$ and $\gamma=\gamma_m\cdots \gamma_1$ is a path, we set 
\[
\deg_c(\gamma):=\#\{i\in\{1,\ldots,m\}\ |\ \gamma_i=c \},
\]
the number of times the path $\gamma$ passes through $c$, and if $(\gamma,\alpha)$ is an element of $\Gamma\| \calB$ we let
\[
\deg_c(\gamma,\alpha):=\deg_c(\alpha)-\deg_c(\gamma).
\]
\end{Not}

We may use  Proposition~\ref{prop: operation circ i} to compute the operation~$\circ$ between any two generators of the algebra~$\HH^*(A)$. This yields the following three corollaries, which determine the table of the operation~$\circ$. The proof is obtained by a straightforward, step-by-step calculation from the formula in Proposition~\ref{prop: operation circ i}, and we give a short proof.

\begin{Cor}\label{cor: circ op 1}
    Let $c\in Q_1\setminus(T\cup \Sp)$, and $(\eta,\beta)\in \Gamma_n\| \calB$ with $n\ge 0$, then 
    \[
    (c,c)\circ (\eta,\beta)=\deg_c(\beta) (\eta,\beta).
    \]
\end{Cor}

\begin{Proof}
    It follows from  Proposition~\ref{prop: operation circ i} (1).
\end{Proof}

\begin{Cor}\label{cor: circ op 2}
    Let $(Q,R,\Sp)$ be a graded skew-gentle triple, $(\gamma,\alpha)\in \Gamma_m\| \calB$, $m\ge 1$, with $\gamma$ $\Gamma$-maximal and $\gamma$ and $\alpha$ neither beginning nor ending with the same arrow. Assume  $(\eta,\beta)\in \Gamma_n\| \calB$ with $n\ge 0$. Then
    \begin{align*}
        (\gamma,\alpha)\circ (\eta,\beta)=
        \left\{
        \begin{array}{ll}
        (\eta,\alpha)  & \text{if $l(\gamma)=1$ and $\beta=\gamma,$} \\
        (\gamma,\alpha)     & \text{if $l(\gamma)>1$ and  $\eta=\beta=\gamma_i $ for some $1\le i\le m,$}\\
       (\gamma_m\dots\gamma_2,\alpha) (\text{coboundary })
       &\text{if $l(\gamma)>1$, $l(\eta)=0$ and $\beta=\gamma_1$ or $\beta=\gamma_m$,}\\
        0  &\text{otherwise}.
        \end{array}
        \right.
    \end{align*}
\end{Cor}

\begin{Proof}
    Assume that $l(\gamma)=1$. By Proposition~\ref{prop: operation circ i} (1), we have 
    \[
    (\gamma,\alpha)\circ_1 (\eta,\beta)=\sum\limits_{j=1}^{l(\beta)} (-1)^{(|\alpha|-|\gamma|)|\beta_{l(\beta)}\cdots \beta_{j+1}|} (\eta,\pi(\beta\vee_j^{\gamma} \alpha)),
    \]
   The $\Gamma$-maximality of $\gamma$ implies that $\pi(\beta\vee_j^{\gamma}\alpha)=0$ when $l(\beta)>1$. When $l(\beta)=1$, we have $\pi(\beta\vee_1^{\gamma}\alpha)=\alpha$ if and only if $\beta=\gamma$. In this case, we obtain
   \[
   (\gamma,\alpha)\circ (\eta,\beta)=(\eta,\alpha).
   \]

   Assume $l(\gamma) > 1$. 
If $l(\beta) > 1$, then by the $\mathcal{A}$-maximality of $\gamma$ and the fact that  $\gamma$ and $\alpha$ neither begin nor end with the same arrow, we have 
\[
\pi(\beta \vee_{l(\beta)}^{\gamma_1} \alpha) = 0 = \pi(\beta \vee_1^{\gamma_m} \alpha).
\]
Therefore, $(\gamma, \alpha) \circ (\eta, \beta) = 0$.

Next, let $l(\beta) = 1$. Under this assumption, if $n = 0$, then $\beta$ is a loop. 
It is clear that $(\gamma, \alpha) \circ_i (\eta, \beta) = 0$ for $2 \le i \le m - 1$. 
For $(\gamma, \alpha) \circ_1 (\eta, \beta)$, we have
\[
(\gamma, \alpha) \circ_1 (\eta, \beta) =
\begin{cases}
(\gamma_m \cdots \gamma_2, \alpha) & \text{if } \gamma_1 = \beta, \\[4pt]
0 & \text{otherwise.}
\end{cases}
\]
If $\gamma_1 = \beta$, then by  definition of skew-gentle algebras, we have $\alpha_1 = \gamma_2$, 
which implies that $(\gamma_m \cdots \gamma_2, \alpha)$ is a coboundary. 
A similar argument applies to $(\gamma, \alpha) \circ_m (\eta, \beta)$.

Hence, we may assume that $l(\eta) \ge 1$. 
If there is no $i \in \{1, \ldots, m\}$ such that $\eta = \beta = \gamma_i$, 
then $\gamma \vee_i^{\beta} \eta = 0$, and thus $(\gamma, \alpha) \circ (\eta, \beta) = 0$. 
Conversely, if there exists $i \in \{1, \ldots, m\}$ such that $\eta = \beta = \gamma_i$, 
then there is only one such $i$, $\gamma \vee_i^{\beta} \eta = \gamma$, and therefore $(\gamma, \alpha) \circ (\eta, \beta) = (\gamma, \alpha)$.
    \end{Proof}

\begin{Cor}\label{cor: circ op 3}
    Let $(Q,R,\Sp)$ be a graded skew-gentle triple and suppose that $Q$ is not the quiver consisting of one vertex and one loop. Let $\llangle{D}_{gr}\in \oC^{\basic}(\Gamma)\setminus\Sp^{\basic}$ and $v\in \mathscr{G}$, then
    \begin{align*}
        \llangle{D}_{gr} \circ v=
        \left\{
        \begin{array}{ll}
        \deg_c(D)\llangle{D}_{gr}     & \text{if $v=(c,c)$ for some $c\in Q_1\setminus (T\cup \Sp)$}, \\
        0     & \text{otherwise}.
        \end{array}
        \right.
    \end{align*}
    Moreover, when $l(D)=1$, $v=(\gamma,\alpha)$ with $\gamma$ a $\Gamma$-maximal element of $\Gamma$ and $\gamma$ and $\alpha$ neither begin nor end with the same arrow, and $\alpha$ begins or ends with $D$, the cocycle $\llangle{D}_{gr}\circ v$ is a coboundary. 
\end{Cor}

\begin{Proof}
We will only provide the proof of the last statement, that is, for the case where $l(D) = 1$ and $v = (\gamma, \alpha)$ 
with $\gamma$ an $\Gamma$-maximal element of $\Gamma$, 
and $\gamma$ and $\alpha$ neither begin nor end with the same arrow. 
The proofs of the other cases are simpler and will not be discussed in detail.

If  $\alpha$ does not pass through $D$, then it is clear that 
$\alpha \vee_i^D s(D) = 0$ for any $i \in \{1, l(\alpha)\}$, 
and hence $\llangle D_{gr} \circ (\gamma, \alpha) = 0$.

If, on the contrary, $\alpha$ passes through $D$, we consider three cases:

\textbf{Case 1:} $\alpha$ begins with $D$, that is, $\alpha_1 = D$. 
   By Proposition~\ref{prop: operation circ i}, we have
   \[
   \llangle D_{gr} \circ (\gamma, \alpha) = \pm (\gamma, \alpha_{l(\alpha)} \cdots \alpha_2).
   \]
   This is a coboundary, since, by definition of skew-gentle algebras, 
   $\gamma$ starts with $\alpha_2$.

   \textbf{Case 2:} There exists $1 < i < l(\alpha)$ such that $\alpha_i = D$. 
   In this case, it is clear that $\pi(\alpha \vee_i^D s(D)) = 0$.

\textbf{Case 3:} $\alpha$ ends with $D$. This case is similar to the first one.
\end{Proof}

We obtain the table of the operation $\circ$ from Corollaries~\ref{cor: circ op 1}, \ref{cor: circ op 2} and \ref{cor: circ op 3}.
\begin{table}[H]
\centering
{\renewcommand{\arraystretch}{1.5}
\begin{tabular}{|c|c|c|c|c|c|}
\hline
\rowcolor{gray!30} 
$\circ$ & $(s(\alpha),\alpha)$ & $\llangle{\alpha}_s$ & $(c,c)$ & $(\gamma,\alpha)$ & $\llangle{C}_{gr}$ \\ \hline
\cellcolor{gray!30} $(s(\alpha'),\alpha')$ & $0$ & $0$ & $0$ & $0$ & $0$ \\ \hline
\cellcolor{gray!30} $\llangle{\alpha'}_s$ & $0$ & $0$ & $0$ & $0$ & $0$ \\ \hline
\cellcolor{gray!30} $(d,d)$ & \cellcolor{blue!15} $\deg_d(\alpha)(s(\alpha),\alpha)$ & \cellcolor{blue!15} $\deg_d(\alpha)\llangle{\alpha}_s$ & \cellcolor{blue!15} $\deg_d(c)(c,c)$ & \cellcolor{blue!15} $\deg_d(\alpha)(\gamma,\alpha)$ & $0$ \\ \hline
\cellcolor{gray!30} $(\gamma',\alpha')$ & $0$ & $0$ &\cellcolor{blue!15}  $\deg_c(\gamma')(\gamma',\alpha')$ & $0$ & $0$ \\ \hline
\cellcolor{gray!30} $\llangle{D}_{gr}$ & $0$ & $0$ &\cellcolor{blue!15} $\deg_c(D)\llangle{D}_{gr}$ & $0$ & $0$ \\ \hline
\end{tabular}}
\caption{The operation $\circ$ between the generators of $\HH^*(A)$, for $Q$ neither the quiver consisting of one loop nor the Kronecker quiver.}
\end{table}

\begin{Thm}\label{thm: bracket of HH}
    Suppose that $Q$ is neither the quiver consisting of  one loop nor the Kronecker quiver. 
    \begin{itemize}
    
        \item[$(1)$] If $(c,c)\in \HH^{1,0}(A)$ corresponds to an arrow $c$ in the complement of $T\cup\Sp$ and $v\in \HH^n(A)$, then $[(c,c),v]=\deg_c(v)v$.

        \item[$(2)$] All other brackets among elements of $\mathscr{G}$ are zero.
        
    \end{itemize}
    Moreover, in $(1)$ if $v$ is also in $\HH^1(A)$, then $\deg_c(v)$ is always $0$, $1$ or $-1$, hence the Lie algebra structure of $\HH^1(A)$ does not depend on $\Char(\bbK)$.
\end{Thm}

Now we treat the excluded cases.

\begin{Rmk}
    Suppose that the quiver $Q$ has exactly one vertex and one loop $a$. From Remark~\ref{remark: cup product of one cycle quiver} we obtain:
    \begin{itemize}
        \item  If $a\in \Sp$, then $\mathscr{G}=\{(s(a),a)\}$ and all brackets between elements of $\mathscr{G}$ are zero.

    \end{itemize}
Now we assume $a\not\in\Sp$. 
    \begin{itemize}    
        \item If either ($a^2\in R$ and either $\Char(\bbK)=2$ or $|a|$ is odd) or ($a^2\not \in R$ and either $\Char(\bbK)=2$ or $|a|$ is even), 
        then $\mathscr{G}=\{(s(a),a), (a,s(a))\}$ and the only non zero bracket between elements of $\mathscr{G}$ is
        \[
        [(a,s(a)),(s(a),a)]=(s(a),s(a)).
        \]

        \item  If $a^2\in R$, $\Char(\bbK)\neq 2$ and $|a|$ is even, then $\mathscr{G}=\{(s(a),a),(a,a), (a^2,s(a))\}$ and the non zero brackets between elements of $\mathscr{G}$ are
        \[
        [(a,a),(s(a),a)]=(s(a),a), \quad [(a,a),(a^2,s(a))]=-2(a^2,s(a)).
        \]

        \item  If $a^2\not \in R$, $\Char(\bbK)\neq 2$ and $|a|$ is odd, then $\mathscr{G}=\{(s(a),a^2),(a,a), (a,s(a))\}$ and the only non zero brackets between elements of $\mathscr{G}$ are
        \[
        [(s(a),a^2),(a,a)]=-2(s(a),a^2), \quad [(a,s(a)),(a,a)]=(a,s(a)).
        \]
            \end{itemize}
\end{Rmk}

\begin{Rmk}
    Suppose now that $Q$ is the Kronecker quiver with arrows $a$ and $b$. Then $\mathscr{G}=\{(a,a),(a,b),(b,a)\}$.
    \begin{itemize}
        \item If $|a|-|b|$ is odd or $\Char(\bbK)=2$, the only non zero brackets between elements of $\mathscr{G}$ are
        \[
        [(a,a),(a,b)]=-(a,b),\quad [(a,a),(b,a)]=(b,a).
        \]

        \item If $|a|-|b|$ is even and $\Char(\bbK)\neq 2$, the only non zero brackets between elements of $\mathscr{G}$ are
        \[
        [(a,a),(a,b)]=-(a,b),\quad [(a,a),(b,a)]=(b,a), \quad [(a,b),(b,a)]=-2(a,a).
        \]
    \end{itemize}
\end{Rmk}

We now give an example.

\begin{Ex}\label{Example:HH^*}
Let $A=\bbK Q/\langle R\cup\{\varepsilon^2-\varepsilon\mid \varepsilon\in\Sp\} \rangle$ be a graded skew-gentle algebra corresponding to the following skew-gentle triple $(Q,R,\Sp)$ with $R=\{ba,cb\}$ and $\Sp=\{\varepsilon_2,\varepsilon_3\}$.
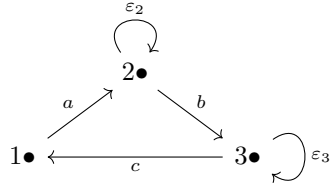
\[
\begin{tikzcd}
    & 2\bullet \arrow[rd, "b"] \arrow["\varepsilon_2", loop, distance=2em, in=55, out=125] & \\ 1\bullet \arrow[ru, "a"] &  & 3\bullet \arrow[ll, "c"] \arrow["\varepsilon_3", loop, distance=2em, in=325, out=35]
\end{tikzcd}
\]
There is no $\calB$-maximal cycle, and $cba$ is the unique $\Gamma$-maximal element. Set $\alpha:=c\varepsilon_3b\varepsilon_2a$. We have  

\[
\oC^{\basic}(\calB)=\begin{cases}
		\{\alpha\} & \text{if $\Char(\bbK)=2$ or $|\alpha|$ even,} \\  \{\alpha^2 \} & \text{otherwise,}
	\end{cases} \text{ and }
    \oC^{\basic}(\Gamma)=\Sp^{\basic}=\begin{cases}
		\Sp & \text{if $\Char(\bbK)=2$,}\\
		\Sp^2 & \text{otherwise.}
	\end{cases}
\]
    
    If we choose the spanning tree $T=\{a,b\}$, then
	the generating set $\mathscr{G}$ in Theorem \ref{Thm: basis-as-alg} of $\HH^*(A)$ as an algebra consists of: 
    \begin{cenum}[start=2]
		\item $\begin{cases}
			\llangle{\alpha^2}_s\in \HH^{0,2|\alpha|}(A)\subset \HH^{2|\alpha|}(A)  & \text{if $\Char(\bbK)\neq 2$ and $|\alpha|$ odd,} \\
			\llangle{\alpha}_s\in \HH^{0,|\alpha|}(A)\subseteq \HH^{|\alpha|}(A) & \text{otherwise.}
		\end{cases}$
		\item $(c,c)\in \HH^{1,0}(A)\subseteq \HH^1(A)$.
		\item $(cba,e_1)\in \HH^{3,-|\alpha|}(A)\subseteq\HH^{3-|\alpha|}(A)$.
	\end{cenum}

    \smallskip
    
    We recall that the enumeration above  corresponds to Theorem \ref{Thm: basis-as-alg}.
	The relation set $\scrR$ in Theorem \ref{Thm: relation of HH} is  $\scrR=\left\{\cnum{2}\smile \cnum{4}, \cnum{3}\smile \cnum{3}, \cnum{3}\smile \cnum{4}, \cnum{4}\smile \cnum{4} \right\}$. Therefore, there is an isomorphism of algebras. \\
   \[\HH^*(A)\cong \bbK\mathscr{G}/\langle\scrR\rangle.\]
    In addition, the all non zero brackets in  Theorem \ref{thm: bracket of HH} between elements of $\scrG$ are 
    \[[\cnum{3},\cnum{4}]=-\cnum{4}\quad \text{and}\quad [\cnum{3},\cnum{2}]=\begin{cases}
        2\cdot\cnum{2} & \text{if $\Char(\bbK)\neq 2$ and $|\alpha|$ odd,}\\ \cnum{2} & \text{otherwise.}
    \end{cases} \]

\end{Ex}

\bigskip

\section{Geometric interpretation}\label{Sec: geometric}

In  \cite{A21, AB22, LSV22} it is shown that—up to small cases—there is a bijection between the set of isomorphism classes of  skew-gentle algebras and certain dissections of orbifold surfaces with orbifold points of order 2 and that this gives a geometric model for the bounded derived category of skew-gentle algebras. This has been extended to the perfect derived categories of graded skew-gentle algebras in \cite{AP24, CK24, BSW24, QZZ}. More precisely, in \cite{AP24, CK24, BSW24}, the partially wrapped Fukaya category of an orbifold surface with orbifold points of order 2 has been given in three different (equivalent) ways, namely by an explicit $A_\infty$-algebra from the surface dissection resulting from the embedding of a ribbon graph, as the orbit category by a $\mathbb{Z}_2$-action of a (surface dissection of a) smooth double cover of the orbifold surface and as the global sections of a cosheaf of $A_\infty$-categories on a ribbon graph of the orbifold surface. It is then shown in \cite{AP24, BSW24, CK24} that for any of these constructions there is a particular (formal) generator of the category in the form of an orbifold dissection corresponding to a graded skew-gentle algebra. Furthermore, all associative graded algebras arising from different formal generators of the same graded orbifold surface are derived equivalent, and it is conjectured in \cite{BSW24} that any graded associative algebra derived equivalent to a graded skew-gentle algebra arises in this way. 
Since the Hochschild cohomology and its structure are invariant under derived equivalence \cite{Ri91}, it follows that our calculations give the Hochschild cohomology and its structure for any of the graded (not necessarily skew-gentle) derived equivalent algebras associated to an orbifold via different (formal) dissections of the orbifold surface and more generally, for the partially wrapped Fukaya category of the graded orbifold surface. We will show 
that, just as in the case of an ungraded gentle algebra, see \cite{CSSS24} and also \cite{BSW25} for the second Hochschild cohomology in the graded gentle case, for a graded skew-gentle algebra the generators of the Hochschild cohomology, the cup product and the Gerstenhaber bracket can be `read off' the orbifold surface associated to the skew-gentle algebra.

The approach that we take here is to construct a graded marked ribbon graph directly from the skew-gentle algebra, following the ideas  for the ungraded case in \cite{LSV22}. This gives rise to a graded marked orbifold surface giving a geometric model of the  perfect derived category of the graded skew-gentle algebra. Furthermore, as described in Section \ref{section: gentle vs skew-gentle}, see also \cite{BSW25} for a complete treatment of the deformation theory of graded gentle algebras, any graded skew-gentle algebra can be obtained as a deformation of a graded gentle algebra. In particular, the special loops are deformations of loops $\varepsilon$ in the graded gentle algebra with $\varepsilon^2 = 0$. In terms of the associated graded marked surface of the graded gentle algebra such a loop (which then necessarily has grading 0) gives rise to a boundary component with one marked point and winding number 1.  It is shown in \cite{BSW25} that the $A_\infty$-deformation along the corresponding Hochschild 2-cocycle corresponds in the surface to replacing that boundary component by an orbifold point of order 2. 

\subsection{Graded orbifold surfaces from graded skew-gentle algebras}
We now introduce the geometric model of a graded skew-gentle algebra following the ideas in \cite{LSV22} and \cite{BSW24}.  
This is based on a graded marked ribbon graph with a set of special vertices. This ribbon graph for skew-gentle algebras is an adaptation of the ribbon graph constructed in \cite{S15} for gentle algebras. 

Given a graded skew-gentle triple $(Q,R,\Sp)$, let $A=\bbK Q/I$ be the corresponding graded skew-gentle algebra where $I$ denotes the two-sided ideal of $\bbK Q$ generated by $R\cup \{ \varepsilon^2-\varepsilon\ \mid \varepsilon\in\Sp \}$. Let $A':=\bbK Q/\langle R\cup \{\varepsilon^2\mid \varepsilon\in \Sp\} \rangle$ be the associated gentle algebra as in Section~\ref{section: gentle vs skew-gentle}.   We define the set of \emph{finite length Sp-maximal paths} as follows: we first determine all maximal paths in the associated graded gentle quiver  $(Q,R \cup \{\varepsilon^2\mid \varepsilon\in \Sp\} )$, that is the set of finite length non zero  paths $p$ such that for any arrow $x \in Q_1 $, $xp=px=0$ in $\bbK Q/\langle R\cup \{\varepsilon^2 \mid \varepsilon\in \Sp\} \rangle$. The set of finite length Sp-maximal paths is then given by the subpaths of the finite length  maximal paths where we have deleted all special loops.  The set of \emph{infinite Sp-maximal paths} is the rotation class of subcycles of primitive cycles without relations where we have deleted all special loops.   An example of an infinite Sp-maximal path is the path $cba$ in Example \ref{Example:HH^*}, which is a subpath of the (primitive) cycle with no relations $c\varepsilon_3b\varepsilon_2a$. Note that a Sp-maximal path contains length two relations, namely everywhere where we have deleted a special loop. Thus it is a non zero path in $Q$ but not in $A$ or $A'$. We then  construct a graded marked  ribbon graph $G$ associated to $(Q,R,\Sp)$ as follows:
\begin{itemize}
    \item the edges of $G$ are in bijection with  the vertices of $Q$,
    \item the vertices of $G$ are given by the set consisting of 
    \begin{itemize}
     \item the set of finite length Sp-maximal paths, 
     \item the set of infinite Sp-maximal paths,
     \item the trivial paths $e_i$ such that $i$ is either the source or the target of only one arrow, or $i$ is the source of exactly one arrow $\alpha$ and the target of exactly one arrow $\beta $, and $\alpha \beta \notin R$,
     \item the set of special loops $\varepsilon_i$ where $i$ is a special vertex in $(Q,R,\Sp)$. We call this type of vertex a \emph{special vertex} of $G$ and denote it by $\times$,
    \end{itemize}
    \item an edge $i$ of $G$ is incident to a vertex $q$ of $G$ if the corresponding path $q$ in $Q$ passes through the vertex $i$ in $Q$, 
    \item an edge adjacent to a special vertex is called a \emph{special edge}, 
    \item the cyclic ordering $\sigma$ of the edges incident to a vertex $q$ of $G$ is induced by linear (or cyclic) order in which $q$ as a path visits the corresponding vertices, 
    \item all vertices of $G$ apart from those corresponding to a special loop or an infinite Sp-maximal path are marked.
    In order to define the marking, first consider the case that $q$ is a finite Sp-maximal path that is not a cycle, then we mark the pair $(t(q), \sigma(t(q)))$, where $\sigma(t(q))=s(q)$. If $q=q_m \ldots q_1$ is a cycle, then since $q$ is a finite Sp-maximal path there necessarily is a relation $q_1q_m$ at vertex $s(q)=t(q)$. By construction, the edge of $G$ corresponding to $s(q)$ is a loop and we can consider it as a pair of distinct half-edges corresponding to $t(q)$ and $s(q)$. In the cyclic ordering at $q$, we have $\sigma(t(q)) = s(q)$ and we mark the pair $(t(q), \sigma(t(q)))$ as before.   
    If $q$ corresponds to a trivial path,  then it is a leaf vertex and we simply set a mark next to it,
    \item each unmarked pair   $(i,\sigma(i))$ corresponds to an arrow  $\alpha :i\rightarrow \sigma(i)$ and we grade such a pair by associating to it the integer $|\alpha |$. 
\end{itemize}

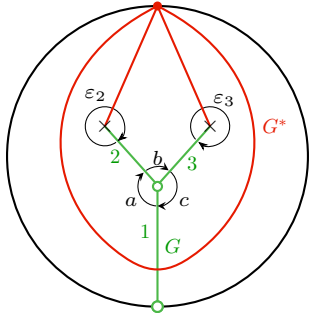
\begin{figure}[H]
    \centering
\begin{tikzpicture}[scale=2,transform  shape, every node/.style={scale=0.55}]

\def\R{1}          
\def\rnodeA{0.035} 

\draw[thick] (0,0) circle (\R);

\draw[yellow!30!green, thick] (0,-1) circle (\rnodeA); 
\fill[white] (0,-1) circle (\rnodeA);                 

\fill[green!10!red] (0,1) circle (0.03);
\coordinate (P) at (-0.35,0.2);
\coordinate (P') at (0.35,0.2);

\draw[yellow!30!green, thick]
  (0,-0.97) -- (0,-0.225);              

\draw[yellow!30!green, thick]
  (-0.01,-0.18) -- (-0.35,0.2);         

\draw[yellow!30!green, thick]
  (0.01,-0.18) -- (0.35,0.2);           

\coordinate (V) at (0,-0.2);

\def\V{(0,-0.2)}    
\def\r{0.13}         

\draw[-stealth] 
  (0,-0.2) ++(48:\r) 
  arc[start angle=48, end angle=-90, radius=\r]; 

\draw[-stealth] 
  (0,-0.2) ++(125:\r)
  arc[start angle=125, end angle=55, radius=\r]; 

\draw[-stealth] 
  (0,-0.2) ++(-94:\r)
  arc[start angle=-94, end angle=-225, radius=\r]; 

\def\rloop{0.13}

\draw[-{Stealth[scale=0.8]}] 
  (P) ++(-55:\rloop)
  arc[start angle=-55, delta angle=-355, radius=\rloop];
  
\draw[-{Stealth[scale=0.8]}] 
  (P') ++(-138:\rloop)
  arc[start angle=-138, delta angle=-350, radius=\rloop];

\node[font=\scriptsize] at (280-95-50:1.65em) {$\varepsilon_2$};
\node[font=\scriptsize] at (96-6-50:1.65em) {$\varepsilon_3$};

\node[font=\scriptsize] at (0-60-60:1.0em) {$a$};
\node[font=\scriptsize] at (0-30-30:1.0em) {$c$};
\node[font=\scriptsize] at (0-30-60:0.030em) {$b$};

\node[left,color=green!60!black][scale=0.7] at (-0.2,0) {$2$};
\node[left,color=green!60!black][scale=0.7]  at (-0.0,-0.5) {$1$};
\node[right,color=green!60!black][scale=0.7] at (0.15,-0.05) {$3$};
\node[right,color=green!60!black][scale=0.7] at (0.0,-0.6) {$G$};
\node[right,color=green!10!red][scale=0.7] at (0.65,0.2) {$G^*$};

\def\rnodeA{0.035} 
\def\rnodeB{0.03}  

\draw[yellow!30!green, thick] (0,-1) circle (\rnodeA);
\draw[yellow!30!green, thick] (0,-0.2) circle (\rnodeB);

\draw[green!10!red, thick]
  (0,1)
    .. controls (-0.8,0.6) and (-0.8,-0.2) .. (-0.3,-0.61)
    .. controls (-0.07,-0.8) and (0.07,-0.8) .. (0.3,-0.61)
    .. controls (0.8,-0.2) and (0.8,0.6) .. (0,1);

\draw[green!10!red, thick]
  (0,1) -- (-0.35,0.2);
\draw[green!10!red, thick]
  (0,1) -- (0.35,0.2);

\node[scale=0.8] at (P)  {$\mathbf{\times}$};
\node[scale=0.8] at (P') {$\mathbf{\times}$};

\end{tikzpicture}
\caption{ This is the surface with two orbifold points of order 2 (depicted as 
    $\mathbf{\times}$) and ribbon graphs $G$ and $G^*$ associated to the graded skew-gentle algebra in Example~\ref{Example:HH^*}.}
\label{Fig:ribbon graphs}
\end{figure}

For example, in Example \ref{Example:HH^*}, the vertices of $G$ are  $ \{cba, e_1, \varepsilon_2, \varepsilon_3\}$ and the edges of $G$ are $\{1,2,3\}$, all of which are adjacent to the vertex $cba$.  The edge $1$ is also adjacent to $e_1$ and $2$ and $3$ are special edges that are, respectively, adjacent to the special vertices $\varepsilon_2$ and $\varepsilon_3$. The only marking is at vertex $e_1$. The corresponding surface is a disk with one marked point on the boundary, one marked point in the interior and two orbifold points of order 2. Note that the interior marked point corresponds to the vertex $cba$ of $G$ which is an infinite Sp-maximal path, see Figure~\ref{Fig:ribbon graphs}.

\begin{Rmk} (1) 
We note that in the above, if a pair  $(i,\sigma(i))$ of cyclically consecutive edges in $G$ corresponds  to a special loop in $(Q, R, \rm{Sp})$ (that is when $i = \sigma(i)$), then the grading of $(i,\sigma(i))$ is zero.

 (2) Note that each edge of $G$ is incident to exactly two vertices of $G$. This follows directly from the facts that  $(Q,R\cup\{\varepsilon^2\mid \varepsilon\in \Sp\})$ is a gentle pair and the close relationship of the associated ribbon graphs. 
 
 (3) When $\Sp=\emptyset$, the graded skew-gentle algebra $A$ is a graded gentle algebra, and the graded marked ribbon graph  coincides with the graded marked ribbon graph of the graded gentle algebra in \cite{S15}.

\end{Rmk}

A graded marked ribbon graph $G$ (which is not the trivial graph with one vertex and no edges), gives rise to a compact oriented surface $S$ with boundary, such that $G$ is a deformation retract of $S$ and the faces of $G$ are in bijection with the boundary components of $S$. The marking of $G$ then  gives a unique way of gluing the marked vertices of $G$ to the boundary of $S$. Note that all unmarked vertices stay in the interior of $S$. In particular, this holds for the special vertices of $G$.  We view the latter as orbifold points of order 2 and denote the set of these orbifold points by $\mathcal{O} $. The embedding of $G$ and the gluing of its vertices to the boundary gives rise to a generalised surface dissection $(S,M,\mathcal{O} ,G)$ where $M$ denotes the set of non-special vertices of $G$ which will refer to as  marked points.

It follows from \cite{AB22, A21, LSV22} in the ungraded case and from \cite{AP24, BSW24, CK24} in the graded case that such a surface dissection corresponds to a (formal) generator of the partially wrapped Fukaya category $\mathcal W(S,M,\mathcal{O} ,\eta)$, and that as generator it gives rise to the graded skew-gentle algebra from which we constructed the graded marked ribbon graph $G$. Following \cite{LP20} the graded marked ribbon graph $G$ induces a line field $\eta$. Note that in \cite{LP20} this construction was given in terms of  ribbon graphs for graded gentle algebras. But every graded skew-gentle algebra is a deformation of a graded gentle algebra induced by a Hochschild 2-cocycle of the gentle algebra corresponding to a loop which squares to zero. These correspond exactly to the squares of the loops $\varepsilon \in \rm{Sp}$.  On the other hand, these 2-cocycles correspond to fully stopped  boundary  components with winding number 1 with respect to the line field $\eta$. The resulting line field on the surface without orbifold points induces the `same' line  field on the  orbifold surface such that each orbifold point has winding number 1. This follows from the fact that every graded skew-gentle algebra is a deformation of a gentle algebra along a 2-cocycle corresponding to loop squaring to zero and where the corresponding boundary component in the surface has winding number one. 
In geometric terms, see \cite{BSW24}, the deformation of the graded (skew-)gentle algebra to a graded skew-gentle algebra by such a 2-cocycle corresponds in the surface to  replacing a fully stopped boundary component with winding number 1 with respect to $\eta$ by an orbifold point which then also has winding number 1 with respect to $\eta$.

Furthermore, it follows from \cite{AB22, A21, LSV22} in the ungraded case and from \cite{AP24, BSW24, CK24} in the graded case that homotopy classes of graded curves connecting boundary components and graded closed curves of winding number zero together with a local system, correspond respectively to string and band objects in the perfect derived category of the graded skew-gentle algebra. 

The Koszul dual $A^!$ of a graded skew-gentle algebra is again a graded skew-gentle algebra and as such $A^!$ has an associated graded ribbon graph $G^!$. In \cite{LSV22}, given a skew-gentle algebra, a dual graph graph $G^*$ is constructed. Namely, let $A$ 
be a skew-gentle algebra with associated graded marked ribbon graph $G$ and let   $(S,M,\mathcal{O}, G)$ be the orbifold surface it generates, that is $G$, when embedded in $S$, is a deformation retract of $S$ with $M$ the set of vertices of $G$ (possibly glued to the boundary as described above) and let $\mathcal{O}$ be the set of orbifold points corresponding to the special loops of $A$. Then $G^*$ (as embedded graph in the surface) is defined as follows: 
in each boundary segment of the dissection $G$, we place exactly one vertex of $G^*$, so that the vertices of $G$ and $G^*$ alternate on the boundary. Furthermore, any unmarked boundary component in $(S,M,\mathcal{O} ,G)$ is replaced by a vertex of $G^*$. 
 Any orbifold point for $G$ is also an orbifold point for $G^*$. The non-special edges of $G$ are in bijection with the non-special edges of $G^*$ where the bijection is given
 by sending an edge $v$ in $G$ to the unique edge in $G^*$ connecting two vertices of $G^*$ and crossing $v$ exactly once.
Every special edge of $G$ corresponds to a special edge in $G^*$ connecting the orbifold point with the unique vertex of $G^*$ such that the edge does not cross any other edge of $G$.
Note that $G$ cuts the surface $S$ into polygons, each of which contains exactly one vertex of $G^*$ either on the boundary or in its interior. Similarly, $G^*$ cuts $S$ into polygons such that each polygon contains exactly one vertex of $G$ either in its boundary or in its interior. 

\begin{Rmk} We note that in terms of the language coming from symplectic geometry and Fukaya categories, the  marked $G$- and $G^*$-points in the interior of the surface correspond to boundary components  as follows. In this paper an unmarked boundary component is replaced by a $G^*$-point in the interior which in terms of surfaces associated to partially wrapped Fukaya categories corresponds to a fully stopped boundary component, see for example \cite{HKK17, BSW24, BSW25}. A  marked point in the interior corresponding to a vertex of $G$  in this paper corresponds to a fully marked, equivalently an unstopped boundary component in \cite{HKK17, BSW24, BSW25}.   
\end{Rmk}

\subsection{Geometric interpretation of the Hochschild cohomology}
We consider both open and closed curves. A curve in $S$ is a
 continuous map $\mathsf{C} : [0,1] \rightarrow S$. We say the curve is \emph{open} if for all $x\in (0,1)$, $\mathsf{C}(x)\in S\setminus \partial S$ and $\mathsf{C}(0)$ and $\mathsf{C}(1)$ are marked points associated to  vertices of $G$. If $\mathsf{C}(0)=\mathsf{C}(1)$ and  for all $x\in [0,1]$, $\mathsf{C}(x)\in S\setminus \partial S$, we call   $\mathsf{C}$ a  \emph{closed} curve.
  We consider a second kind of  open curves. Namely, open curves $\mathsf{C^*}$  whose   endpoints $\mathsf{C}^*(0)$ and $\mathsf{C}^*(1)$ are marked points associated to vertices of $G^*$. 
  Note that when talking about a curve, we always identify it with its image in $S$ and consider it up to homotopy. Furthermore, we always assume that the curves we are working with  are in
 minimal position with respect to $G$ and $G^*$.

\subsubsection{Definition of the combinatorial winding number}\label{Subsubsec:winding number} We now  define the combinatorial winding number $w (\mathsf{C})$ of a curve
 $\mathsf{C} $. A grading on a curve $\mathsf{C} $ is given by a function
 $f_n: \mathsf{C} \cap G^*\rightarrow \mathbb{Z} $, where the set $\mathsf{C} \cap G^*$ is ordered by the direction of $\mathsf{C} $. Let $x_1, x_2 \in [0,1]$ be such that $\mathsf{C} (x_1)$ and $\mathsf{C} (x_2)$ are, respectively,  the first and second
 crossing of $\mathsf{C} $ with $G^*$, and $\gamma_1$ be the segment of $\mathsf{C} $ that lies between $x_1$ and $x_2$. Then $\gamma_1$ lies in a  polygon $P_f$ cut out by $G^*$ and cuts $P_f$ into two subpolygons, where one of them contains a marked point corresponding to a $G$-vertex $v$. Denote by $\alpha _1,\dots,\alpha _m$ the arrows in the subpolygon of $P_f$ which does not contain $v$. Now, the grading $f_n$ of
 $\mathsf{C} $ is defined by assigning the integer $n$ to $\mathsf{C} (x_1)$, then we have 
 \[
f_n(\mathsf{C}(x_2)):=
\begin{cases}
    n+1-\sum_{i=1}^m|\alpha_i| & \text{if } v \text{ lies to the right of } \mathsf{C},  \\ n-1+\sum_{i=1}^m|\alpha_i| & \text{if } v \text{ lies to the left of } \mathsf{C}.
\end{cases}
\]
Propagating this along the whole of $\mathsf{C} $, defines the grading $f_n$.
The \emph{combinatorial winding number} $w(\mathsf{C})$ is defined to be $w(\mathsf{C})=f_n(\mathsf{C}(x_m))-n$
  where the ordered set of $\mathsf{C} \cap G^*$  is given by $\{\mathsf{C} (x_1),\dots,\mathsf{C} (x_m)\}$. If $\mathsf{C} $ is a closed curve, we can fix $x_1$ to be any one of
 the intersections of $\mathsf{C}$ and $G^*$. Without loss of generality we can assume that $x_1=0$ and $x_m =1$. We note, that $w(\mathsf{C})$ is independent of $n$.

For every connected boundary component $B$ in $S$, we denote by $\mathsf{C}_{B}$  a representative  of the homotopy class of curves homotopic to $B$. Furthermore, we refer to  the $G$-vertices in the interior of $S$ as \emph{$G$-punctures} and to the $G^*$-vertices in the interior as \emph{$G^*$-punctures}.  For a $G$-puncture $P$ we denote by $\mathsf{C}_P$  a representative of the homotopy class of simple closed curves around $P$.  We note that $P$ as a vertex of $G$ corresponds to a  cycle $p_C$ without relations in $(Q,I)$ and  the winding number 
of the curve $\mathsf{C}_P$  is given by 
$w(\mathsf{C}_P)=|p_{\mathsf{C}}|$. Similarly, for a $G^*$-puncture $P^*$ we denote by $\mathsf{C}_{P^*}$  a representative of the homotopy class of simple closed curves around $P^*$. Note that $P^*$ as a vertex of $G^*$ corresponds to a  cycle $p_C$ with full relations in $(Q,I)$ and  the winding number 
of the curve $\mathsf{C}_{P^*}$  is given by 
$w(\mathsf{C}_{P^*})=l(p_{\mathsf{C}})-|p_{\mathsf{C}}|$. 
More generally, by the construction of $(S,M,\mathcal{O},G)$, every curve $\mathsf{C}$ corresponds to a sequence  of arrows and formal inverses of arrows in $Q$ and we will denote this sequence  by $p_{\mathsf{C}}$. Conversely, note that,  for every $(\gamma,\alpha)\in \Gamma\|\calB$, we can associate a curve $\mathsf{C}$ to the sequence of arrows and inverse arrows $p_{\mathsf{C}}:=\gamma\alpha^{-1}$ and that  the winding number  $w(\mathsf{C})$ is equal to $l(\gamma)-|\gamma|+|\alpha|-1$.
Note that, for an  element $\alpha=\pi(\alpha_m\cdots\alpha_1)\in \calB$,  any of the $\alpha_i$ might be a concatenation of an arrow $\beta$ with a special loop on one or both ends
, depending on whether $s(\alpha_i)$ and $t(\alpha_i)$ is a special vertex or not. Since the special loops are zero graded, their presence has no bearing on the winding number we just defined.
Finally, we note that it can be checked that for closed curves the combinatorial winding number coincides with the geometric winding number from the line field $\eta$ which is induced by the grading of $A$, see for example \cite{LP20}. 

Recall that $A=\bbK Q/I$ where $I = \langle R \cup \{\varepsilon^2 -\varepsilon \mid \varepsilon \in \Sp\}\rangle$. We now consider the graded gentle algebra $\Lambda \cong \bbK Q_\Lambda/\langle R\rangle$  which is obtained from $A$ by deleting all special loops, that is $\Lambda = A/\langle \Sp \rangle$. In particular,   $(Q_\Lambda)_0=Q_0, (Q_\Lambda)_1=Q_1\setminus \Sp$. Let $(S_{\Lambda},M_\Lambda,G_\Lambda)$ be the graded marked surface obtained from the graded marked  ribbon graph $G_\Lambda$ associated to the graded gentle algebra $\Lambda$. It follows from the constructions, see also \cite{LSV22} for the ungraded case, that  $(S_\Lambda,M_\Lambda,G_\Lambda)$ is obtained from $(S, M, \mathcal{O}, G)$ by replacing a disk around each orbifold point by a smooth disk and so the surface $S_\Lambda$ is a smooth surface. 

 \subsubsection{The fundamental group of a quiver} For a quiver $Q$, we now recall the definition of the fundamental group of the graph underlying  $Q$. Let $\bar Q$ be the double quiver of $Q$,  where $\bar Q_1=Q_1\cup Q_1^{-1}$ and $Q_1^{-1}$ is the set of formal inverses of the arrows in quiver $Q$, that is, for any $a \in Q^{-1}$, the formal inverse $a^{-1}$ of $a$ is a new arrow such that 
$s(a^{-1})=t(a)$ and $t(a^{-1})=s(a)$. 
Now let $v\in Q_0$ be a vertex.
\begin{itemize}
    \item A \emph{path} in $\bar Q$ is $a_n\cdots a_1$ where $a_i\in \bar Q_1,t(a_i)=s(a_{i+1})$, 
    for $1\leq i\leq n-1$, and a \emph{cycle} at $v$ is a path that starts and ends at $v$.
    \item A \emph{spur} is a subpath of the form $aa^{-1}$ where $a\in \bar Q_1$ and $a^{-1}$ is its inverse.
    \item An \emph{elementary homotopy} is the insertion or deletion of a spur in a path.
    \item Two loops are said to be \emph{homotopic} if one can be transformed into the other by a finite sequence of elementary homotopies. 
    \end{itemize}
The set of homotopy classes of cycles with base point $v$, 
under the operation of path concatenation, forms a group. It is called the 
\emph{fundamental group} of $Q$
based at $v$ and denoted by $\pi_1(Q,v)$. Since $Q$ is connected, we can omit the base point $v$ and refer to  the fundamental group of $ Q$ up to isomorphism without referring to its basepoint.

\subsubsection{Geometric interpretation of the Hochschild cohomology} Let $\mathscr{G}$ be the generating set of $\HH^*(A)$
as described in Theorem \ref{Thm: basis-as-alg} and let $\mathcal{F}\subseteq\mathscr{G}$ be a fixed set of derivations arising from the  complement of a chosen  spanning tree $T$ of $Q\setminus \Sp$. Recall from Proposition \ref{1th-graded-HH(sg)} that the elements in $\mathcal{F}$ are elements of $\HH^{1,\sbl}(A)$. 

\begin{Lem}\label{Lem-F}
Let $A$ be the graded skew-gentle algebra associated to a graded skew-gentle triple $(Q,R,\Sp)$ and let $\Lambda=\bbK Q_\Lambda/\langle R\rangle \cong A/\langle \Sp \rangle$ be the associated graded gentle algebra as above with surface model $(S_{\Lambda},M_\Lambda,G_\Lambda)$. Then  the free group generated by the set $\mathcal{F}$ in $\HH^{1,\sbl}(A)$ is isomorphic to the fundamental group $\pi_1(S_\Lambda)$. Furthermore, the rank of $\pi_1(S_\Lambda)$ is $1-|Q_0|+|Q_1|$.
\end{Lem}
\begin{Proof}
    Let $T$ be a spanning tree of $Q$, then $T$ is also  a spanning tree of the quiver $Q_\Lambda$. The free group on the set $(Q_\Lambda)_1\setminus T$ and $\pi_1(Q_\Lambda)$ are isomorphic, see, for example,  \cite[Theorem 1.2.10]{Laz14}. More precisely, if we take an arrow $c\in (Q_\Lambda)_1\setminus T=Q_1\setminus (T\cup \Sp)$, the corresponding element in $\pi_1(Q_\Lambda,v)$ is $T[t(c),v]\cdot c\cdot T[v,s(c)]$, where $T[u,w]$ is the unique shortest path in $T$ from vertex $u$ to vertex $w$. 

On the other hand, by \cite[Proposition 1.22]{OPS} we have $\pi_1(S_\Lambda)\cong\pi_1(Q_\Lambda)$. Therefore, there is an isomorphism between the free group on the set $\mathcal{F}$ in $\HH^{1,\sbl}(A)$ and the fundamental group $\pi_1(S_\Lambda)$.

\end{Proof}

The following theorem gives the geometric interpretation of the Hochschild cohomology in terms of the associated surface, which we also summarise in Table~\ref{table: Hochschild}. For notational purposes we will say that a boundary component contains one marked point if it contains  exactly one marked point corresponding to a vertex of $G$ and hence also one marked point corresponding to a vertex of $G^*$. Also recall from the end of Subsection~\ref{Subsubsec:winding number} the definition of the graded gentle algebra $\Lambda = \bbK Q_\Lambda / \langle R \rangle \cong A/\langle \Sp \rangle$ associated to a graded skew-gentle algebra $A = \bbK Q/I$ where $I = \langle R \cup \{ \varepsilon^2-\varepsilon \mid \varepsilon \in \Sp \}\rangle$ and $Q_\Lambda$ is obtained from the quiver of $Q$ by deleting all the special loops.

\begin{Thm}\label{geom.cohom}
Let $A$ be the graded skew-gentle algebra associated to a graded skew-gentle triple $(Q,R,\Sp)$, and $\Lambda= \bbK Q_\Lambda/\langle R\rangle \cong A/ \langle \Sp \rangle$ the associated graded gentle algebra. Let $(S,M,\mathcal{O},G)$ be the orbifold dissection of $A$ and $(S_{\Lambda},M_\Lambda,G_\Lambda)$ the surface dissection of $\Lambda$, respectively.
Then there is a one to one correspondence 
\[
\{\text{generators of }\HH^*(A)\text{ in }\mathscr{G}\}\;\longleftrightarrow\;
\begin{aligned} 
&\{\text{boundary components with 1 marked point}\text{ in }(S,M,G)\} \\
&\qquad \qquad\quad\quad \cup \;\; \{\text{$G$-punctures} \text{ in }(S,M,G)\}\\
& \qquad \qquad\quad\quad\cup \;\; \{\text{$G^*$-punctures}\text{ in }(S,M,G)\}\\
& \qquad \qquad\quad\quad\quad \cup \;\; \{\text{generators of }\pi_1(S_\Lambda)\},
  \end{aligned}
\]

where the bijection is broken down as follows:  
\begin{itemize}

\item[$(1)$] the pairs $(c,c)\in \HH^{1,0}(A)$ where $c\in Q_1\setminus (T\cup \Sp)$ $\longleftrightarrow$ generators of $\pi_1(S_\Lambda)$,

\item[$(2)$] the pairs $(s(\alpha),\alpha)\in \HH^{0,|\alpha|}(A)$  where $\alpha$ is a $\calB$-maximal path in $(Q,I)$ $\longleftrightarrow$ the boundary components $B$ with 1 marked point in $S$ such that $w(\mathsf{C_B})+1= 0+|\alpha|$,

\item[$(3)$] the pairs $(\gamma,\alpha)\in \HH^{l(\gamma),|\alpha|-|\gamma|}(A)$ 
where $\gamma$ is a $\Gamma$-maximal element of $\Gamma$ and $\gamma$ and $\alpha$ neither begin nor end with the same arrow $\longleftrightarrow$ the boundary components $B$ with 1 marked point in $S$ such that 
$w(\mathsf{C_B})+1=
l(\gamma)+|\alpha|-|\gamma|$,

\item[$(4)$] the sums $\llangle {\alpha}_s\in \HH^{0,|\alpha|}(A)$ where $\alpha$ is a primitive cocomplete cycle and $|\alpha|$ even or $\Char (\bbK)=2$
$\longleftrightarrow$ the $G$-punctures $P$ in $S$ such that  $w(\mathsf{C_P})=|\alpha|$,

\item[$(4')$] the sums $\llangle {\alpha^2}_s\in \HH^{0,2|\alpha|}(A)$ where $\alpha$ is a primitive cocomplete cycle and $|\alpha|$ is odd and $\Char (\bbK)\neq2$  $\longleftrightarrow$ the $G$-punctures $P$ in $S$ such that  $w(\mathsf{C_P})=|\alpha|$,

\item[$(5)$] the sums $\llangle{C}_{gr}\in \HH^{l(C),-|C|}(A)$ where $C$ is a primitive complete cycle in $Q_{l(C)}\setminus \Sp^{l(C)}$ and $l(C)-|C|$ is even or $\Char (\bbK)=2$
$\longleftrightarrow$ 
the $G^*$-punctures $P^*$ such that  $w(\mathsf{C}_{P^*})=l(C)-|C|$,

\item[$(5')$] the sums $\llangle{C^2}_{gr}\in \HH^{2l(C),-2|C|}(A)$ where $C$ is a primitive complete cycle in $Q_{l(C)}\setminus \Sp^{l(C)}$ and $l(C)-|C|$ is odd and $\Char (\bbK)\neq2$
$\longleftrightarrow$
the $G^*$-punctures $P^*$ such that $w(\mathsf{C}_{P^*}) =l(C)-|C|$.
\end{itemize} 
\end{Thm}

\begin{Proof}
    
Given a spanning tree $T$ of the quiver $Q$, let $\mathcal{F}\subseteq\mathscr{G}$ be a fixed set of derivations arising from the complement of $T\cup \Sp$. By Lemma \ref{Lem-F}, there is a one to one correspondence between the set $\mathcal{F}$ and the generators of $\pi_1(S_\Lambda)$.
Thus we only need to show the correspondence between the generators in $\mathscr{G}\setminus \mathcal{F}$ and the set of boundary components with 1 marked point and $G$-punctures and $G^*$-punctures in $(S,M,\mathcal{O},G)$.
    We begin by showing that every generator $g\in \mathscr{G}\setminus \mathcal{F}$ as described in Theorem \ref{Thm: basis-as-alg}
can be naturally associated to either a boundary component with 1 marked point, or a $G$-puncture or a $G^*$-puncture in $(S,M,\mathcal{O},G)$.

Suppose $g=(s(\alpha),\alpha)\in \HH^{0,|\alpha|}(A)$ with $\alpha$ a $\calB$-maximal path, then $\alpha=\alpha_n\cdots\alpha_1$ is a cycle with $\alpha_1\alpha_n=0$ in $A$. Note that since $\alpha$ is a $\calB$-maximal path, no other arrow starts or ends at $s(\alpha)=t(\alpha)$ and, moreover,  $s(\alpha)$ is not a special vertex. 
Moreover, $\alpha$ is also a maximal path in $A'$  thus it will give rise to a finite length $\Sp$-maximal path $\alpha'=\alpha'_m\cdots\alpha'_1$  by deleting all special loops in $\alpha$. By construction of $G$, $\alpha'$ corresponds to a vertex of $G$, and the edge in $G$ which corresponds to the vertex $s(\alpha)=s(\alpha')$ is a loop in $G$. Moreover, the marked pair at the vertex $\alpha'$ in $G$ is the pair $(t(\alpha'_m), s(\alpha'_1))$. Thus, by  construction of $S$, there is a boundary component $B$ that contains only one marked point $\alpha'$. Locally, in $S$ we have the following diagram: 

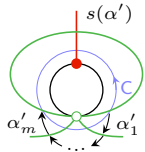
\begin{figure}[H]
    \centering
\begin{tikzpicture}
    \begin{scope}

\draw[line width=.5pt] (0,0) circle(1em);
\draw[draw=blue!50, line width=.3pt] (0,0) circle(1.5em);
\node[font=\scriptsize, color=blue!50] at (3:1.9em) {$\mathsf{C}$};

\draw[-stealth, draw=blue!50, line width=.3pt]
    (0:1.5em) arc[start angle=0, end angle=20, radius=1.5em];

\node[circle,fill=green!10!red,inner sep=0,outer sep=0,minimum size=.4em] at (90:1em) {};

\draw[green!10!red, line width=0.6pt] (90:1em) -- ++(0,2em);
\node[font=\scriptsize, above] at (60:2.5em) {$s(\alpha')$};

\draw[draw=yellow!30!green, line width=.6pt, opacity=0]
  (-190:1.6em) to[out=-90, in=-180] (-90:1em)
   (10:1.6em) to[in=-30, out=-135] (-90:1em);

\draw[draw=yellow!30!green, line width=.6pt]
  (-135:2.5em) to[out=0, in=-120] (-90:1em)
  (-45:2.5em) to[out=-180, in=-60] (-90:1em);

\draw[draw=yellow!30!green, line width=.6pt]
  (-90:1em)
    arc[x radius=2.5em, y radius=1.6em,
        start angle=-90, end angle=270];

\begin{scope}[yshift=-1em]
\draw[-stealth, overlay] (-147:1.25em) arc[start angle=-147, end angle=-180, radius=1.6em];
\draw[-stealth, overlay] (5:1.25em) arc[start angle=-4, end angle=-32, radius=1.65em];
\draw[-stealth, overlay] (-36:1.25em) arc[start angle=-36, end angle=-64, radius=1.25em];
\node[font=\scriptsize, overlay] at (-78:1.25em) {.};
\node[font=\scriptsize, overlay] at (-90:1.25em) {.};
\node[font=\scriptsize, overlay] at (-102:1.25em) {.};
\draw[-stealth, overlay] (-116:1.25em) arc[start angle=-116, end angle=-144, radius=1.25em];
\node[font=\scriptsize, right] at (1.15em,-.3em) {$\alpha_1'$};
\node[font=\scriptsize, left]  at (-1.1em,-.25em) {$\alpha_m'$};
\end{scope}

\node[circle,fill=white,draw=yellow!30!green,inner sep=0,outer sep=0,minimum size=.35em] at (-90:1em) {};

\end{scope}
\end{tikzpicture}
\end{figure}

It follows from the diagram that $g$ naturally corresponds to the boundary curve $\mathsf{C}_B$ with $w(\mathsf{C}_B)=-1+\sum_{i=1}^m|\alpha_i|=|\alpha|-1$ and $p_{\mathsf{C}_B}=\alpha$. 
Note that the total degree of $g$ is exactly $w(\mathsf{C}_B)+1$.

Suppose $g=(\gamma,\alpha)\in \HH^{n,|\alpha|-|\gamma|}(A)$ with $\gamma=\gamma_n\cdots\gamma_1$ a $\Gamma$-maximal element, and $\gamma$ and $\alpha\in \calB$ neither beginning nor ending with the same arrow. 
Since $\alpha$ is a $\calB$-maximal element, it will give rise to a finite  length $\Sp$-maximal path $\alpha'=\alpha'_m\cdots\alpha'_1$. 
By the construction of $G$, $\alpha'$ corresponds to a vertex $v$ of $G$ and the marked pair at $v$ is $(t(\alpha'_m),s(\alpha'_1))$. 
Thus, there is a maximal fan $F$ at the vertex $v$ with m+1 edges on some boundary component of $S$. Analogously, $\gamma$ is a $\Gamma$-maximal element, so it corresponds to a vertex $v^*$ of $G^*$ and there is a maximal fan $F^*$ at the vertex $v^*$ with $n+1$ edges on some boundary component of $S$. Since $\gamma$ and $\alpha$  are parallel paths,
the first edge of $F$ - corresponding to $s(\alpha)=s(\alpha')$- intersects  the first edge of $F^*$ - corresponding to $s(\gamma)$, and the last edge of $F$ - corresponding to $t(\alpha)=t(\alpha')$- intersects the last edge of $F^*$ - corresponding to $t(\gamma)$. By the construction of $(S,M,\mathcal{O},G)$, the vertices $v$ and $v^*$ are on the same boundary component $B$ and there are no other marked points on $B$. The boundary curve $\mathsf{C}_B$ is exactly $\mathsf{C}_{\gamma\alpha^{-1}}$, with $w(\mathsf{C}_{\gamma\alpha^{-1}})=n-1+|\alpha|-|\gamma|$ and $p_{\mathsf{C}_{\gamma\alpha^{-1}}}=\gamma\alpha^{-1}$. Note that the total degree of $g$ is $w(\mathsf{C}_B)+1$. Locally in $S$ we have the following diagram:
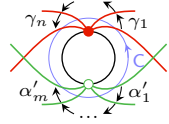
\begin{figure}[H]
    \centering
\begin{tikzpicture}
\draw[line width=.5pt] (0,0) circle(1em);
\draw[draw=blue!50, line width=.3pt] (0,0) circle(1.5em);
\node[font=\scriptsize, color=blue!50] at (-3:1.9em) {$\mathsf{C}$};

\draw[-stealth, draw=blue!50, line width=.3pt]
    (0:1.5em) arc[start angle=0, end angle=20, radius=1.5em];

\node[circle,fill=green!10!red,inner sep=0,outer sep=0,minimum size=.4em] at (90:1em) {};

\draw[draw=green!10!red, line width=.6pt]
  (190:3em) to[out=45, in=150] (90:1em)
  (135:2.5em) to[out=0, in=120] (90:1em)
  (-10:3em) to[in=30, out=135] (90:1em)
  (45:2.5em) to[out=180, in=60] (90:1em);

\draw[draw=yellow!30!green, line width=.6pt]
  (-190:3em) to[out=-45, in=-150] (-90:1em)
  (-135:2.5em) to[out=0, in=-120] (-90:1em)
  (10:3em) to[in=-30, out=-135] (-90:1em)
  (-45:2.5em) to[out=-180, in=-60] (-90:1em);

\begin{scope}[yshift=1em]
\draw[-stealth, overlay] (147:1.25em) arc[start angle=147, end angle=175, radius=1.25em];
\draw[-stealth, overlay] (4:1.25em) arc[start angle=4, end angle=32, radius=1.25em];
\draw[-stealth, overlay] (36:1.25em) arc[start angle=36, end angle=64, radius=1.25em];
\node[font=\scriptsize, overlay] at (78:1.25em) {.};
\node[font=\scriptsize, overlay] at (90:1.25em) {.};
\node[font=\scriptsize, overlay] at (102:1.25em) {.};
\draw[-stealth, overlay] (116:1.25em) arc[start angle=116, end angle=144, radius=1.25em];
\node[font=\scriptsize, right] at (1.15em,.3em) {$\gamma_1$};
\node[font=\scriptsize, left]  at (-1.15em,.4em) {$\gamma_n$};
\end{scope}

\begin{scope}[yshift=-1em]
\draw[-stealth, overlay] (-147:1.25em) arc[start angle=-147, end angle=-175, radius=1.25em];
\draw[-stealth, overlay] (-4:1.25em) arc[start angle=-4, end angle=-32, radius=1.25em];
\draw[-stealth, overlay] (-36:1.25em) arc[start angle=-36, end angle=-64, radius=1.25em];
\node[font=\scriptsize, overlay] at (-78:1.25em) {.};
\node[font=\scriptsize, overlay] at (-90:1.25em) {.};
\node[font=\scriptsize, overlay] at (-102:1.25em) {.};
\draw[-stealth, overlay] (-116:1.25em) arc[start angle=-116, end angle=-144, radius=1.25em];
\node[font=\scriptsize, right] at (1.15em,-.3em) {$\alpha_1'$};
\node[font=\scriptsize, left]  at (-1.1em,-.25em) {$\alpha_m'$};
\end{scope}

\node[circle,fill=white,draw=yellow!30!green,inner sep=0,outer sep=0,minimum size=.35em] at (-90:1em) {};
\end{tikzpicture}
\end{figure}

Given $g=\llangle {\alpha}_s\in \HH^{0,|\alpha|}(A)$ with $\alpha\in \oC^{\basic}(\calB)$, the cycle  $\alpha$  gives rise to an 
$\Sp$-maximal path $\alpha'=\alpha_m'\cdots \alpha_1'$, which corresponds to a $G$-puncture $P$.
If $\alpha$ is a primitive cycle, then we have $p_{\mathsf{C}_{P}}=\alpha$, $w(\mathsf{C}_{P})=|
\alpha|$.  Otherwise, $\alpha$ is the square of a primitive cycle $D$ and the path traced out by $D$ corresponds to the boundary curve $\mathsf{C}_{P}$, that is $p_{\mathsf{C}_{P}}=D$. In this situation $\alpha$ corresponds to $\mathsf{C}_{P}^2$ with $p_{\mathsf{C}^2_{P}}=
\alpha$ and $w(\mathsf{C}^2_{P})=|\alpha|$. Locally in $S$ we have the following diagram:
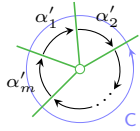
\begin{figure}[H]
    \centering
\begin{tikzpicture}[scale=1.5]
\begin{scope}
\draw[draw=blue!50, line width=.3pt] (0,0) circle(1.35em);
\node[font=\scriptsize, color=blue!50] at (-45:1.8em) {$\mathsf{C}$};

\draw[-stealth, draw=blue!50, line width=.3pt]
    (0:1.35em) arc[start angle=0, end angle=20, radius=1.5em];

\draw[draw=yellow!30!green, line width=.6pt] (0,0) to (180-20:1.7em);
\draw[draw=yellow!30!green, line width=.6pt] (0,0) to (180-85:1.7em);
\draw[draw=yellow!30!green, line width=.6pt] (0,0) to (180-150:1.7em);
\draw[draw=yellow!30!green, line width=.6pt] (0,0) to (180+45:1.7em);

\foreach \a in {110, 45, -20, -85, -150} {
    \draw[-stealth] (180+\a-2:1em) arc[start angle=180+\a-2, end angle=180+\a-61, radius=1em];
}

\node[font=\scriptsize] at (-51.5:1em) {.};
\node[font=\scriptsize] at (-51.5+12:1em) {.};
\node[font=\scriptsize] at (-51.5-12:1em) {.};
\node[font=\scriptsize] at (180-20+32.5:1.5em) {$\alpha_m'$};
\node[font=\scriptsize] at (180-85+32.5-5:1.55em) {$\alpha_1'$};
\node[font=\scriptsize] at (180-150+32.5:1.5em) {$\alpha_2'$};

\node[circle,fill=white,draw=yellow!30!green,inner sep=0,outer sep=0,minimum size=.35em] at (0,0) {};
\end{scope}
\end{tikzpicture}
\end{figure} 

Now suppose $g=\llangle{C}_{gr}\in \HH^{n,-|C|}(A)$ with $C=c_n\cdots c_1\in \oC^{\basic}(\Gamma)\setminus \Sp^{\basic}.$ Here,  $C$ gives rise to a $G^*$-puncture $P^*$ in $S$. If $C$ is a primitive cycle, then we have $p_{\mathsf{C}_{P^*}}=C$, $w(\mathsf{C}_{P^*})=n-|C|$. Otherwise, $C$ is the square of a primitive cycle $D$ and the path traced out by $D$ corresponds to the curve $\mathsf{C}_{P^*}$, that is $p_{\mathsf{C}_{P^*}}=D$ and $p_{\mathsf{C}^2_{P^*}}=C$ with $w(\mathsf{C}^2_{P^*})=n-|C|$. Locally in $S$ we have the following diagram:
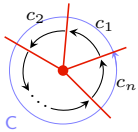
\begin{figure}[H]
    \centering
\begin{tikzpicture}[scale=1.5]
\begin{scope}

\draw[draw=blue!50, line width=.3pt] (0,0) circle(1.35em);
\node[font=\scriptsize, color=blue!50] at (225:1.85em) {$\mathsf{C}$};

\draw[-stealth, draw=blue!50, line width=.3pt]
    (0:1.35em) arc[start angle=0, end angle=20, radius=1.5em];

\node[circle,fill=green!10!red,inner sep=0,outer sep=0,minimum size=.4em] at (0,0) {};

\draw[draw=green!10!red, line width=.6pt] (0,0) to (20:1.7em);
\draw[draw=green!10!red, line width=.6pt] (0,0) to (85:1.7em);
\draw[draw=green!10!red, line width=.6pt] (0,0) to (150:1.7em);
\draw[draw=green!10!red, line width=.6pt] (0,0) to (-45:1.7em);

\foreach \a in {-110, -45, 20, 85, 150} {
    \draw[-stealth] (\a+2:1em) arc[start angle=\a+2, end angle=\a+61, radius=1em];
}

\node[font=\scriptsize] at (51.5+180:1em) {.};
\node[font=\scriptsize] at (51.5+180+12:1em) {.};
\node[font=\scriptsize] at (51.5+180-12:1em) {.};
\node[font=\scriptsize] at (20-32.5:1.65em) {$c_n$};
\node[font=\scriptsize] at (85-32.5-5:1.55em) {$c_1$};
\node[font=\scriptsize] at (150-32.5:1.5em) {$c_2$};

\end{scope}
\end{tikzpicture}
\end{figure}

Conversely, suppose that $B$ is a boundary component with one marked point corresponding to a vertex $v$ of $G$. By the construction of $S$, this boundary component also contains exactly one vertex $v^*$ of $G^*$. The maximal fan $F$ in $G$ at $v$ corresponds to a finite length $\Sp$-maximal path, which gives rise to a $\calB$-maximal path $\alpha$ in $A$. Similarly, the maximal fan $F^*$ in $G^*$ at $v^*$ corresponds to a $\Gamma$-maximal path $\gamma$ where there are $l(\gamma)+1$ edges of $F^*$. Then the boundary curve $\mathsf{C}_B$ is homotopic to the curve $\mathsf{C}_{\gamma\alpha^{-1}}$, and $p_{\mathsf{C}_B}=\gamma\alpha^{-1}$, $w(\mathsf{C}_B)=l(\gamma)-1+|\alpha|-|\gamma|$. We see that $g=(\gamma,\alpha)$ is a generator in $\mathscr{G}\setminus \mathcal{F} $ that belongs to $\HH^{l(\gamma),|\alpha|-|\gamma| }(A)$. In particular, we note that when $l(\gamma)=0,$ the maximal fan $F^*$ might contain only a single edge $s(\alpha)$.

 Let $P$ be a $G$-puncture in $S$, by definition of $G$, $P$ corresponds to a rotation class of a primitive cycle without relations in $(Q,R,\Sp)$. Thus $p_{\mathsf{C}_P}$ is a primitive cyclic path in $\calB$ and $w(\mathsf{C}_P)=|p_{\mathsf{C}_P}|$. By Theorem \ref{Thm: basis-as-alg}, if $w(\mathsf{C}_P)$ is even or $\Char (\bbK)=2$ then $\llangle{p_{\mathsf{C}_P}}_s$ is a generator in $\HH^{0,|p_{\mathsf{C}_P}|}(A)$, otherwise $\llangle{p^2_{\mathsf{C}_P}}_s$ is a generator in $\HH^{0,2|p_{\mathsf{C}_P}|}(A)$.

Finally, consider a 
$G^*$-puncture $P^*$ in $S$ which corresponds to a boundary component $B$ with no marked points. Assume the valency of $P^*$ is $n$, then the curve 
 $\mathsf{C}_{P^*}$ corresponds to a complete cycle $p_{\mathsf{C}_{P^*}}$ in $\Gamma\setminus \Sp^n$ such that the winding number $w(\mathsf{C}_{P^*})$ equals $n-|p_{\mathsf{C}_{P^*}}|$. Thus, if $w(\mathsf{C}_{P^*})$ is even or $\Char (\bbK)=2$ then $\llangle{p_{\mathsf{C}_{P^*}}}_{gr}$ is a generator in $\HH^{n,-|p_{\mathsf{C}_{P^*}}|}(A)$, otherwise $\llangle{p^2_{\mathsf{C}_{P^*}}}_{gr}$ is a generator in $\HH^{2n,-2|p_{\mathsf{C}_{P^*}}|}(A)$.
\end{Proof}

\medskip

\begin{table}[H]
    \centering
\begin{tikzpicture}[x=1em, y=1em]

\def\WL{12} 
\def\WR{30} 

\def\Ha{2}  
\def\Hb{8}  
\def\Hc{7}  
\def\Hd{7}  
\def\He{7}   

\def\Yb{-\Ha}
\def\Yc{-\Ha - \Hb}
\def\Yd{-\Ha - \Hb - \Hc}
\def\Ye{-\Ha - \Hb - \Hc - \Hd}
\def\Yf{-\Ha - \Hb - \Hc - \Hd - \He} 

\draw[line width=.4pt] (0,0) rectangle (\WL+\WR, \Yf);

\draw[line width=.4pt] (\WL,0) -- (\WL, \Yf);

\draw[line width=.4pt] (0,\Yb) -- (\WL+\WR,\Yb); 
\draw[line width=.4pt] (0,\Yc) -- (\WL+\WR,\Yc);
\draw[line width=.4pt] (0,\Yd) -- (\WL+\WR,\Yd);
\draw[line width=.4pt] (0,\Ye) -- (\WL+\WR,\Ye);

\node at (\WL/2, -\Ha/2) {Boundary components};
\node at (\WL + \WR/2, -\Ha/2) {Generators of $\HH^*(A)$};

\begin{scope}[shift={(\WL/2, \Yb - \Hb/2)}]

\draw[line width=.5pt] (0,0) circle(1em);
\draw[draw=black!20, line width=.3pt] (0,0) circle(1.4em);
\node[font=\scriptsize, color=black!20] at (3:1.9em) {$\mathsf{C}$};

\node[circle,fill=green!10!red,inner sep=0,outer sep=0,minimum size=.4em] at (90:1em) {};

\draw[green!10!red, line width=0.6pt] (90:1em) -- ++(0,2em);
\node[font=\scriptsize, above] at (60:2.5em) {$s(\alpha')$};

\draw[draw=yellow!30!green, line width=.6pt, opacity=0]
  (-190:1.6em) to[out=-90, in=-180] (-90:1em)
   (10:1.6em) to[in=-30, out=-135] (-90:1em);

\draw[draw=yellow!30!green, line width=.6pt]
  (-135:2.5em) to[out=0, in=-120] (-90:1em)
  (-45:2.5em) to[out=-180, in=-60] (-90:1em);

\draw[draw=yellow!30!green, line width=.6pt]
  (-90:1em)
    arc[x radius=2.5em, y radius=1.6em,
        start angle=-90, end angle=270];

\begin{scope}[yshift=-1em]
\draw[-stealth, overlay] (-147:1.25em) arc[start angle=-147, end angle=-180, radius=1.6em];
\draw[-stealth, overlay] (5:1.25em) arc[start angle=-4, end angle=-32, radius=1.65em];
\draw[-stealth, overlay] (-36:1.25em) arc[start angle=-36, end angle=-64, radius=1.25em];
\node[font=\scriptsize, overlay] at (-78:1.25em) {.};
\node[font=\scriptsize, overlay] at (-90:1.25em) {.};
\node[font=\scriptsize, overlay] at (-102:1.25em) {.};
\draw[-stealth, overlay] (-116:1.25em) arc[start angle=-116, end angle=-144, radius=1.25em];
\node[font=\scriptsize, right] at (1.15em,-.3em) {$\alpha_1'$};
\node[font=\scriptsize, left]  at (-1.1em,-.25em) {$\alpha_m'$};
\end{scope}

\node[circle,fill=white,draw=yellow!30!green,inner sep=0,outer sep=0,minimum size=.35em] at (-90:1em) {};

\end{scope}

\node[align=center] at (\WL + \WR/2, \Yb - \Hb/2) {
   $(s(\alpha),\alpha)\in\HH^{0,|\alpha|}(A)\subseteq \HH^{ w_{\eta}(\mathsf{C})+1}(A)$, \\
   where $\alpha$ is a $\calB$-maximal path and $\alpha'=\alpha_m'\cdots\alpha_1'$  \\ the corresponding  $\Sp$-maximal path.
};

\begin{scope}[shift={(\WL/2, \Yc - \Hc/2)}]

\draw[line width=.5pt] (0,0) circle(1em);
\draw[draw=black!20, line width=.3pt] (0,0) circle(1.5em);
\node[font=\scriptsize, color=black!20] at (-3:1.9em) {$\mathsf{C}$};

\node[circle,fill=green!10!red,inner sep=0,outer sep=0,minimum size=.4em] at (90:1em) {};

\draw[draw=green!10!red, line width=.6pt]
  (190:3em) to[out=45, in=150] (90:1em)
  (135:2.5em) to[out=0, in=120] (90:1em)
  (-10:3em) to[in=30, out=135] (90:1em)
  (45:2.5em) to[out=180, in=60] (90:1em);

\draw[draw=yellow!30!green, line width=.6pt]
  (-190:3em) to[out=-45, in=-150] (-90:1em)
  (-135:2.5em) to[out=0, in=-120] (-90:1em)
  (10:3em) to[in=-30, out=-135] (-90:1em)
  (-45:2.5em) to[out=-180, in=-60] (-90:1em);

\begin{scope}[yshift=1em]
\draw[-stealth, overlay] (147:1.25em) arc[start angle=147, end angle=175, radius=1.25em];
\draw[-stealth, overlay] (4:1.25em) arc[start angle=4, end angle=32, radius=1.25em];
\draw[-stealth, overlay] (36:1.25em) arc[start angle=36, end angle=64, radius=1.25em];
\node[font=\scriptsize, overlay] at (78:1.25em) {.};
\node[font=\scriptsize, overlay] at (90:1.25em) {.};
\node[font=\scriptsize, overlay] at (102:1.25em) {.};
\draw[-stealth, overlay] (116:1.25em) arc[start angle=116, end angle=144, radius=1.25em];
\node[font=\scriptsize, right] at (1.15em,.3em) {$\gamma_1$};
\node[font=\scriptsize, left]  at (-1.15em,.4em) {$\gamma_n$};
\end{scope}

\begin{scope}[yshift=-1em]
\draw[-stealth, overlay] (-147:1.25em) arc[start angle=-147, end angle=-175, radius=1.25em];
\draw[-stealth, overlay] (-4:1.25em) arc[start angle=-4, end angle=-32, radius=1.25em];
\draw[-stealth, overlay] (-36:1.25em) arc[start angle=-36, end angle=-64, radius=1.25em];
\node[font=\scriptsize, overlay] at (-78:1.25em) {.};
\node[font=\scriptsize, overlay] at (-90:1.25em) {.};
\node[font=\scriptsize, overlay] at (-102:1.25em) {.};
\draw[-stealth, overlay] (-116:1.25em) arc[start angle=-116, end angle=-144, radius=1.25em];
\node[font=\scriptsize, right] at (1.15em,-.3em) {$\alpha_1'$};
\node[font=\scriptsize, left]  at (-1.1em,-.25em) {$\alpha_m'$};
\end{scope}

\node[circle,fill=white,draw=yellow!30!green,inner sep=0,outer sep=0,minimum size=.35em] at (-90:1em) {};
\end{scope}

\node[align=center] at (\WL + \WR/2, \Yc - \Hc/2) {
 $(\gamma,\alpha)\in\HH^{n,|\alpha|-|\gamma|}(A)\subseteq \HH^{w_{\eta}(\mathsf{C})+1}(A)$,   \\
where $\gamma$ is a $\Gamma$-maximal path,  $\gamma$ and $\alpha$ \\
neither begin nor end with the same arrow, and \\ $\alpha'=\alpha_m'\cdots\alpha_1'$ is the  $\Sp$-maximal path corresponding to $\alpha$. 
};

\begin{scope}[shift={(\WL/2, \Yd - \Hd/2)}, scale=1.4]

\draw[draw=black!20, line width=.3pt] (0,0) circle(1.35em);
\node[font=\scriptsize, color=black!20] at (-45:1.8em) {$\mathsf{C}$};
\draw[draw=yellow!30!green, line width=.6pt] (0,0) to (180-20:1.5em);
\draw[draw=yellow!30!green, line width=.6pt] (0,0) to (180-85:1.5em);
\draw[draw=yellow!30!green, line width=.6pt] (0,0) to (180-150:1.5em);
\draw[draw=yellow!30!green, line width=.6pt] (0,0) to (180+45:1.5em);

\foreach \a in {110, 45, -20, -85, -150} {
    \draw[-stealth] (180+\a-2:1em) arc[start angle=180+\a-2, end angle=180+\a-61, radius=1em];
}

\node[font=\scriptsize] at (-51.5:1em) {.};
\node[font=\scriptsize] at (-51.5+12:1em) {.};
\node[font=\scriptsize] at (-51.5-12:1em) {.};
\node[font=\scriptsize] at (180-20+32.5:1.5em) {$\alpha_m'$};
\node[font=\scriptsize] at (180-85+32.5-5:1.55em) {$\alpha_1'$};
\node[font=\scriptsize] at (180-150+32.5:1.5em) {$\alpha_2'$};

\node[circle,fill=white,draw=yellow!30!green,inner sep=0,outer sep=0,minimum size=.35em] at (0,0) {};
\end{scope}

\node[align=center] at (\WL + \WR/2, \Yd - \Hd/2) 
{
$\llangle{\alpha^i}_s\in \HH^{0,i|\alpha|}(A)\subseteq
\HH^{i w_{\eta}(\mathsf{C})}(A)$,\\
   where $\alpha$ is a primitive cocomplete cycle and
\\$\alpha'=\alpha_m'\cdots\alpha_1'$  is the corresponding Sp-maximal path 
   \\
  with $i=1$ if $|\alpha|$ is even  or $\Char(\bbK)=2$, and $i=2$ otherwise. 
  };

\begin{scope}[shift={(\WL/2, \Ye - \He/2)}, scale=1.4] 

\draw[draw=black!20, line width=.3pt] (0,0) circle(1.35em);
\node[font=\scriptsize, color=black!20] at (225:1.85em) {$\mathsf{C}$};
\node[circle,fill=green!10!red,inner sep=0,outer sep=0,minimum size=.4em] at (0,0) {};

\draw[draw=green!10!red, line width=.6pt] (0,0) to (20:1.5em);
\draw[draw=green!10!red, line width=.6pt] (0,0) to (85:1.5em);
\draw[draw=green!10!red, line width=.6pt] (0,0) to (150:1.5em);
\draw[draw=green!10!red, line width=.6pt] (0,0) to (-45:1.5em);

\foreach \a in {-110, -45, 20, 85, 150} {
    \draw[-stealth] (\a+2:1em) arc[start angle=\a+2, end angle=\a+61, radius=1em];
}

\node[font=\scriptsize] at (51.5+180:1em) {.};
\node[font=\scriptsize] at (51.5+180+12:1em) {.};
\node[font=\scriptsize] at (51.5+180-12:1em) {.};
\node[font=\scriptsize] at (20-32.5:1.65em) {$c_n$};
\node[font=\scriptsize] at (85-32.5-5:1.55em) {$c_1$};
\node[font=\scriptsize] at (150-32.5:1.5em) {$c_2$};

\end{scope}

\node[align=center] at (\WL + \WR/2, \Ye - \He/2) {
 $\llangle{C^i}_{gr}\in \HH^{in,-i|C|}(A)\subseteq \HH^{i w_{\eta}(\mathsf{C})}(A)$,  \\
  where $C=c_n\cdots c_1$ is a primitive complete cycle in $\Gamma_n\setminus \Sp^n$ \\
   with $i=1$ if $n-|C|$ is even or $\Char(\bbK)=2$, and $i=2$ otherwise.
};

\end{tikzpicture}

\medskip

\caption{Correspondence of the generators of $\HH^*(A)$ with boundary components with 1 marked point, $G^*$-punctures and $G$-punctures. We note that in terms of the partially wrapped Fukaya category of an orbifold surface with stops as considered, for example in \cite{AP24, CK24, BSW24}, a $G$-puncture (as in the third row of the table) corresponds to a boundary component with no stops and a $G^*$-puncture (as in the fourth row of the table) corresponds to a fully stopped boundary component.} 
\label{table: Hochschild}
\end{table}

 In order to interpret the cup product and the Gerstenhaber bracket in the  Hochschild cohomology of a graded skew-gentle algebra $A$, we need to give another interpretation of the generators in $\mathcal{F}$. From the construction of $G$ and $G^*$ on $(S,M,\mathcal{O},G)$, we immediately see that an element $(c,c)$ in $\mathcal{F} $ corresponds to a pair of curves $(\mathsf{C}_{c},\mathsf{C}^*_{c})$ where $p_{\mathsf{C}_{c}}=c=p_{\mathsf{C}^*_{c}}$, and that $w(\mathsf{C}_{c})=1-|c|$.

\subsection{Geometric interpretation of the cup product}
Note that any basis element in Theorem \ref{Thm: coho-basis-sg} that is not a generator is a cup product of generators.
It follows from Table \ref{table: cup product}, that there are exactly four  non trivial cup products between the generators of $\HH^*(A)$:
\begin{itemize}
    \item [(i)] $\llangle{\alpha}_s\smile \llangle{\alpha}_s=\llangle{\alpha^2}_s$ with $\alpha\in \oC^{\basic}(\calB)$. The curve associated to $\llangle{\alpha}_s$ is a primitive cycle or a square of a primitive cycle around a $G$-puncture.
Thus the cup product can be described by a curve obtained from 
concatenating the curves of the corresponding generators,

\item [(ii)] $\llangle{C}_{gr}\smile \llangle{C}_{gr}=(-1)^{l(C)}\llangle{C^2}_{gr}$ with $C\in \oC^{\basic}(\Gamma)\setminus \Sp^{\basic}$.
The curve associated to $\llangle{C}_{gr}$ is a primitive cycle or a square of a primitive cycle around a $G^*$-puncture, thus the interpretation of the cup product is
similar to case (i),

\item[(iii)]$(c,c)\smile \llangle{\alpha}_s=(c,c\alpha)$ with $(c,c)\in \mathcal{F}$, $\alpha\in \oC^{\basic}(\calB)$, and $c$ appears in $\alpha$. The generator $(c,c)$ corresponds to a pair of curves $(\mathsf{C}_{c},\mathsf{C}^*_{c})$ with $p_{\mathsf{C}_{c}}=c=p_{\mathsf{C}^*_{c}}$, and $\llangle{\alpha}_s$ corresponds to a cycle $\mathsf{C}_P$ around a $G$-puncture $P$ with $p_{\mathsf{C}_P}=\alpha$,
so the element $(c,c\alpha)$ corresponds to the pair of curves $(\mathsf{C}_c, \mathsf{C}_c^*\cdot \mathsf{C}_P )$ where $\mathsf{C}_c^*\cdot \mathsf{C}_P$ is obtained by the concatenation of the curves $\mathsf{C}^*_c$ and $\mathsf{C}_P$ at their intersection point,

\item[(iv)]$(c,c)\smile \llangle{C}_{gr}=(-1)^{|c|\cdot |C|}(cC,c)$ with $(c,c)\in \mathcal{F}$, $C\in \oC^{\basic}(\Gamma)\setminus \Sp^{\basic}$, and $c$ appears in $C$. The generator $\llangle{C}_{gr}$ corresponds to a cycle $\mathsf{C}_{P^*}$ around a $G^*$-puncture $P^*$ with $p_{\mathsf{C}_{P^*}}=C$,
and the element $(cC,c)$ corresponds to the pair of curves $(\mathsf{C}_c\cdot \mathsf{C}_{P^*}, \mathsf{C}_c^* )$ where $\mathsf{C}_c\cdot \mathsf{C}_{P^*}$ is  
obtained as the concatenation of 
 the curves $\mathsf{C}_c$ and $\mathsf{C}_{P^*}$ at their intersection point.
\end{itemize}

\subsection{Geometric interpretation of the Gerstenhaber bracket}
Recall from
Theorem \ref{thm: bracket of HH} that the Gerstenhaber bracket of the generators of the Hochschild
cohomology is almost always zero and that the only non trivial ones arise from
brackets of the form $[\mathcal{F},\mathcal{G}\setminus \mathcal{F}]$.
 More precisely, the bracket 
 $[(c,c),v]$
  for $(c, c)\in \mathcal{F}$
and $v=(\gamma,\alpha)\in \mathcal{G}\setminus \mathcal{F}$ is non zero if $\deg_c
(v)$ is non zero, that is $c$ appears in $v$ and in this
case, $[(c, c), v] = \deg_c
(v)v=(\deg_c(\alpha)-\deg_c(\gamma))v$. In terms of the geometric model, suppose that $(c, c)$ is
given by a pair of open curves $(\mathsf{C}_c,\mathsf{C}_c^*)$, and the closed curve associated to $v$ is $\mathsf{C}_v$. The fact that $\deg_c(\alpha)$ is non zero corresponds to the curves $\mathsf{C}_c^*$ and $\mathsf{C}_v$ running parallel between two consecutive edges in a fan in the ribbon graph $G$, and $\deg_c(\alpha)$ counts the number of times the curves run parallel in such a way.
Similarly, $\deg_c(\gamma)$ being non zero corresponds to the curves $\mathsf{C}_c$ and $\mathsf{C}_v$ running parallel between two consecutive edges in a fan in the dual graph $G^*$, and $\deg_c(\gamma)$ counts the number of times this happens.

\section*{ Acknowledgements}
The authors would like to thank Severin Barmeier and Kai Wang for many helpful discussions. The first and the last authors are supported by the China Scholarship Council (CSC). The
second author is supported by the DFG through the project SFB/TRR 191 (No. 281071066-TRR 191). The fourth author is supported by the Outstanding Doctoral Students Overseas Study Program of the University of Science and Technology of China. Part of this work was completed while  the first, fourth, and last authors were visiting the  University of Cologne. They sincerely thank the University of Cologne for its hospitality.



\begin{thebibliography}{88}

\bibitem{Abb15}
H.~Abbaspour, \emph{On algebraic structures of the Hochschild complex}, in: Free loop spaces in geometry and topology, IRMA Lec. Math. Theor. Phys. \textbf{24}, Eur. Math. Soc., Zürich, (2015), 165-222.

\bibitem{A21}
C. Amiot, \emph{Indecomposable objects in the derived category of a skew-gentle algebra using orbifolds}, Representations of algebras and related structures, EMS Ser. Congr. Rep, 2021: 1-24.

\bibitem{AB22}
C. Amiot, T. Br\"ustle,  
\emph{Derived equivalences between skew-gentle algebras using orbifolds}, Doc. Math., \textbf{27} (2022), 933-982.

\bibitem{AP24}
C.~Amiot, P.~G.~Plamondon, \emph{Skew-group $A_{\infty}$-categories as Fukaya categories of orbifolds}, arXiv: 2405.15466, (2024).

\bibitem{BSW24}
S. Barmeier, S. Schroll, Z.~F.~Wang, \emph{Partially wrapped Fukaya categories of orbifold surfaces}, arXiv: 2407.16358, (2024).

\bibitem{BSW25}
S.~Barmeier, S.~Schroll, Z.~F.~Wang, \emph{Deformations of partially wrapped Fukaya categories of surfaces}, arXiv: 2512.16354, (2025).

\bibitem{CSS20}
C.~Chaparro, S.~Schroll and A.~Solotar, \emph{On the Lie algebra structure of the first Hochschild cohomology of gentle algebras and Brauer graph algebras}, J. Algebra, \textbf{558} (2020), 293--326, Special Issue in honor of Michel Broué.

\bibitem{CSSS24} 
C.~Chaparro, S.~Schroll, A.~Solotar and M.~Suárez-Álvarez, \emph{The Hochschild cohomology and the Tamarkin--Tsygan calculus of gentle algebras}, arXiv: 2311.08003, (2023).

\bibitem{CE56}
H.~Cartan, S.~Eilenberg, \emph{Homological algebra}, Princeton Univ. Press, Princeton, NJ, (1956) MR 0077480 (17,1040e).

\bibitem{CK24}
C.~H.~Cho, K.~Kim, \emph{Topological Fukaya category of tagged arcs}, arXiv: 2404.10294, (2024). 

\bibitem{CS2015} 
S.~Chouhy, A.~Solotar, \emph{Projective resolutions of associative algebras and ambiguities}, J. Algebra, \textbf{432} (2015), 22-61.

\bibitem{C19} 
S.~Chouhy, \emph{On geometric degenerations and Gerstenhaber formal deformations}, Bull. London Math. Soc., \textbf{51} (2019), 787-797.

\bibitem{Cib90} 
C.~Cibils, \emph{Rigidity of truncated quiver algebras}, Adv. Math., \textbf{79} (1) (1990), 18-42.

\bibitem{FMT02}
Y.~Félix, L.~Menichi, and J.~C.~Thomas, \emph{Duality in Gerstenhaber algebras}, arXiv:math/0211229 [math.AT], (2002).

\bibitem{Ganatra1}
S.~Ganatra, \emph{Symplectic cohomology and duality for the wrapped Fukaya category}, 
Thesis (Ph.D.)–Massachusetts Institute of Technology
ProQuest LLC, Ann Arbor, MI, (2012). 

\bibitem{Ganatra2}
 S. Ganatra, \emph Symplectic cohomology from Hochschild (co)homology, available at \texttt{https://sheelganatra.com/wrapcy2.pdf}.

\bibitem{GP99}
C.~Geiß, J.~A.~de la Pe\~{n}a, \emph{Auslander--Reiten components for clans}, Bol. Soc. Mat. Mex., \textbf{5} (1999), 307-326.

\bibitem{GP95}
C.~Geiß, J.~A.~de la Pe\~{n}a, \textit{On the deformation theory of finite-dimensional algebras}, Manuscripta Math., \textbf{88} (2) (1995), 191-208.

\bibitem{G63}
M. Gerstenhaber, \emph{The homology structure of an associative ring}, Ann. of Math., \textbf{4} (1963), 267-288.

\bibitem{HKK17}
F.~Haiden, L.~Katzarkov and M.~Kontsevich, \emph{Flat surfaces and stability structures}, Publ. Math. Inst. Hautes \'Etutes Sci., \textbf{126} (2017), 246-318.

\bibitem{HLW23} 
Y.~Han, X.~Liu and K.~Wang, \emph{Hochschild (co)homologies of DG $K$-algebras and their koszul duals}, Front. Math., \textbf{18} (5) (2023), 1113-1155.

\bibitem{H45}
G.~Hochschild, \emph{On the cohomology groups of an associative algebra}, Ann. of Math., \textbf{2} (46)
(1945), 58-67.

\bibitem{Ri91}
J.~Rickard, \textit{Derived equivalences as derived functors}, J. London Math. Soc., \textbf{43} (1991), 37-48.

\bibitem{LSV22}
D. Labardini-Fragoso, S. Schroll, Y. Valdivieso, \emph{Derived categories of skew-gentle algebras and orbifolds}, Glasg. Math. J., \textbf{64} (3) (2022), 649-674.

\bibitem{Laz14}
F. Lazarus, \emph{Combinatorial graphs and surfaces from the
computational and topological viewpoint,} Habilitation  thesis, Université de Grenoble,  2014, available at https://pagesperso.g-scop.grenoble-inp.fr/~lazarusf/Documents/hdr-Lazarus.pdf.

\bibitem{LP20}
Y. Lekili, A. Polishchuk, \emph{Derived equivalences of gentle algebras via Fukaya categories}, Math. Ann., \textbf{376} (1) (2020), 187-225.

\bibitem{ML63} 
S.~Mac~Lane, \emph{Homology}, Grundlehren Math. Wiss., \textbf{114}, Academic Press, New York, (1963).

\bibitem{OPS}
S.~Opper, P. G. Plamondon, S. Schroll, \emph{A geometric model for the derived category of gentle algebras}, arXiv: 1801.09659, (2025).

\bibitem{QZZ}
Y. Qiu, C. Zhang, Y. Zhou, \emph{Two geometric models for graded skew-gentle algebras}, arXiv: 2212.10369, (2022).

\bibitem{S15}
S. Schroll, \emph{Trivial extensions of gentle algebras and Brauer graph algebras},
J. Algebra, \textbf{444} (2015), 183-200.

\bibitem{Sol}
A. Solotar,
\emph{The Gerstenhaber bracket in Hochschild cohomology: methods and examples},
Representation theory and beyond. Workshop and 18th international conference on representations of algebras, ICRA 2018, Prague, Czech Republic, August 13--17, 2018 2020, 287-298,

\bibitem{VW}
Y. Volkov, S. Witherspoon, 
\emph{Graded {Lie} structure on cohomology of some exact monoidal categories},
Homology Homotopy Appl., \textbf{26} (2) (2024), 79-98.

\end{thebibliography}
\end{document}